\documentclass[letterpaper, 10pt,conference]{ieeeconf}

\IEEEoverridecommandlockouts 
\usepackage{amsmath,amssymb,amsfonts}
\usepackage{graphicx}
\usepackage{textcomp}
\usepackage{graphicx}      % include this line if your document contains figures
\usepackage{algorithm}
\usepackage{algorithmicx}
\usepackage{bm}            % bold greek letters
\usepackage{amsfonts}
\usepackage{mathtools}      
\usepackage{mathrsfs}
\usepackage{array,multicol,multirow}

\usepackage{url}
\usepackage{hyperref}  
\usepackage{lettrine}
\usepackage{comment}

\newcommand{\p}{\bm{\mu}}      % Scheduling signal process
\newcommand{\pdim}{n_{\mu}}         % Dimension of the scheduling signal
\newcommand{\nx}{n_x}     
\newcommand{\ny}{n_y}
\newcommand{\nn}{n_n}
\newcommand{\un}{n_u}
\newcommand{\udim}{n_\mathrm{u}}
\newcommand{\ydim}{n_\mathrm{y}}
\newcommand{\xdim}{n_\mathrm{x}}

\newcommand{\msig}{\p_{\sigma}}
\newcommand{\psig}{p_{\sigma}}
\newcommand{\Psig}{P_{\sigma}}

\newcommand{\sigSet}{\sigma \in \Sigma}
\newcommand{\word}{\sigma_{1}\sigma_{2}\cdots \sigma_{k}}
\newcommand{\wordSet}[1]{w \in \Sigma^{#1}}
\newcommand{\expect}[1]{E\left[#1 \right]}

\newcommand{\hankredu}{\mathcal{H}_{ \alpha, \beta }}

\newcommand{\uw}[1]{\p_{#1}}
\newcommand{\zwy}{\mathbf{z}^{\mathbf{y}}_{w} }
\newcommand{\zwyd}{\mathbf{z}^{\mathbf{y}^{d}}_{w} }
\newcommand{\Htplusu}{\mathcal{H}_{t,+}^{\ub}}
\newcommand{\zwys}
{\mathbf{z}^{\mathbf{y}^{s}}_{w} }

\newcommand{\Tsigys}{T^{\y^s,\y^s}_{\sigma,\sigma}}
\newcommand{\Tsigy}{T^{\y,\y}_{\sigma,\sigma}}
\newcommand{\Tsigyd}{T^{\y^d,\y^d}_{\sigma,\sigma}}

\newcommand{\Lwys}{\Lambda_{w}^{\y^s,\y^s}}
\newcommand{\Lwy}{\Lambda_{w}^{\y,\y}}
\newcommand{\Lwyd}{\Lambda_{w}^{\y^d,\y^d}}

\newcommand{\zwu}{\mathbf{z}^{\mathbf{u}}_{w} }
\newcommand{\zwr}{\mathbf{z}^{\mathbf{r}}_{w} }

\newcommand{\zvr}{\mathbf{z}^{\mathbf{r}}_{v} }

\newcommand{\zsigy}{\mathbf{z}^{\mathbf{y}}_{\sigma}}
\newcommand{\zsigys}{\mathbf{z}^{\y^s}_{\sigma}}
\newcommand{\zsigyd}{\mathbf{z}^{\y^d}_{\sigma}}

\newcommand{\timeset}{t \in \mathbb{Z}}
\newcommand{\yb}{\mathbf{y}}

\newcommand{\vb}{\mathbf{v}}
\newcommand{\eb}{\mathbf{e}}
\newcommand{\xb}{\mathbf{x}}
\newcommand{\ub}{\mathbf{u}}
\newcommand{\z}{\mathbf{z}}
\renewcommand{\r}{\mathbf{r}}
\renewcommand{\b}{\mathbf{b}}

\newcommand{\varu}{Q_{u}}

\newcommand{\Bs}{G}   % Stochastic B matrix

\newcommand{\covseq}{\Psi_{\mathbf{y}}}
\newcommand{\incov}{\Psi_{\mathbf{u},\mathbf{y}}}
\newcommand{\alphas}{\bar{\alpha}}
\newcommand{\betas}{\bar{\beta}}

\newtheorem{Notation}{Notation}
\newtheorem{Definition}{Definition}

\newtheorem{Theorem}{Theorem}
\newtheorem{Corollary}{Corollary}
\newtheorem{Lemma}{Lemma}
\newtheorem{Remark}{Remark}
\newtheorem{Assumption}{Assumption}
\newenvironment{pf}[1]{\begin{proof}{#1}}{\end{proof}}

\let\labelindent\relax
\usepackage{enumitem}
\renewcommand{\theenumi}{\arabic{enumi}}

\renewcommand{\theenumii}{\arabic{enumii}}

\renewcommand{\theenumiii}{\arabic{enumiii}}

\renewcommand{\theenumiv}{\arabic{enumiv}}

\newcommand{\btheta}{\boldsymbol{\theta}}

\newcommand{\GBS}{\textbf{GBS}}

\newcommand{\y}{\textbf{y}}
\newcommand{\bu}{\textbf{u}}
\newcommand{\x}{\textbf{x}}

\newcommand{\e}{\textbf{e}}

\title{Minimal covariance realization and system identification algorithm for a class of stochastic linear switched systems with i.i.d. switching}
\author{Elie Rouphael, Manas Mejari, Mihaly Petreczky,  Lotfi Belkoura%
\thanks{E. Rouphael, M. Petreczky and  L. Belkoura are with
Univ. Lille, CNRS, Centrale Lille, UMR 9189 CRIStAL, F-59000 Lille, France
{\tt\small first.lastname@univ-lille.fr}
}%
\thanks{M. Mejari is with Swiss AI lab IDSIA-SUPSI, Lugano, Switzerland.
{\tt\small manas.mejari@supsi.ch}%
}
}
\begin{document}
\maketitle
\thispagestyle{empty}
\pagestyle{empty}

\begin{abstract}
In this paper, we consider stochastic realization theory of Linear Switched Systems (LSS) with i.i.d. switching. We characterize minimality of stochastic LSSs and show existence and uniqueness (up to isomorphism) of minimal LSSs in innovation form. We present a realization algorithm to compute a minimal LSS in innovation form from output and input covariances. Finally, based on this
realization algorithm, by replacing true covariances with empirical ones, we propose a statistically consistent system identification algorithm.
%and 
%evaluate its performance.
%The later, based on the realization algorithm, uses the collected data in order to compute the covariances of the inputs and outputs. 
\end{abstract}

\section{Introduction}
In this paper we consider the stochastic realization problem for \emph{stochastic linear switched systems (abbreviated as LSS)}, i.e. for state-space representations:
 \begin{equation}
	\label{eq:aslpv}
            \mathcal{S} \left\{     \begin{aligned}
                 &  \xb(t+1)=A_{\btheta(t)}\xb(t)+B_{\btheta(t)}\ub(t)+K_{\btheta(t)}\vb(t) \\
                 & \y(t)=C\xb(t)+D\ub(t)+F\vb(t)
		\end{aligned}\right.
               \end{equation}
           where 
	$A_{\sigma} \!\! \in \!\! \mathbb{R}^{\nx \times \nx}, B_{\sigma} \! \in \! \mathbb{R}^{\nx \times \un}, K_{\sigma} \! \in\!  \mathbb{R}^{\nx \times \nn}$,
    $\sigma \!\in\! \Sigma=\{1,\ldots,\pdim\}$,
	$C\! \in \! \mathbb{R}^{\ny \times \nx}$, $F \! \in \! \mathbb{R}^{\ny \times \nn}$,
    $D \! \in \! \mathbb{R}^{\ny \times \un}$
	and $\xb$, $\vb$, $\bu$, $\y$, $\btheta$ 
    are the stochastic state, noise, input, output 
    and switching processes.
    %process process, $\vb$ is the noise process and $\y$ is the output process, 
    %$\p=\begin{bmatrix} \p_1,\ldots,\p_{\pdim} \end{bmatrix}^T$ is the 
    %switching process defined as follows: there exists a stochastic process 
    %$\btheta$ is the switching process taking values 
    %in $\Sigma$.
    The process $\btheta$ takes values in the set of discrete states $\Sigma$.
    %and we refer to the values of $\btheta$ as 
    %discrete states or modes. 
    %such that
    %and $\p \in \begin{Bmatrix} b_i \end{Bmatrix}_{i=1}^{\pdim}$, where $b_i = \begin{bmatrix} 0, \ldots, 0,1,0, \ldots, 0 \end{bmatrix}$ is the $\pdim$-tuple with all components equal to 0, except the \emph{i}th, which is 1. Then, the switching process can be represented as follows:
   %$\p_{i}(t)=\chi(\btheta(t)=i)$, where $\chi$ is the indicator µ%function, i.e. 
  %$\chi(\btheta(t)=i)=\begin{cases}
  %1 & \text{if $\btheta(t)=i,$}\\ 
  %0 & \text{otherwise.}
  %\end{cases}$ for all $i \in \{1, \ldots, \pdim\}$, $t \in \mathbb{Z}$. 
   %In this paper, for the sake of simplicity, we assume that $\btheta$ is an i.i.d. process.
   We call $\mathcal{S}$ from \eqref{eq:aslpv} a \emph{realization}
   of the tuple  $(\tilde{\y},\bu,\btheta)$, if $\tilde{\y}=\y$.
   We call $\nx$ the \emph{dimension} of $\mathcal{S}$.
   %The dimension $n$ of the state-space of $\mathcal{S}$
   %is called the dimension of $\mathcal{S}$.
   %$\mathcal{S}$ is said to be
   %µ\emph{minimal}, if it has the smallest possible dimension among all
   %LPV

   %Intuitively, an LSS can be viewed as a system which
   %switches between several
   %linear time-invariant (LTI) systems. 
   LSSs of the form \eqref{eq:aslpv} contain LTI systems as a special case,
   when $\pdim=1$. 
   Moreover, they correspond to \emph{jump-Markov linear systems} if $\btheta$ is a Markov process \cite{CostaBook}. Both switched and jump-Markov systems have a rich literature and a wide variety of applications, e.g., \cite{Sun:Book,CostaBook}.   
   We  refer to \eqref{eq:aslpv} as switched rather than jump-Markov systems, because 
   in contrast to the latter, in the formulation of the system identification problem, we treat the switching process as an external input, whose role is similar to that of $\bu$.
   %unlike the latter we include the case when the switching is a control input \cite{Sun:Book}. 
   %In contrast, for jump-Markov systems the switching is not a control input \cite{CostaBook}. We think of \eqref{eq:aslpv} either as a jump-Markov system, or a switched system, for which the switching is a control input and wich for purposes of system identification was chosen to be a sample of an i.i.d. (binary noise) process. In the literature, the definition of both switched and jump-Markov systems allows for the matrices of the output equation to depend on the discrete mode. 
   In contrast to the standard definition \cite{CostaBook,Sun:Book}
   for the sake of simplicity, we assume that the matrices of the output equation do not depend on the current discrete mode.

  \textbf{Motivation and context:}
    %In this paper we study stochastic realization theory of LSSs. In particular, we
    %investigate how to compute
    %minimal LSS in innovation form
    %both from output covariances and from data, and we
    %show  uniqueness of such systems up to isomoprhism. 
    %In particular we will address uniqueness of LSS in innovation form,
    %and its application to system identification. 
    In order to motivate the contribution of the paper, we firs define the system identification problem of LSSs. To this end, 
    let us define the \emph{deterministic behavior} $\mathcal{B}$ of 
    a tuple $S=(\{A_i,B_i,K_i\}_{i=1}^{\pdim},C,D,F)$ of matrices
    as the set of all deterministic signals $(y,u,q)$,
    all defined on $\mathbb{Z}$, taking values in 
    $\mathbb{R}^{\ny},\mathbb{R}^{\un},\Sigma$
    such that 
    there exists state a trajectory $x$ and noise realization $v$ satisfying:
    \begin{equation}
	\label{eq:aslpv:def}
            \begin{aligned}
                 &  x(t+1)=
               %\sum_{i=1}^{\pdim} 
               A_{q(t)}x(t)+B_{q(t)}u(t)+K_{q(t)}v(t) \\
                 & y(t)=Cx(t)+Du(t)+Fv(t)
		\end{aligned}
    \end{equation}
  That is, $\mathcal{S}$ can be viewed   as a stochastic version of \eqref{eq:aslpv:def}, %where the signals, $y,x,u,v,q$ are sample paths
  %µof the processes $\y,\x,\bu,\vb,\btheta$.
  in particular, all samples paths of $(\y,\bu,\btheta)$ are 
  elements of the deterministic behavior of $S$. 
  %The system identification problem is as follows: 
  
  Assume that $\mathcal{S}$ from \eqref{eq:aslpv}, referred to as the true system, is a realization of $(\y,\bu,\btheta)$, and we observe a finite portion $\{(y(t),u(t),q(t))\}_{t=0}^{N}$ of a sample path of $(\y,\bu,\btheta)$. The \emph{identification} task is to find matrices
 $\hat{S}=(\{\hat{A}_i,\hat{B}_i,\hat{K}_i\}_{i=1}^{\pdim},\hat{C},\hat{D},\hat{F})$ 
 %is the outcome of a system identification
 %algorithm applied to $N$ data points, 
 from $\{(y(t),u(t),q(t))\}_{t=0}^{N}$ 
  % of a sample path $(y,u,q)$ such that the following holds:  
  such that in the limit, as $N \rightarrow \infty$, %we would like the 
  the deterministic behavior of $\hat{S}$ equals (approximates) the deterministic behavior of the true systems. 
  The motivation for this problem formulation is that in many control problems, the switching signal is viewed either as a control input or as an arbitrary external input/disturbance.
  That is, the output of the estimated model should be close to the that of the true one for switching and inputs different from $\btheta$ and $\bu$. For jump-Markov systems, this can be relaxed by requiring that 
   the estimated model has the same output response as the true one for arbitrary inputs but for the fixed switching process $\btheta$.
  %the designated the stochastic input $\bu$ and 
  %switching $\btheta$.%Instead, this should hold for any input, switching and any choice of
  %the noise process. 
  In particular, we think of the stochasticity of 
  the observed signals as
  %are samples of stochastic processes as
  %that the input and switching were sampled from $\bu$ and $\btheta$, and that the true system was driven by a stochastic noise $\vb$ 
  a persistency of excitation assumption, not an assumption on the desired  model.
  %However, the estimated model
  %should approximate the true one for other choices of inputs and switching
  %and noise signals. 
  %while for identification experiment the input, noise and switching signals are sampled from 
  %$\bu,\bv,\btheta$, the output estimated
  % system should approximate the true one for any input, noise and switching. 
   %That is, the estimated LSS should not only recreate (approximately)
   %the behavior
   %of the true system for the switching and input which are sampled from
   %$(\y,\bu,\btheta)$, but for all inputs and switching signals.
   %problems assume that the switching and noise
   %signals are sampled from fixed stochastic processes. 
   %However, even in this case
   %the estimated model should approximate the true one for 
   %all inputs.
   
   %the estimated and true system should have (approximately)
   %the same output response for all inputs but only for switching signals sampled by $\btheta$
   %and noise signals sampled from some process
   %$\vb$. 

   One way to make the problem above well-posed is to ensure that
   both the estimated and true systems belong to a class of
   LSSs with the following property:
   if two elements of this class generate the same output for
   some stochastic input, switching and noise process, then their
   matrices differ only  by a change of basis, i.e., they are isomorphic. 
   Then if the output response of the estimated system is approximately the same as that of the true one for the designated processes $\bu$ and $\btheta$, then 
   the it is isomorphic to the true one, and hence the two systems have the same deterministic behavior. 
   %We then need identification algorithms which generate
   %LSSs of this class, and we should  bassume that
   %the true system belongs to this class without
   
    For LTI systems, minimal systems in innovation form \cite{LindquistBook}
    represent such a class, as any two minimal LTI systems in innovation form realizing the same output are isomorphic.
    Moreover, under suitable conditions, any stochastic LTI system can be transformed into a minimal one in innovation form while preserving its output \cite{LindquistBook}. Hence, it is reasonable to assume that the true system is of this class. In addition,are also several
    algorithms which return minimal systems in innovation form,
    e.g., subspace identification methods \cite{BookOverschee,KatayamaBook}
    and some parametric methods \cite{LjungBook}. 
       For the latter, the use of minimal 
    systems in innovation form also justifies viewing LTI systems as 
    optimal predictors of the  current output based on past outputs and inputs.
    The uniqueness (up to isomorphism)
    of minimal systems in innovation form is also used for 
    proving consistency of parametric \cite{HannanBook2012,HanzonBook}
    and subspace methods \cite{LindquistBook,CHIUSO2005377}.
    %and any two such LTI systems generating the same output are isomorphic. This observation is central for justifying classical system identification  setups \cite{LjungBook}. First, it allows us to 
    %view LTI systems as and try to estimate the underlying system matrices 
    %by minimizing the prediction error. Second, if we assume that the data generator was a minimal LTI system in innovation form, and we know that the estimated model 
    %µgenerates (approximately) the same output as the true one for the input process use in the identification experiment, then  they are (approximately) isomorphic. Therefore, the true and the estimated system generate (approximately) the
   %same output for any input and any noise realization. This allows proper definition and proofs  of consistency \cite{HannanBook2012,HanzonBook}
   %and also justifies the use of the estimated models for control. 
     %it means that estimating the data generator matrices is related to finding
   %which is central 
   %ustifies the use of prediction error minimization algorithms. 
   In fact, the covariance realization algorithm for LTI systems in innovation form \cite{LindquistBook} is the basis for subspace identification algorithms. %In fact, for subspace identification uniqueness of minimal LTI systems in innovation form was used implicitly in the problem formulation, as it guaranteed that finding the system matrices up to a change of basis is a well-posed problem. 
    %Finally, realization theory
   %Mhelped to analyze consistency of these algorithms %\cite{LindquistBook,CHIUSO2005377}.
    %First, they allow us to view the data generating system as
    %an optimal predictor, which uses past inputs and outputs to predict
    %future states. This then allows us to view the problem of estimating
    %the matrices of the data generator as the problem of finding the best
    %such predictor. This point of view is central to traditional system
    %identification \cite{LjungBook}. 

    %It is not a-priori clear how to define innovation form for LSSs so
    %that the favourable properties known for LTI system would carry over. 
    %In addition, 

  \textbf{Contribution:}
   We extend the concept of
   minimal (dimensional) systems in innovation form to a subclass 
  LSSs and  apply it
   to system identification.
   We assume that both $\btheta$
   and $\bu$ are i.i.d. processes. These assumptions represent the simplest case of persistently exciting data, and 
   there is little hope to tackle
   the more general case without solving this one first.
   %Moreover, in any case, the estimated models
   %are valid for all  inputs and switching signals. 
   
   We present a necessary and sufficient
   condition for minimality LSS\ realization of 
   $(\y,\bu,\btheta)$ in innovation form, and we show 
   that all such realizations of an output are isomorphic.  
   In addition, we present a realization
   algorithm for computing such LSSs from 
   covariances of inputs and outputs. 
   %Moreover, it is shown that minimal LSS realization of the same
   %tuple $(\y,\bu,\btheta)$ are related by a linear state-space transformation.
   
   Finally, we present a statistically consistent 
   system identification algorithm which is based
   on the latter realization algorithm, and which returns, in the limit, a minimal LSS in innovation form.
   %his identification algorithm
   %is statistically consistent, and in the limit it returns a
   %minimal LSS in innovation form which realizes $(\y,\bu,\p)$. 
   In particular, if the true system
   is a minimal LSS in innovation form, then in the limit
   this identification
   algorithm returns a LSS which is isomorphic to the true LSS,
   hence, it has the same deterministic behavior.
   We further improve
   %number of methods to improve
   this identification
   algorithm  %combining the algorithm 
   by combining it with Gradient-Based (GB) methods \cite{cox2018towards}.

   \textbf{Related work:}
     Realization theory of deterministic
    switched systems was studied in \cite{MPLBJH:Real}, 
    for stochastic switched  systems with no inputs in \cite{CDC2007,PetreczkyBilinear,RouphaelPetreczkyArxiv}. 
    %In contrast to \cite{CDC2007,PetreczkyBilinear,RouphaelPetreczkyArxiv}, the current paper covers the case of 
    %systems with input, but it has more restrictive assumptions for the switching signal and output equations. 
    To the best of our knowledge, our results on 
    realization theory of stochastic LSS with inputs are completely new. 

    Identification of switched and jump-Markov systems is an active research area, see the \cite{LauerBook,SurveyPaoletti,JumpMarkov1,JumpMarkov4,JumpMarkov5,bako2009line}, and the references therein. However, most of the literature  assumes
    that the switching signal is unobserved, which makes the problem more challenging but also leads to lack of consistency results for state-space representations. A consistent subspace
    identification algorithms was presented in \cite{PETRECZKY2023101308} for noiseless LSSs. In \cite{sarkar2019nonparametric} a consistent identification algorithm for noisy LSSs was presented,
    but the noise gain matrix and the noise covariance 
    were not estimated, and the identification was based on several i.i.d. time series.

    %For jump-Markov systems there is a long history of applying
    %expectation-maximization based algorithms (EM). However, in general it
    %is unclear to which extent these methods are consistent, and more,
    %importantly, if the input-output behavior of the estimated system
    %is close to the true one for inputs and switching signals which are
    %different from the ones used for identification.
    %Moreover, most of the methods assume a certain parameterization
    %of the state-space representation, with the exception of \cite{bako2009line}. However, \cite{bako2009line} assumes only output noise and does not propose consistency guarantees.

    LSSs can be viewed as a subclass of linear-parameter varying (LPV)
    systems, if the switching signal viewed as a discrete scheduling. There is a wealth
    of literature on system identification for LPV systems, including subspace methods, \cite{Toth2010SpringerBook,CoxTothSubspace,MejariLPVS2019,Verhaegen:CDC04,Verdult02} and the references therein.
     %However, most of the literature LPV susbspace identification provides no proofs of consistency of subspace  methods, and most of them do not estimate the noise gains and noise variance. 
     Stochastic
     realization theory of LPV systems was invetsigated in \cite{MejariLPVS2019,CoxLPVSS}.
     In \cite{CoxLPVSS} the existence of LPV systems in innovation form
     was investigated, 
     %but the resulting state-space representation 
     %had a dynamical dependence on scheduling variables, i.e., 
     but the noise
     gain matrices of the obtained system depended on the current and past scheduling. 
     Moreover, \cite{CoxLPVSS} did not address minimality and uniqueness 
     of systems in innovation form. 
     Unlike the identification algorithm of this  paper, the existing subspace
     methods, see the overview \cite{CoxTothSubspace}, 
     are not proven to be  consistent or to return a minimal system in innovation form.
     %for consistency or for returning a minimal system innovation form. 
     
    A statistically consistent identification algorithm was
    presented in \cite{MejariCDC2019} for noisy LPV systems with no inputs, and in \cite{MejariLPVS2019}
    for noisy LPV systems with inputs. This paper is an extension of \cite{MejariLPVS2019},
    the main difference w.r.t. \cite{MejariLPVS2019} is as follows:
    \textbf{(1)} the conditions on the scheduling signals from  \cite{MejariLPVS2019} do not directly apply for i.i.d. switching signals (see Remark \ref{rem:manas1});
    \textbf{(2)} we characterize minimality and uniqueness (up to isomorphism)
    of systems in innovation form; \textbf{(3)} the proposed identification algorithm
    returns a minimal system in innovation form, while the one of \cite{MejariLPVS2019}
    returns a non-minimal one;  that is, the system returned by the algorithm from 
    \cite{MejariLPVS2019} need not have the 
    deterministic behavior as the true system;  \textbf{(4)} we present novel modifications to improve the empirical behavior of the identification algorithm.
    %which were absent in \cite{MejariLPVS2019}. 
    %The results of the paper
    %can be extended to include the LPV systems covered by \cite{MejariLPVS2019}.

    To sum up, the identification algorithm of this paper 
    is the first algorithm for noisy LSSs that has all of  the following properties: \textbf{(a)} it estimates the noise gain and the noise covariance, \textbf{(b)} returns a minimal LSSs in innovation form,  and \textbf{(c)} applies to one single long time series. 
    %In particular,
    %\cite{MejariCDC2019}
    %does not apply LSSs with inputs considered in this paper.
    %, i.e., to LPV systems with binary white noise scheduling and control inputs. 
    %MThe results of 
    %re the closest one to the one of this paper. However, \cite{MejariLPVS2019} does not apply to
    %switched systems.
    %as it assumes zero mean i.i.d. scheduling with a diagonal covariance matrix, and  i.i.d switching
    %does not satisfy these assumptions.
    %and there seem to be no obvious way to adapt it to fit the assumptions of \cite{MejariLPVS2019}. 
    %Moreover, minimality of innovation forms  was not investigated in  %\cite{MejariLPVS2019}, and the identification algorithm of \cite{MejariLPVS2019}
    %does not return a minimal system. % in innovation form. 
    %In turn, it means that even if it was possible to extend the algorithm \cite{MejariLPVS2019} to the switched case easily, 
    %it is unclear if 
    %However, similarly to \cite{MejariLPVS2019} the current paper
    %relies heavily on \cite{PetreczkyBilinear}, and 
    %the results of this paper could probably be extended to 
    %MLPV systems considerd in \cite{MejariLPVS2019}. 
    
    This paper uses technical results from \cite{LCSSArxive}
    on the decomposition of outputs of LSSs into stochastic and deterministic components.
    %which 
    %The system class considered in \cite{LCSSArxive}
    %considered both LPV  systems and LSSs. 
    %which treats the decomposition of outputs of 
    %generalized LSS systems into a stochastic and deterministic components, and it uses this decompositon 
    The latter is used in \cite{LCSSArxive} 
    to show existence of a realization in innovation form and to provide \emph{sufficient}
    conditions for minimality and uniqueness of LSSs in innovation form.
    In contrast, we formulate 
    \emph{necessary and sufficient} conditions for minimality, and we show isomorpĥism of a much wider class of LSSs in innovation forms, see Remark \ref{lcss:rem} for a detailed
    explanation of the gap between the two results. 
      % both the LPV systems of \cite{MejariLPVS2019} and LSSs.
    %The main difference w.r.t.  \cite{LCSSArxive} is that 
    %\textbf{(1)} instead of sufficient we propose %conditions which are 
    %\emph{necessary and sufficient} conditions  for minimality, 
    %\textbf{(2)} 
    %of all minimal 
    %LSSs in innovation form, not only those which satisfy the sufficient condition of \cite{LCSSArxive},  
    Moreover, in contrast to this paper, \cite{LCSSArxive} proposes no realization or identification algorithm.
    %An  version of this paper with detailed proofs is availabe in the report \cite{rouphael2024minimal}.
    %is an extended version of this paper and it contains detailed proofs.
    
    %we propose a covariance realization algorithm
    %for computing a minimal LSS in innovation form, \textbf{(4)} we formulate a consistent identification algorithm.
    %In particular, our results demonstrate that the sufficient conditions \cite{LCSSArxive} are far from being necessary. 
    %Furthermore, the results of this paper show that
    %the sufficient conditions of \cite{LCSSArxive} are far from being necessary. 

%\input{Technical_assumptions_on_the_scheduling_process}
\section{Technical preliminaries}
\label{sect:prelim}
%We will present some technical notations and definitions regarding the involved processes, in order to help prepare the exposition of the main results.
In this section we recall some concepts from \cite{PetreczkyBilinear} which we use
to define a suitable subclass of LSSs. 
\\
\textbf{Probability theory:}
We use the usual notation and terminology of probability theory \cite{Bilingsley}. All the random variables and stochastic processes are understood w.r.t. to a probability space $\left(\Omega, \mathcal{F}, \mathcal{P}\right)$. %where $\mathcal{F}$ is a $\sigma$-algebra over $\Omega$.
	%(\emph{i.e.}, $\mathcal{F}$ is a collection of subsets of $\Omega$, that includes $\Omega$ itself, is closed under complement, is closed under countable unions and is closed under countable intersections) and $\mathcal{P}$ is a probability measure on $\mathcal{F}$. For two $\sigma$-algebras $\mathcal{F}_i$, $i=1,2$, $\mathcal{F}_1 \lor \mathcal{F}_2$ denotes the smallest $\sigma$-algebra generated by the $\sigma$-algebras $\mathcal{F}_1,\mathcal{F}_2$. 
We use $E$ for expected values.
%and denoted by $E[\mathbf{h}]$ and conditional expectation w.r.t. $\sigma$- algebra $\mathcal{F}$ is denoted by $\expect{\mathbf{h} \mid \mathcal{F}}$. 
All the stochastic processes in this paper are discrete-time ones defined over the time-axis $\mathbb{Z}$ of the set of integers.
%a stochastic process $\r$ is a collection of random variables
%µ$\{\r(t)\}_{t \in \mathbb{Z}}$ taking values in some set $X$. 
We use bold letters for random variables and stochastic processes.
\\
\textbf{Switching process:}
 The  process
 $\btheta$ is independent, identically 
 distributed taking values in $\Sigma=\{1,\ldots,\pdim\}$,
 and 
 %i.e., $\btheta(t)$ and $\btheta(s)$ are independent for all $s, t \in \mathbb{Z}$ and
 %the distribution of  $\btheta(t)$ does not depend on $t \in \mathbb{Z}$. Moreover, in the rest of
% the paper, we use the following notation 
 \begin{equation}
 \label{psigmas}
 p_{\sigma} = \mathcal{P}(\btheta(t)=\sigma), \quad  \sigma \in \Sigma. \end{equation}
\textbf{Sequences of discrete modes:}
A \emph{non empty word} over $\Sigma$ is a finite sequence of letters, i.e., $w = \sigma_{1}\sigma_{2}\cdots \sigma_{k}$, where $0 < k \in \mathbb{Z}$, $\sigma_{1}, \sigma_{2}, \ldots, \sigma_{k} \in \Sigma$. The set of \emph{all} nonempty words is denoted by $\Sigma^{+}$. We denote an \emph{empty word} by $\epsilon$. Let $\Sigma^{*} = \epsilon \cup \Sigma^{+}$. 
We define the concatenation of words in the standard manner.
%The concatenation of two nonempty words $v = a_{1}a_{2}\cdots a_{m}$ and  $w= b_{1}b_{2}\cdots b_{n}$ is defined as $vw = a_{1}\cdots a_{m} b_{1} \cdots b_{n}$ for some $m,n > 0$. Note that if $w = \epsilon$ or $v= \epsilon$,  then $v\epsilon = v$ and $\epsilon w = w$, moreover, $\epsilon \epsilon = \epsilon$. 
The length of the word $w \in \Sigma^{*}$ is denoted by $|w|$, and $|\epsilon| =0$. 
\\
\textbf{Matrices and matrix products:}
%\begin{Notation}[Matrix Product]\label{not:product}
 Denote by  $I_n$  $n \times n$
    identity matrix.
	Consider $n \times n$ square matrices $\{A_{\sigma}\}_{\sigSet}$. For  the empty word $\epsilon$, let $A_{\epsilon} = I_n$.
 for any word %$\wordSet{+}$ of the form
 $w = \sigma_{1}\sigma_{2}\cdots\sigma_{k} \in \Sigma^{+}$, $k\!>\!0$ and $\sigma_{1}, \ldots, \sigma_{k} \in \Sigma$, we define 
	\begin{equation}
    \label{not:product}
	A_{w} = A_{\sigma_{k}}A_{\sigma_{k-1}}\cdots A_{\sigma_{1}}.
	\end{equation}
	%	For example, $w = 123$, $A_{w} = A_{123} = A_{3}A_{2}A_{1}$
%\end{Notation}
\textbf{ZMWSII processes:}
 Next, we recall the concept of zero-mean wide-sense stationary (ZMWSII) process w.r.t. $\btheta$
 from \cite{PetreczkyBilinear}. To this end,
 we introduce the following indicator process
 $\p(t)=\begin{bmatrix} \p_1(t), & \ldots, & \p_{\pdim}(t) \end{bmatrix}^T$ where 
 \begin{equation}
 \label{scheduling:def}
    \p_{\sigma}(t)=\chi(\btheta(t)=\sigma)=\left\{ 
      \begin{array}{ll}
       1 & \mbox{ if } \btheta(t)=\sigma \\
      0 & \mbox{ otherwise}
    \end{array}\right.
 \end{equation}
 where $\chi$ is the indicator function. 
 The process $\p$ is admissible
 in the terminology of \cite{PetreczkyBilinear}.
 \begin{Remark}[Relationship with \cite{MejariLPVS2019}]
 \label{rem:manas1}
 If we view $\p$ as a scheduling signal, then it does not satisfy the assumptions of \cite{MejariLPVS2019}: while 
 it is an i.i.d. process , it is not zero mean, $\p_i(t),\p_j(t)$
 are not uncorrelated for $i,j \in \Sigma$, and $\p_1(t) \ne 1$. By applying a suitable normalizing affine transformation
 to $\p$ we can bring it to a form
 %$\p$ which normalizes it  which involves substracting the mean of $\p$
 %and mutlitplying it by the square root of its covariance,
 %then the resulting process 
 which satisfies \cite{MejariLPVS2019}; for example, for 
 $\pdim=2$, we can take $\tilde{\p}(t)=(\p_1(t)+\p_2(t),  \p_2(t)-p_2)^{T}=(1,\p_2(t)-p_2)^T$.
 The system matrices of the thus arising LPV system are affine combinations of the system matrices of original LSS. However, this transformation does not appear to 
 simplify the derivation of the results. 
 To the contrary, it would require handling the affine transformation between system matrices when analyzing minimality.
 \end{Remark}
 
 Next, we need the following
 products of components of $\p$.
 For every word $\wordSet{+}$ where $w=\word$, $k \geq 1$, $\sigma_{1},\ldots, \sigma_{k} \in \Sigma$, we define the  process $\uw{w}$ as follows:
\begin{equation}
%	\begin{align}
\label{eqn:uw}
		\uw{w}(t) = \uw{\sigma_{1}}(t-k+1)\uw{\sigma_{2}}(t-k+2)\cdots\uw{\sigma_{k}}(t)  %\forall t \in \mathbb{Z} \nonumber \\
	%	p_w&=p_{\sigma_1}p_{\sigma_2} \cdots p_{\sigma_k}.
%	\end{align}
\end{equation}
For an empty word $w= \epsilon$, we set $\uw{\epsilon}(t)=1$.
For the constants $\{p_{\sigma}\}_{\sigma \in \Sigma}$ from \eqref{psigmas}
define the products
\begin{equation}
	\begin{aligned}\label{eqn:uw1}
%		\uw{w}(t) &= \textcolor{red}{\uw{\sigma_{1}}(t-k+1)\uw{\sigma_{2}}(t-k+2)\cdots\uw{\sigma_{k}}(t) } \\ %\forall t \in \mathbb{Z} \nonumber \\
		p_w&=p_{\sigma_1}p_{\sigma_2} \cdots p_{\sigma_k}
	\end{aligned}
\end{equation}
for any $w \in \Sigma^{+}$, $\sigma_i \in \Sigma$. For an empty word $w= \epsilon$, we set $p_{\epsilon}=1$.
For a stochastic process $\mathbf{r} \in \mathbb{R}^{r}$ and for each $\wordSet{*}$ we define the stochastic process $\zwr$ as
\begin{equation}\label{eqn:zwu}
	\zwr(t) = \mathbf{r}(t-|w|) \uw{w}(t-1)\frac{1}{\sqrt{p_{w}}},
\end{equation}
where $\uw{w}$ and $p_w$ are as in \eqref{eqn:uw} and \eqref{eqn:uw1}. For $w=\epsilon$, $\zwr(t)=\r(t)$.
The process $\zwr$ in \eqref{eqn:zwu} is interpreted as the product of the \emph{past} of $\r$ and $\p$.
The process $\zwr$ will be used as predictors for future values of $\r$ for various choices of $\r$.
\\
We are now ready to recall from \cite{PetreczkyBilinear} the definition of ZMWSII process w.r.t. $\btheta$
(w.r.t; $\p$ in the terminology of \cite{PetreczkyBilinear}).
	% \subsection{Properties of the input and output processes}
%\label{Subsec:assumptions1}
 %Below we define the concept of 
 %admissible input processes and 
 %ZMWSII processes w.r.t. switching process. These concepts will be used to define the nature of the processes and the subclass of the LSSs we will
 \begin{Definition}[\cite{PetreczkyBilinear}]\label{def:ZMWSSI}
	A stochastic process $\r$ is \emph{Zero Mean Wide Sense Stationary} (ZMWSSI) if 
	\begin{enumerate}
		\item For $\timeset$, the $\sigma$-algebras generated by $\{\mathbf{r}(k) \}_{k \leq t}$, and $\{\msig(k) \}_{k \geq t, \sigSet}$ respectively are conditionally independent w.r.t $\sigma$-algebra generated by 
        $\{\msig(k) \}_{k < t, \sigSet}$, %and $ denoted by $\filtr$, $\filt$ and $\filtp$ respectively, are such that $\filtr$ and $\filtp$ are conditionally independent w.r.t. $\filt$.
		%i.e., intuitively, future values of $\p$ do not directly depend on the past values of $\mathbf{r}$, however only via dependence of past $\mathbf{r}$ on past values of $\p$.
		\item The processes $\{\zwr \}_{\wordSet{*}}$  are zero mean, square integrable and jointly wide sense stationary, i.e. the covariances 
        %$\expect{\r(t)(\zwr(t))^{T}}$,
        $\{\expect{\zwr(t)(\zvr(t))^{T}}\}_{w,v \in \Sigma^{*}}$ are 
        independent of $t$.
         %\expect{\r(t) (\zwr(s))^{T}}
		%$\forall t,s,k \in \mathbb{Z},w,v \in \Sigma^{+}$, 
		%\begin{align*}
		%&   \expect{\mathbf{r}(t)} = 0, ~  %\expect{\zwr(t)} =0  \\
		%&	\expect{\r(t+k)(\zwr(s+k))^{T}} = \expect{\r(t) (\zwr(s))^{T}}, \\
		%&	\expect{\r(t+k)(\r(s+k))^{T}} = \expect{\r(t)(\r(s))^{T}}, \\
		%µ&	\expect{\zwr(t+k)(\zvr(s+k))^{T}} = %\expect{\zwr(t) (\zvr(s))^{T}}. 
		%		\end{align*}
	\end{enumerate}
\end{Definition} 
The concept of ZMWSII 
is an extension of  wide-sense stationarity, if 
$\Sigma^{+}$ is viewed as a time axis.
%Indeed, 
%Assume that  $\r$ is ZMWSII w.r.t. $\p$, then $\zwr(t)$ is wide-stationary and 
%square-integrable for all $w \in \Sigma^{+}$. Moreover,
%the covariances $\expect{\zwr(t)(\zvr(t))^{T}}$ do not depend on $t$
%and depend only on the common
%prefix of $w$ and $v$. 
%More precisely,
%these covariances are zero, unless 
%$w=sv$ or $v=sw$ for some $s \in \Sigma^{*}$. In the latter case
%$\expect{\zwr(t)(\zvr(t))^{T}}$  equals
%$\expect{\r(t)(\zsr(t))}$
%(resp. $\expect{\r(t)(\zsr(t))^{T}}$).
More precisely, 
%consider the following covariances defined 
for any two  ZMWSII 
processes $\r$ and $\mathbf{b}$
%let us define the following covariances
and for sequences $w,v\in \Sigma^{*}$,
define the covariances:
\begin{equation}
\label{cov:not}
%\Lambda_{w}^{\y}=E[\yb(t)(\z^{\yb}_{w}(t))^T] \\
\Lambda_{w}^{\r,\mathbf{b}}=E[\r(t)(\z^{\mathbf{b}}_{w}(t))^T], \quad 
T_{w,v}^{\r,\b}=E[\z^{\r}_{w}(t)(\z^{\b}_{v}(t))^T] 
\end{equation}
Then, $T_{w,v}^{\r,\r}=\Lambda^{\r,\r}_{s}$ if
$w=sv$, and $T_{w,v}^{\r,\r}=(\Lambda^{\r,\r}_{s})^T$
if $v=sw$, for some $s \! \in \Sigma^{*}$, and $T_{w,v}^{\r,
\r}\!\!=\!0$ otherwise.
I.e., the covariance $T_{w,v}^{\r,\r}$ depends
on the difference between $w$ and 
$v$.

\section{Stationary LSSs}

Below we present the definition of the class of LSSs which is studied in this paper. To this end, we recall from
\cite{LCSSArxive} the notion of a white noise process
w.r.t. $\btheta$.
\begin{Definition}[White noise w.r.t. $\btheta$, \cite{LCSSArxive}]
	A ZMWSII process $\mathbf{r}$ is a white noise w.r.t. $\btheta$, if for all $v,w \in \Sigma^{+}$, %$v \in \Sigma^{*}$, 
     %$\sigma \in \Sigma$,
     $T^{\r,\r}_{w,v}=0$ if $w \ne v$,
     and $T^{\r,\r}_{w,w}$ is non-singular.
	%and 
	%\[ E[\mathbf{r}(t)(\zwr(t))^T]=0, ~ 
    %     E[\z_{\sigma v}^{\r}(t)(\z_{\sigma w}^{\r}(t))^T]=E[\z^{\r}_{\sigma}(t)(\z^{\r}_{\sigma}(t))^T] > 0$, for all $w \in \Sigma^{+}$, $\sigma \in \Sigma$. 
	%\[ E[\zwr(t)(\zsvr(t))^T] \!\!=\!\! \left\%{\begin{array}{ll}
        %0 & \mbox {if } w \ne \sigma v \\
        %E[\z^{\r}_{\sigma}(t)(\z^{\r}_{\sigma}(t))^T]  & \mbox{if } w=\sigma v %\mbox{ and } \\ 
    %& \sigma \mbox{ first letter of } w         
    %\end{array}\right., 
    % \]
      %and $E[\z^{\r}_{\sigma}(t)(\z^{\r}_{\sigma}(t))^T]$ is
      %nonsingular for all $\sigma \in \Sigma$.
     %$E[\z_{\sigma w}^{\r}(t)(\z_{\sigma w}^{\r}(t))^T]=E[\z^{\r}_{\sigma}(t)(\z^{\r}_{\sigma}(t))^T] > 0$, for all $w,v \in \Sigma^{+}$. 
\end{Definition}
The notion of a white noise process w.r.t. $\btheta$ (w.r.t. $\p$
in the terminology of \cite{LCSSArxive}) is an extension of the usual concept of the white noise process if $\Sigma^{+}$ is
viewed as an additional time axis. 
In particular, if $\r$ is a white noise w.r.t. $\btheta$, the collection $\{\z_{w}^{\r}(t)\}_{w \in \Sigma^{+}}$ is a sequence of uncorrelated random variables and the covariance of $\z_w^{\r}(t)$ does not depend on $t$ and it depends only on the first letter of $w$.  

Next, we state the assumptions on the input and output.
%which is a special case of that of \cite{LCSSArxive}.
\begin{Assumption}[Inputs and outputs]
	\label{asm:main}\textbf{(1)}
	$\mathbf{u}$ is a white noise  w.r.t. $\btheta$
    and the covariance $T^{\ub,\ub}_{\sigma,\sigma}=Q_{\bu}$ 
    %$E[\z^{\ub}_{\sigma}(t)(\z^{\ub}_{\sigma}(t))^{T}]= E[\ub(t-1)(\ub(t-1))^{T}] =  \varu > 0$ 
    does not depend on
	$\sigma \in \Sigma$. 
	\textbf{(2)} The process
	$\begin{bmatrix} \mathbf{y}^T, \!\! & \!\! \mathbf{u}^T \end{bmatrix}^T$ is a ZMWSSI.
\end{Assumption}
Now, we are ready to define the class of stationary LSS, which is a special case of stationary generalized switched linear systems defined in \cite{LCSSArxive}.
\begin{Definition}[Stationary LSS]
\label{defn:LPV_SSA_wo_u}
%\label{def:Stationary}
 %The  GLSS 
  A  \emph{stationary LSS} (abbreviated 
  sLSS) %of $(\y,\bu,\theta)$ 
  is a system \eqref{eq:aslpv}, such that
%a process $\r$ taking values in $\mathbb{R}^{p}$, 
%is a system of the form \begin{equation}
%	\label{eq:aslpv}
%            \mathcal{S} \left\{     \begin{aligned}
%                 &  \tilde{\xb}(t+1)=\sum_{i=1}^{\pdim} (\tilde{A}_i\tilde{\xb}(t)+\tilde{K}_i\tilde{\vb}(t))\p_i(t) \\
%                 & \r(t)=\tilde{C}\tilde{\xb}(t)+\tilde{F}\tilde{\vb}(t)
%		\end{aligned}\right.
%               \end{equation}
%           where 
%	$\tilde{A}_{\sigma} \in \mathbb{R}^{\tilde{n} \times \tilde{n}}, \tilde{K}_{\sigma} \in \mathbb{R}^{\tilde{n} \times \tilde{m}}$,
%	$\tilde{C} \in \mathbb{R}^{p \times \tilde{n}}$, $\hat{F} \in \mathbb{R}^{p \times \tilde{m}}$, 
%        and $\vb$ is a process taking values in $\mathbb{R}^{\tilde{m}}$, 
%         and the following holds:
\begin{itemize}

	\item[\textbf{1.}]
		$\mathbf{w}=\begin{bmatrix} \vb^T, \!\! & \!\! \bu^T \end{bmatrix}^T$ is a white noise process w.r.t. $\btheta$. 

    \item[\textbf{2.}]
  The process
  $\begin{bmatrix} \xb^T\!\!, & \!\! \mathbf{w}^T \end{bmatrix}^T$  is a ZMWSSI, and $T^{\xb,\mathbf{w}}_{\sigma,\sigma}=0$,
  $\Lambda^{\xb,\mathbf{w}}_{w}=0$
     %$    
        %\forall (\sigma,w )\in \Sigma \times \Sigma^{+}:
      %  E[\z^{\xb}_{\sigma}(t)(\z^{\mathbf{w}}_{\sigma}(t))^T]=0, ~
      %  E[\xb(t)(\z^{\mathbf{w}}_w(t))^T]=0,
    %$
  for all $\sigSet$, $w \in \Sigma^{+}$.

  % and $E[\tilde{\vb}(t)\tilde{\vb}^T(t)\p_i^2]=Q_i$.
	
	\item[\textbf{3.}]
		The eigenvalues of the matrix $\sum_{\sigSet} \psig {A}_{\sigma} \otimes {A}_{\sigma}$ are inside the open unit circle.

    %\item[\textbf{4.}]  For all $\sigma_1,\sigma_2 \in \Sigma$, if $(\sigma_1,\sigma_2) \notin \mathcal{E}$, then
  %$A_{\sigma_2}A_{\sigma_1}=0$ and $A_{\sigma_2}\begin{bmatrix} B_{\sigma_1}^T & K_{\sigma_1}^T \end{bmatrix}^T  E[\z^{\mathbf{w}}_{\sigma_1}(t) (\z_{\sigma_1}^{\mathbf{w}}(t))^T]=0$.
		%	\end{enumerate} 
	%We call $\tilde{\xb}$ the state process and $\tilde{\vb}$ the noise process. 
       %We say that a asLPV-SSA is a realization of a process $\tilde{\r}$, if its output process equals $\r$. 
       \end{itemize}
If $B_{\sigma}=0$, $\sigma \in \Sigma$, and $D=0$  the we call \eqref{eq:aslpv} an \emph{autonomous stationary} LSS (asLSS), and in this case we say it is an asLSS realization 
of $(\y,\btheta)$.
%       (abbreviated as \emph{asGLSS of $(\y,\p)$}).
\end{Definition}

In the terminology of \cite{PetreczkyBilinear}, an sLSS corresponds to a stationary \GBS\  w.r.t. inputs $\{\p_{\sigma}\}_{\sigma \in \Sigma}$ and with noise 
$w=[\vb^T,\bu^T]^T$. In the terminology of \cite{LCSSArxive}, a sLSS is a stationary generalized switched system for which the switching process
is i.i.d.
%From \cite{PetreczkyBilinear}, if a process $\r$ can be represented by an asLSS, then $\r$ is a ZMWSSI process and 
%$\xb$ is uniquely determined by $\vb$ and the matrices $( %\{A_{\sigma},K_{\sigma}\}_{\sigma \in \Sigma}, C,D)$. 
Note that the processes $\xb$ and $\yb$ are ZMWSII, in particular, they
are wide-sense stationary, and that $\xb$ is orthogonal to the future values of the noise process $\vb$. 
%We concentrate on wide-sense stationary processes,
%as they reflect the stationary behavior of the underlying
%system, and they are easier to analyze. 
We need stationarity, as even for LTI case stochastic realization theory exists only for the stationary case \cite{LindquistBook}. 
%it is difficult to estimate the distribution of non-stationary processes. Also, wide-sense stationary processes solve the problem of the initial state conditions.

As it was noted \cite{LCSSArxive}, the state of an sLSS
is uniquely determined by its matrices and noise and input  process, i.e.,
%In order to
%present this relationship, we need the following notation.
%In order to define this notion more precisely, 
%from \cite[Lemma 2]{PetreczkyBilinear} it follows that
% From \cite{PetreczkyBilinear}, it follows that
\begin{equation}
	\label{stat:state:eq1}
	\xb(t)=\sum_{\sigma \in \Sigma, w \in \Sigma^{*}} \sqrt{p_{\sigma w}}  A_w\left(K_{\sigma}\z^{\vb}_{\sigma w}(t)+B_{\sigma} \z^{\bu}_{\sigma w}(t) \right)
\end{equation}
where the right-hand side
%on the right-hand side of \eqref{stat:state:eq1}
%is  absolutely convergent, and convergence 
is absolutely convergent in the mean-square sense. This motivates the following notation.  %We  use the following notation.
\begin{Notation}
\label{Notation:sLPSS}
 We identify a sLSS of $(\yb,\ub,\btheta)$  of the form \eqref{eq:aslpv} 
with the tuple  $\mathcal{S}=(\{A_{\sigma},B_{\sigma},K_{\sigma}\}_{\sigma=1}^{\pdim},C,D,F,\vb)$.
\end{Notation}

\section{Minimality of sLSS in innovation form}
\label{sect:real}
In this section, we present the main results on existence and uniqueness of minimal sLSSs in innovation form.  To this end, first, we recall realization theory of deterministic linear switched systems, and then
relate it with realization 
theory of sLSSs.

\subsection{Review of deterministic realization theory of LSSs}
\label{Sect:Chap3}
%Recall from \ref{Main:results} the deterministic LSS representation.
%In this section we recall from \cite{MPLBJH:Real,PetreczkyRealChapter}
%some elements of realization theory for deterministic linear switched. 
A \emph{deterministic linear switched system (abbreviated as dLSS)}, is a system of the form
\begin{equation}
%\label{petreczky:eq:LSSformA}
\label{eqn:LPV_SSA:det}
\mathscr{S} % \left \{
\begin{cases}
&\mathrm{x}(t+1)=\mathcal{A}_{q(t)}\mathrm{x}(t)+\mathcal{B}_{q(t)}\mathrm{u}(t) \\
& \mathrm{y}(t)=\mathcal{C}\mathrm{x}(t)+\mathcal{D}\mathrm{u}(t)
\end{cases}%\right.
\end{equation}
where $\{\mathcal{A}_{\sigma}, \mathcal{B}_{\sigma}\}_{\sigma \in \Sigma}, \mathcal{C}, \mathcal{D}$ are matrices of suitable dimensions,  
$q:\mathbb{Z} \rightarrow \Sigma$ is the switching signal, 
$\mathrm{x}:\mathbb{Z} \rightarrow \mathbb{R}^{\xdim}$ is the state trajectory
$\mathrm{u}:\mathbb{Z} \rightarrow \mathbb{R}^{\udim}$ is the input trajectory
$\mathrm{y}:\mathbb{Z} \rightarrow \mathbb{R}^{\ydim}$ is the output trajectory with finite support \footnote{a function $g:\mathbb{Z} \longrightarrow \mathbb{R}^p$ has a finite support, then it $\exists \; t_0 \in \mathbb{Z}$, such that $\forall \; t<t_0$, $g(t)=0$.}.
We identify a dLSS 
of the form \eqref{eqn:LPV_SSA:det} with the tuple 
\begin{equation} 
\label{eqn:LPV_SSA:det_not} 
\mathscr{S}=(\{\mathcal{A}_{i},\mathcal{B}_{i}\}_{i=1}^{\pdim},\mathcal{C},\mathcal{D})
\end{equation}
%The number $\nx$ is called the dimension of $\mathscr{S}$. 
%\textcolor{red}{Manas: please make the notation consistent with tha of stochastic LPV-SSA, i.e. use $(\{A_{i},B_{i}\}_{i=0}^{\pdim},C,D)$ instead.} 
We refer to \cite{MPLBJH:Real,PetreczkyRealChapter}
for a detailed definition input-output behavior of dLSSs, their realization theory, minimality, etc. 
%For our purposes, it is enough to view dLSSs as
%representations of Markov-parameters.
In particular, let us call any function $M:\Sigma^{*} \rightarrow \mathbb{R}^{n_{\mathrm y} \times n_{\mathrm u}}$ a \emph{Markov function}.
%is the function $M_{\mathscr{S}}:\Sigma^{*} \rightarrow \mathbb{R}^{n_\mathrm{y} \times \un}$, such
%that for all $w \in \Sigma^{*}$, 
The dLSS $\mathscr{S}$ \emph{realizes} $M$, if 
for all $w \in \Sigma^{*}$,
\begin{equation}\label{eqn:sub_markov}
	%\forall w \in \Sigma^{*}: ~ 
 M(w)=\left\{\begin{array}{ll}
		\mathcal{CA}_s\mathcal{B}_{\sigma}, \ \ & w=\sigma s, \, \sigma \in \Sigma, \, s \in \Sigma^{*} \\
		\mathcal{D}.           \ \  & w=\epsilon
	\end{array}\right.
\end{equation}
If $\mathscr{S}$ is a realization of $M$, then 
$M$ is referred to as the \emph{Markov function of $\mathscr{S}$} and it is denoted by \emph{$M_{\mathscr{S}}$}.
The values $\{M(w)\}_{w \in \Sigma^{*}}$ of $M_\mathscr{S}$ are the \emph{Markov parameters} of $\mathscr{S}$.
From \cite{MPLBJH:Real,PetreczkyRealChapter} it then follows that two dLSSs have the same input-output behavior, if and only if
their Markov functions are equal. 
%Moreover, the values of the sub-Markov function (sub-Markov parameters) can be
%determined from the input-output behavior.
%Sub-Markov parameters are analogous to Markov-parameters of linear systems, and just like in the linear case, a Ho-Kalman-like realization algorithm can be formulated which computes a dLPV-SSA from sub-Markov parameters \cite{PetreczkyLPVSS,RolandAbbas}. 
%For a function $M:\Sigma^{*} \rightarrow \mathbb{R}^{n_\mathrm{y} \times \un}$, we say that the dLSS $\mathscr{S}$ is a \emph{realization} of
%$M$, if $M$ equals the Markov function of %$\mathscr{S}$, i.e., $M=M_{\mathscr{S}}$.
%We call two dLSS input-output equivalent if they realize the same sub-Markov function.
For dLSS \eqref{eqn:LPV_SSA:det}
the integer $\nx$ is called the \emph{dimension} of $\mathscr{S}$, and we say that a 
dLSS is \emph{minimal}, if there exists no other dLSS of smaller dimension which realizes the same Markov-function.
%$\nx$ is larger than the dimension of $\mathscr{S}^{'}$ and their sub-Markov parameters are equal, i.e. $M_{\mathscr{S}}=M_{\mathscr{S}}^{'}$. 
%We call a dLSS $\mathscr{S}$ a \emph{minimal  realization} of a Markov function, if 
%$\mathscr{S}$ is minimal and it is a realization of $M$. 
From \cite{MPLBJH:Real,PetreczkyRealChapter}, it follows that a dLSS is minimal if and only if it is span-reachable and observable, and the latter
properties are equivalent to rank conditions of the extended reachability and observability matrices \cite[Definition 23, Theorem 1, Theorem 2]{PetreczkyRealChapter}.
Furthermore, any dLSS  can be transformed to a minimal one  with the same Markov function, using
a minimization algorithm \cite[Procedure 3 1]{PetreczkyRealChapter}\footnote{\cite[Corollary 1]{PetreczkyRealChapter} should be applied with zero initial state}.
For a more detailed discussion see \cite{PetreczkyRealChapter}. 
Moreover, any two minimal dLSS which are input-output equivalent are isomorphic, i.e. they are 
related by a linear change of coordinates
\cite[Theorem 1]{PetreczkyRealChapter}.

%\subsection{Correlation analysis} % finding an asLSS of  $(\yb^d,\p)$}
%\label{real:alg1}
%In this section, we describe an adaptation of the \emph{correlation
	%Manalysis} (CRA) method \cite{CoxIFAC,CoxLPVSS} for finding a stationary LSS representation of $\yb^d$ with noise process
%$\ub$. 

\subsection{Relationship between sLSSs and dLSSs}
  %As it was mentioned above, sLSSs can be viewed
  %as GBSs in the terminology of \cite{PetreczkyBilinear}, such that one of the %components of the noise process is the control input $\bu$. 
  The idea behind relating sLSSs and dLSSs
  is to express the covariances of outputs of sLSSs
  as Markov parameters of dLSSs. 
%In order to present the relationsip between sLSSs and dLSSs, 
  To this end, we need to recall 
  \cite{LCSSArxive} the 
 decomposition of the output $\y$
 %\subsection{Decomposition of the  output of  LSS representation}
%into deterministic and stochastic components}
%\label{sect:decomp}
%In this section we will address the existence of a LSS representation of $(\yb,\ub,\p)$. To this end, first we show that the output process of  a sLSS  admits a decomposition 
into deterministic and stochastic parts, and
%which depends only on $\bu$ and a \emph{stochastic} part which depends only on the noise of the data generator. 
%The deterministic part depends only on the input process, while the stochastic part depends only on the noise process.
%In order to present what follows,
%for which we need recall 
the notion of orthogonal projection %$E_l$ 
from \cite{LCSSArxive}.
\begin{Notation}[$E_l$]
	\label{hilbert:notation}
	Recall from \cite{Bilingsley} that
    the set $\mathcal{H}_1$ of real valued square integrable
	random variables % denoted by $\mathcal{H}_1$, %taking values in 
    is a Hilbert-space %$\mathcal{H}_1$ 
    with the scalar product $\langle \mathbf{z}_1,\mathbf{z}_2\rangle=E[\mathbf{z}_1\mathbf{z}_2]$.
    If $M$ is a  closed subspace   of $\mathcal{H}_1$, then 
    denote by $E_l[\mathbf{h}\!\! \mid \!\! M]$
    the orthogonal projection, in the usual sense for Hilbert-spaces,  of $\mathbf{h} \in \mathcal{H}_1$
    onto $M$.
    %We denote this Hilbert-space by %$\mathcal{H}_1$. 
	Let $\mathbf{z}=(\z_1,\ldots,\z_k)^T$ 
    be a $k$-dimensional %square integrable %\emph{vector-valued} 
	random variable 
 %taking its values in $\mathbb{R}^k$, 
 such that $\z_i$
    belongs to $\mathcal{H}_1$ for all $i=1,\ldots,k$. 
  %  Let . 
	The  orthogonal projection of $\mathbf{z}$ onto $M$, \emph{denoted by $E_l[\mathbf{z}\!\!\mid\!\! M]$},
	is defined as $k$-dimensional
     random variable
    $(E_l[\z_1\!\! \mid \!\!M],E_l[\z_2 \!\!\mid\!\! M], \ldots, E_l[\z_k\!\!\! \mid \!\! M] )^T$.
 %vector-valued 
    %square-integrable 
     %random variable $\mathbf{z}^{*}=\begin{bmatrix} \mathbf{z}_1^{*},\ldots,\mathbf{z}_k^{*} \end{bmatrix}^T$ such that $\mathbf{z}_i^{*} \in M$ is the orthogonal projection, in the usual sense used for Hilbert-spaces,  of $z_i$  onto $M$.
     %coordinate $\mathbf{z}_i$ of $\mathbf{z}$, viewed as an element of
     %$\mathcal{H}_1$ onto $M$. 
     %viewed
     %usually defined for Hilbert spaces. 
	If $\mathfrak{S}$ is a subset of vector valued random variables, coordinates of which all belong to $\mathcal{H}_1$, and 
	   $M$ is generated by the coordinates of the elements of $\mathfrak{S}$,
    %i.e. 
	%$M$ is the smallest (with respect to set inclusion) closed subspace of %$\mathcal{H}_1$
	%which contains the set $\{ \alpha^Ts %\mid  s \in \mathfrak{S}, \alpha \in %\mathbb{R}^p\}$. 
    then instead of $E_l[\z \mid M]$ we use \( E_{l}[\mathbf{z} \mid \mathfrak{S}] \). 
	%to denote the projection of $z$ to $M$.
	%Notice that $\mathbf{z}^{*}$ is uniquely determined by the following property: $E[(\mathbf{z}-\mathbf{z}^{*})\x^T]=0$ for all $\x \in Z$.
	% That is, $M$ is the smallest (with respect to set inclusion) closed linear subset such that for any $\mathbf{z}=(\mathbf{z}_1,\ldots,\mathbf{z}_p) \in Z$, $M$ contains the scalar random variables $\mathbf{z}_1,\ldots,\mathbf{z}_p$.
	%Notice also that one can interpret $E_{l}[\mathbf{z} \mid Z]$ as the best approximation (prediction) 
	%of $\mathbf{z}$ in terms of (infinite) linear combinations of elements~of~$Z$. 
\end{Notation}
Intuitively, %the orthogonal projection 
$E_l[\mathbf{z} \mid \mathfrak{S}]$ is \emph{minimal variance linear prediction} of $\mathbf{z}$ based on the elements of $\mathfrak{S}$.

The \emph{deterministic component} $\yb^d$ 
and the \emph{stochastic component} $\yb^s$
of $\yb$ are 
defined as 
	\begin{equation}
		\label{decomp:outp:eq1}
  \begin{split}
		&  {\yb}^d(t)=E_l[\yb(t) \mid \{\z_w^{\ub}(t)\}_{w \in \Sigma^{+}} \cup \{\ub(t)\}], \\
	%\end{equation}
	%and the \emph{stochastic component} of $\yb$ is %defined as 
	%\begin{equation}
		%\label{decomp:outp:eq2}
		& \yb^s(t)=\yb(t)-\yb^d(t).
    \end{split}
	\end{equation}
%\end{Definition}
Intuitively, $\y^d$  is  the best  linear 
prediction of $\y$ 
%which is
based on the present and past  values of $\bu$ multiplied by indicator functions of past discrete modes.
%and non-linear in the past values of $\p$.
%if the coefficients of the linear prediction depend on the past values of the switching signal. 
%From the definition it follows thatv\( \yb(t) = {{\yb}}^{d}(t) + {{\yb}}^{s}(t)\).
%
%i.e., the process $\yb(t)$ can be represented as the sum of its deterministic and  stochastic
%components. 
%From \cite{MejariLPVS2019} 
In fact, by \cite[Theorem 1]{LCSSArxive}, 
 $\y^d$  (resp. $\y^s$)
 is the output of
the asLSS obtained from \eqref{eq:aslpv},
by setting $K_{\sigma}=0$ (resp.
$B_{\sigma}=0$), for all $\sigma \in \Sigma$
%the deterministic part 
%by considering $\vb=0$ and
%viewing $\bu$ as noise,  and 
%$\y^s$ is the output of
%the asLSS obtained from \eqref{eq:aslpv}, by taking $\bu=0$.
%and viewing $\vb$ as noise.

%Let us define the map  $\incov: \Sigma^{*} \rightarrow \mathbb{R}^{\ny \times \un}$ as follows
%\begin{equation}\label{eqn:input_covar}
%	\incov(w) =\
    %left\{\begin{array}{lr} 
 %     \frac{1}{\sqrt{p_w}} E[\yb^d(t)(\zwu(t))^{T} ] \varu^{-1} 
      %\ \ &  \forall w \in \Sigma^{+}  \\
	%	E[\yb^d(t)\ub^T(t)]\varu^{-1} &w=\epsilon 
	%\end{array}\right. 
%\end{equation}
%where $\varu \!=\! \mathrm{var}(\ub)$, i.e., 
%$\incov{w}$ are the covariances of inputs and 
%of the deterministic components of the outputs. 
%Note that $E[\y(t)\u^T(t)]=E[\y^d(t)\u^T(t)]$

%Similarly, we can define \emph{covariance sequence} $\Psi_{\ybs}: \Sigma^{*} \rightarrow \mathbb{R}^{\ny \times \ny}$, where $\Psi_{\ybs}(\epsilon) =I_{\ny}$, and for all $w \in \Sigma^{+}$,
%\begin{equation}\label{eqn:covseq}
%	\covseq(w) = %\left\{\begin{array}{rl}
%	E[\ybs(t) (\zwys(t))^{T} ]  %& w \in \Sigma^{+} \\
	%I_{\ny} &  w=\epsilon 
	%             \end{array}\right. 
%M\end{equation}
%Intuitively, the elements $\covseq{w}$ are autocovariances of 
%the stochastic component of $\y^s$. 

%To make the discussion above more clear, we introduce
%the following combined 
Define the Markov function
$M_{\yb,\ub}: \Sigma^{*} \rightarrow \mathbb{R}^{\ny \times (\ny + \un)}$ by
\begin{equation} 
\label{min:markov:eq1}
   M_{\yb,\ub}(w)=\begin{cases} \begin{bmatrix} \Lambda^{\y^d,\bu}_{w} Q_{\bu}^{-1},  & \Lambda_w^{\y^s} \end{bmatrix} & w \ne \epsilon 
   \vspace{0.2cm} \\
   %& \\
   \begin{bmatrix} \Lambda^{\y^d,\bu}_{\epsilon} Q_{\bu}^{-1},  & I_{\ny} \end{bmatrix}
   & w=\epsilon 
   \end{cases}
\end{equation}
It can be shown that 
%$\Lambda^{\y^{d},\bu}_{\sigma v} \varu^{-1}=\sqrt{p_{\sigma v}} CA_vB_{\sigma}$, 
%$\Lambda^{\y^d,\bu}_{\epsilon}=D$, i.e., 
the first 
components of $M_{\yb,\bu}$
can be expressed as Markov parameters of a suitably defined dLSS.
Since, by \cite{LCSSArxive}, $\y^s$ is the output of an asLSS, and hence of a GBS in terminology of \cite{PetreczkyBilinear}, then by \cite[Lemma 4]{PetreczkyBilinear} $\Lambda_w^{\y^s}$
can be shown to be the Markov function of a suitable dLSS. 
That is, 
%sLSS realizations of $(\y,\bu,\btheta)$ give rise to dLSS realizations of 
$M_{\yb,\ub}$ has a dLSS realization.

%\begin{Definition}[dLSS associated stocLSS]\label{defn:det_asco_stationary}
Formally, define the process $\xb^s(t)=\x(t)-E_l[\x(t) \mid \{\z_w^{\ub}(t)\}_{w \in \Sigma^{+}} \cup \{\ub(t)\}]$
and the matrices 
\begin{equation}
\label{min:eq-1}
\begin{split}
& \Bs_{\sigma}=\sqrt{p_{\sigma}} (A_{\sigma} T^{\xb^s,\xb^s}_{\sigma,\sigma} C^T + K_{\sigma} T^{\vb,\vb}_{\sigma,\sigma} F^T) \\
\end{split}
\end{equation}
Then  define the \emph{dLSS associated with sLSS} as
$$\mathscr{S}_{\mathcal{S}}=(\{\sqrt{p_\sigma} A_\sigma,\begin{bmatrix} \sqrt{p_\sigma} B_\sigma & \Bs_\sigma \end{bmatrix} \}_{\sigma=1}^{\pdim},C,\begin{bmatrix} D & I_{\ny} \end{bmatrix}), $$
%for all $\sigma \in \Sigma$,
%$\Bs_{\sigma}$ is defined as follows.
%$\{\IBs_i\}_{\sigma=1}^{\pdim}$ is defined as follows. Let $P_{\sigma}=p_{\sigma} T^{\x^s,\x^s}_{\sigma},\sigma}$, $\sigma \in \Sigma$, where
	 %$$P_{i}=E[\xb^s(t)(\xb^s(t))^{T}\p^2_i(t)]$$,
Note that by \cite[proof of Theorem 1]{LCSSArxive} $\xb^s(t)$ is the state
of the asLSS obtained from $\mathcal{S}$
by considering $B_{\sigma}=0$, and hence 
$\xb^s$ is a ZMWSII process and 
$T^{\xb^s,\xb^s}_{\sigma,\sigma}$ is well-defined.
%where $Q_i=E[\vb(t)\vb(t)\p_i^2(t)]$ for %M$i=1,\ldots,\pdim$.
%From~\cite[Lemma 4]{PetreczkyBilinear} the following Lemma follows easily.
%\begin{Lemma}[Lemma 4 from \cite{PetreczkyBilinear}]
%\label{lemma:stoch-real}
% The dLSS $(\{\sqrt{p_{\sigma}}A_{\sigma},\Bs_{\sigma}\}_{\sigma \in \Sigma},C,I)$ is a realization of
% $\covseq$.
%such that the sub-Markov function $M_{\mathscr{S}_{det}}$ of $\mathscr{S}_{det}$
% equals $\covseq$, \emph{i.e}, 
%\( \covseq(w) = \sqrt{p_w} CA_{w}\Bs, ~ \forall \wordSet{*} \). 
 %\[ \expect{\zwsigs(t)(\zwsigs(t))^{T}} = \frac{1}{\psig} (C\Psig C^{T} + \Qsig), ~ \forall \sigSet. \]
%\end{Lemma}
 %\end{Definition}
\begin{Lemma}
\label{min:col:lem1}
 If $\mathcal{S}$ is a sLSS realization of $(\yb,\ub,\btheta)$, then the associated dLSS $\mathscr{S}_{\mathcal{S}}$ is a realization of $M_{\yb,\ub}$. 
 \end{Lemma}
 The proof follows from \cite[Lemma 4]{PetreczkyBilinear}, \cite[Theorem 2]{LCSSArxive},
 it is presented in \cite{rouphael2024minimal}.
 %and some tedious computations,
 %and follows the lines of proof for
 %\cite[Lemma %1]{MejariLPVS2019}.
%Since asLSSs are special cases of GBSs from \cite{PetreczkyBilinear}, it follows that $\y^s$ can be realized by GBS and hence, by \cite[Lemma 4]{PetreczkyBilinear}
%that covariances  of  can be expressed through the  the product of  system matrices of \eqref{eq:aslpv}. A somewhat tedious calculation also reveals that 
\begin{Remark}[Computing $\Bs_{\sigma}$]
\label{gi:comp}
 In order to compute  %$\{P_i\}_{i=1}^{\pdim}$ and hence 
 $\Bs_{\sigma}$, we can use
 that by \cite[proof Theorem 1]{LCSSArxive}
 $\xb^s$ is the state 
 of the asLSS obtained from 
 \eqref{eq:aslpv}
 by setting $B_{\sigma}\!\!=\!\!0$, $\sigma \in \Sigma$.
 Hence, by 
\cite[Lemma 5]{PetreczkyBilinear}, the covariance of $\xb^s$ satisfies
$p_{\sigma} T^{\xb^s,\xb^s}_{\sigma,\sigma} = \lim_{N \rightarrow \infty} P_{\sigma}^N$, where 
$P_{\sigma}^0=0$ and 
 \(
    %\begin{split}
    %& P_i = \lim_{N \rightarrow \infty} P_i^N, ~
    P_{\sigma}^{N+1}=p_{\sigma} \sum_{\sigma_1 \in \Sigma}  (A_{\sigma_1}P_{\sigma_1}^NA_{\sigma_1}^T + p_{\sigma_1} K_{\sigma_1}T^{\vb,\vb}_{\sigma_1,\sigma_1} K_{\sigma_1}^T)
   %\end{split}
 \).
  %In particular, $P_{\sigma}=p_{\sigma} \sum_{j \in \Sigma} (A_jP_jA_j^T + p_j K_jT^{\vb,\vb}_{j,j}K_j^T)$ . 
\end{Remark}
Conversely, we can associate with any dLSS realization
of the Markov function $M_{\y,\bu}$ a sLSS realization
of $(\y,\bu,\btheta)$. In order to define the latter,
we need to specify its noise process, which happens to
be the \emph{innovation noise} $\e^s$ 
of $\y^s$ as defined
in \cite[eq. (16)]{PetreczkyBilinear}. By \cite[Theorem 2]{LCSSArxive}, %$\e^s$ can be defined
	\begin{equation} 
		\label{decomp:lemma:innov}
     \begin{split}
		& \eb^s(t) =\yb(t)-\hat{\yb}(t)  \\ %\yb^s(t)-E_l[\yb^s(t) \mid \{\z^{\yb^s}_w (t)\}_{w \in \Sigma^{+}}] = \\
  %\begin{equation} 
	%\label{decomp:lemma:innov2
		 & \hat{\yb}(t)=E_l[\yb(t) \mid \{\z^{\yb}_w (t), \z^{\ub}_w(t) \}_{w \in \Sigma^{+}} \cup \{\ub(t)\}] 
	\end{split}
	\end{equation}
That is, $\eb^s(t)$ is the prediction error
of the best linear predictor $\hat{\y}(t)$ of $\y(t)$
given the products of past outputs, inputs and discrete modes, i.e.,
$\e^s$ is direct extension of
the classical
innovation noise.
We say that
$\yb$ is \emph{full rank},
if for all $\sigma \in \Sigma$,
the covariance $T^{\eb^s,\eb^s}_{\sigma,\sigma}$
%Q_i=E[\eb^s(t)(\eb^s(t))^T\p_i^2(t)]$ 
is invertable. 
%Note that, the process  $\eb^s(t)$. 
This is a direct extension of the classical notion of a full rank process. 
%The definition is a direct extension of
%the standard definition of innovation noise for LTI
%systems. 
%With these definition in mind, let us define the correspondence between dLSSs and sLSSs.
 %Let $\mathcal{S}$ be a sLSS o%f $(\yb,\ub,\p)$ 
 % of the form  \eqref{eq:aslpv}. 
%Conversely, one can associated with a dLSS a
%sLSS, defined as follows.
%\begin{Definition}[sLSS associated with dLSS]\label{defn:stoc_asco_stationary}
 If $\yb$ is full rank,  for an observable dLSS $$\mathscr{S}=(\{\hat{A}_\sigma,\begin{bmatrix}  \hat{B}_\sigma & \hat{G}_\sigma \end{bmatrix}\}_{\sigma=1}^{\pdim}, \hat{C},\begin{bmatrix} D & I_{\ny} \end{bmatrix})$$ 
 for which $\sum_{\sigma=1}^{\pdim} \hat{A}_\sigma \otimes \hat{A}_\sigma$ is a Schur matrix\footnote{all its eigenvalues are inside the unit disk}, 
 define the \emph{sLSS $\mathcal{S}_{\mathscr{S}}$ associated with $\mathscr{S}$ as } 
 $$ \mathcal{S}_{\mathscr{S}}=(\{\frac{1}{\sqrt{p_i}} \hat{A}_\sigma,\frac{1}{\sqrt{p_\sigma}} \hat{B}_\sigma, \hat{K}_\sigma\}_{\sigma=1}^{\pdim},\hat{C},\hat{D},I_{\ny},\eb^s), $$
where $\hat{K}_\sigma=\lim_{\mathcal{I} \rightarrow \infty} \hat{K}_\sigma^{\mathcal{I}}$, and $\{\hat{K}_\sigma^{\mathcal{I}}\}_{\sigma \in \Sigma,\mathcal{I} \in \mathbb{N}}$ satisfies 
	\begin{equation}
		\label{statecov:iter}
		\begin{split}
			%\hat{P}^{i+1}_{\sigma} \!&=\! \sum \limits_{\sigma_{1} \in \Sigma, \sigma_{1} \sigma \notin L } \psig \left( \hat{A}_{\sigma_{1}} \hat{P}^{i}_{\sigma_{1} }  \hat{A}^{T}_{\sigma_{1}}  + \hat{K}_{\sigma_{1}} \hat{Q}^{i}_{\sigma_{1} }  \hat{K}^{\mathcal{I}})^{T}_{\sigma_{1}}  \right)\\
			& 	\hat{P}^{\mathcal{I}+1}_{\sigma} \!\!=\!\! p_{\sigma}\!\! \sum \limits_{\sigma_{1} \in \Sigma} 
    \frac{1}{p_{\sigma_1}}\left( \hat{A}_{\sigma_{1}} \hat{P}^{\mathcal{I}}_{\sigma_{1}}  \hat{A}_{\sigma_{1}}^{T}  + \hat{K}_{\sigma_{1}}^{\mathcal{I}} \hat{Q}^{\mathcal{I}}_{\sigma_{1} }  (\hat{K}^{\mathcal{I}}_{\sigma_{1}})^T  \right)\\
			& \hat{Q}^{\mathcal{I}}_{\sigma} = \psig T^{\y^s,\y^s}_{\sigma,\sigma} - \hat{C}  \hat{P}^{\mathcal{I}}_{\sigma } (\hat{C})^{T}  \\ 
			& \hat{K}^{\mathcal{I}}_{\sigma} = \left( \hat{\Bs}_{\sigma} \sqrt{\psig} - \frac{1}{\sqrt{\psig}}\hat{A}_{\sigma} \hat{P}^{\mathcal{I}}_{\sigma}  (\hat{C})^{T} \right) \left(  \hat{Q}^{\mathcal{I}}_{\sigma} \right)^{-1}
		\end{split}
	\end{equation}
	with $\hat{P}^{0}_{\sigma} =0$. 
%\end{Definition} 
  Note that the noise process of $\mathcal{S}_{\mathscr{S}}$
  is the innovation noise, and
   its noise and  state covariances satisfy
   %$\hat{\xb}$  
   %and the noise processes $\eb^s$ 
   %of $\mathcal{S}_{\mathscr{S}}$ satisfy 
        \begin{equation}
	\label{lemma:stoch-real1:eq1}
	  \begin{split}
      p_{\sigma} T^{\eb^s,\eb^s}_{\sigma,\sigma}=
      \lim_{\mathcal{I} \rightarrow \infty} \hat{Q}_{\sigma}^{\mathcal{I}}, ~ \quad 
	    %& E[\eb^s(t)(\eb^s(t))^T\p^2_{\sigma}(t)]=\hat{Q}_{\sigma}=\lim_{\mathcal{I} \rightarrow \infty} \hat{Q}_{\sigma}^{\mathcal{I}} \\
	  %  & 
     %E[\hat{\xb}(t)\hat{\xb}^T(t)\p^2_{\sigma}]=\hat{P}_{\sigma}=
       p_{\sigma} T^{\hat{\xb},\hat{\xb}}_{\sigma,\sigma}=\lim_{\mathcal{I} \rightarrow \infty} \hat{P}_{\sigma}^{\mathcal{I}}
	\end{split}
      \end{equation}
%\eqref{statecov:iter}
	%and $\eb^s$ is as in \eqref{decomp:lemma:innov2}.
	%Moreover, from \cite{} it follows that
 where $\hat{\xb}$ is the unique state process of $\mathcal{S}_{\mathscr{S}}$.
The convergence of \eqref{statecov:iter} and \eqref{lemma:stoch-real1:eq1} follow from the proof of \cite[Theorem 3]{PetreczkyBilinear}.
 \begin{Lemma}
 \label{min:col:lem2}
If $\yb$ is full rank and $(\yb,\ub,\btheta)$ has a realization by sLSS, and $\mathscr{S}$ is a minimal dLSS realization of $M_{\yb,\ub}$, then the associated sLSS $\mathcal{S}_{\mathscr{S}}$ is a sLSS realization of $(\yb,\ub,\btheta)$. 
\end{Lemma}
The proof of Lemma \ref{min:col:lem2} is presented in \cite{rouphael2024minimal}, it relies on
\cite[Theorem 2]{LCSSArxive} and \cite[Theorem 3]{PetreczkyBilinear}.

\subsection{Main result on minimal sLSSs in innovation form}
\label{sec:min}
%the extended $\nx-1$-step reachability and observability matrices of $\mathscr{S}_{\mathcal{S}}$ satisfy $\rank \mathcal{R}_{\nx-1}=\nx$ and $\rank \mathcal{O}_{\nx-1}$.
A sLSS of the form \eqref{eq:aslpv} is said to be in \emph{innovation form}, if its noise process $\vb$ equals $\eb^s$ from  
\eqref{decomp:lemma:innov} and $F=I_{\ny}$. 
Note that the sLSS associated with a dLSS is in innovation form. % then it is a generator
%is driven by past outputs
%and inputs
 The sLSS \eqref{eq:aslpv} in innovation form is a predictor 
 %which generates 
%past outputs and inputs to the
%t
%$\hat{\y}(t)=E_l[\yb(t) \mid \{\z^{\yb}_w (t), \z^{\ub}_w(t) \}_{w \in \Sigma^{+}} \cup \{\ub(t)\}]$
%of the output $\y(t)$ 
%using past outputs and inputs, 
\begin{equation*}
\label{eq:alpsv:innov}
\begin{split}
&\xb(t\! +\! 1)\!\!=%\hat{A}_{\btheta(t)} %=
\!\!(A_{\btheta(t)}\! -\! K_{\btheta(t)}C)\xb(t)+B_{\btheta(t)}\ub(t)+K_{\btheta(t)}\y_u(t) \\
%[B_{\btheta(t)},K_{\btheta(t)}] \begin{bmatrix} \bu(t) \\ \y(t) \end{bmatrix} \\
& \hat{\y}(t)\! = \! C\xb(t)+D\ub(t), ~ \y_u(t)=\y(t)-D\bu(t) 
\end{split}
\end{equation*}
%i.e., agenerator 
%where $\hat{\y}(t)=E_l[\yb(t) \mid \{\z^{\yb}_w (t), \z^{\ub}_w(t) \}_{w \in \Sigma^{+}} \cup \{\ub(t)\}]$
%is the prediction of $\y(t)$.
which generates the best linear prediction $\hat{\y}(t)$ 
%\eqref{decomp:lemma:innov} 
of the output $\y(t)$,
see \eqref{decomp:lemma:innov},
while being driven by past outputs and inputs.

We say that a sLSS realization of $(\y,\bu,\btheta)$ is \emph{minimal},
if it has the smallest dimension among all sLSS realizations of $(\y,\bu,\btheta)$. 
In order to study minimality, we need the concepts of reachability and observability.
We call a sLSS $\mathcal{S}$ \emph{reachable} (resp. \emph{observable}), if the associated 
dLSS $\mathscr{S}_{\mathcal{S}}$ is span-reachable, (resp. observable) according to the definition
of \cite[Definition 19]{PetreczkyRealChapter}, applied to zero initial state. 

In order to investigate uniqueness of
minmal sLSS in innovation form, we need the concept
of isomorphism
Let  $\mathcal{S}$ and $\widetilde{\mathcal{S}}$ be sLSSs realization of $(\yb,\ub,\btheta)$ in innovation form
such that $\mathcal{S}$ is of the form \eqref{eq:aslpv}
and
%(\{A_i,K_i,B_i\}_{i=0}^{\pdim},C,D,\eb^s)$ and  
%Let 
%\begin{equation}
%\label{lemma:min1:eq2}
$\widetilde{\mathcal{S}}=(\{\tilde{A}_{\sigma},\tilde{B}_{\sigma}, \tilde{K}_{\sigma}\}_{\sigma=1}^{\pdim},\tilde{C},\tilde{D},I_{\ny},\eb^s)$. We say that $\mathcal{S}$ and $\widetilde{S}$ are \emph{isomorphic}, if there exists a nonsingular matrix $T$ such that
\begin{equation}
\label{lemma:min1:eq1}
\begin{split}
  &\tilde{A}_{\sigma}=TA_{\sigma}T^{-1}, \tilde{K}_{\sigma}=TK_{\sigma}, ~ \tilde{B}_{\sigma}=TB_{\sigma}, \quad \sigma \in \Sigma  \\
  &\tilde{C}=CT^{-1}, D=\tilde{D}
\end{split}
\end{equation} 
%\end{equation} 
%be a sLSS of $(\yb,\ub,\p)$ in forward innovation form. 
%Finally, we need the following technical condition: 
%where $Q_i=E[\eb^s(t)(\eb^s(t))^T\p_i^2(t)]$ and $i=1,\ldots,\pdim$
%nd if $\yb$ is full rank, then in \eqref{lemma:min1:eq1}, $\tilde{K}_iQ_i=TK_iQ_i$ can be replaced by $\tilde{K}_i=K_i$

%We can now  characterize minimality of sLSSs in innovation form follows.
%\begin{Procedure}[Constructing a minimal sLSS]
%\label{proc:slpv:min}
%Assume $\mathcal{S}$ is a sLSS of $(\yb,\ub,\p)$.
%Consider $(\{A_i,K_i,B_i\}_{i=1}^{\pdim},C,D)$ and 
%\begin{enumerate}
%\item
% Construct the dLSS $\mathscr{S}_{\mathcal{S}}$ associated with $\mathcal{S}$
  %be a minimal dLSS realization of $\incov$ and construct the sLSS $\Gamma_d=(\{\hat{A}^{d}_{i}, \hat{B}^{d}_{i} \}_{i=1}^{\pdim}, \hat{C}^{d}, \hat{D}^{d},\ub)$ of $(\yb^d,\p)$ from $\mathscr{S}_d$ as described in Lemma \ref{thm:cra}.
%\item  
%Construct the asLSS $\mathcal{S}_s$from a minimal dLSS realization of $\covseq$ as described \eqref{lemma:stoch_real1:eq}, Lemma \ref{lemma:stoch-real1}.
%\item 
% Construct the sLSS  $\mathcal{S}$
% of $(\yb,\ub,\p)$ from $\Gamma_d$ and $\Gamma_s$ described in Lemma \ref{decomp:lemma:inv}.
%\item 
% Construct a minimal dLSS $\mathscr{S}_m$  from $\mathscr{S}_{\mathcal{S}}$ by using procedure from \cite[Corollary 1]{PetreczkyLPVSS}.
%\item 
%  Return the  sLPVS-SSA $\mathcal{S}_m=\mathcal{S}_{\mathscr{S}_m}$ associated with $\mathscr{S}_m$.
%\end{enumerate}
%\end{Procedure}
\begin{Theorem}[Main result: minimal innovation form]
\label{theo:min}
 Assume that $(\yb,\ub,\btheta)$ has a sLSS realization and that $\yb$ is full rank w.r.t. $\bu$ and $\btheta$. Then the following holds.
\begin{enumerate}
 %$(\yb,\ub,\p)$ has a sLSS in forward innovation form which is span-reachable and observable. 
\item
\label{theo:min:1}
A sLSS realization of $(\yb,\ub,\btheta)$ is minimal if and only if it is reachable and observable. 
\item
\label{theo:min:2}
  If a dLSS $\mathscr{S}_m$ is a minimal realization of $M_{\yb,\ub}$, then the associated sLSS $\mathcal{S}_{\mathscr{S}_m}$ is a minimal sLSS  realization of $(\yb,\ub,\btheta)$ in innovation form.
\item 
\label{theo:min:3}
  There exists a minimal sLSS realization of 
  $(\y,\bu,\btheta)$ in innovation form.
\item
\label{theo:min:4}
Any two minimal sLSS of $(\yb,\ub,\btheta)$ in innovation form are isomorphic.  
\end{enumerate}
\end{Theorem}
The proof is presented in \cite{rouphael2024minimal}, 
it relies on the fact that any dLSS can be transformed to a minimal one, 
that minimal input-output dLSS are isomorphic and on Lemma \ref{min:col:lem1}-\ref{min:col:lem2}. 

Note that observability and reachability of the associated dLSS can be checked using simple rank conditions of suitable matrices \cite{PetreczkyRealChapter}.
That is, our conditions for minimality can be checked numerically. Moreover, a sLSS can be transformed effectively to a minimal one in innovation form as 
follows: first we compute the associated dLSS, we minimize it using \cite[Procedure 3]{PetreczkyRealChapter} to obtain a minimal dLSS realization of $M_{\y,\bu}$, and then we compute the sLSS associted with the latter minimal dLSS. 
\begin{Remark}[Relationship with \cite{LCSSArxive}]
\label{lcss:rem}
For a sLSS \eqref{eq:aslpv}, the minimality condition of \cite{LCSSArxive} corresponds to
reachability and observability
of the dLSS which is obtained
from the dLSS \eqref{eqn:LPV_SSA:det}
associated with 
\eqref{eq:aslpv} by 
deleting certain columns of $\mathcal{B}_{\sigma}$.
% $(\{A_{\sigma},G_{\sigma}\}_{\sigma=1}^{\pdim},C,D)$,
The latter is sufficient, but not
necessary for the reachability and observability of the 
%The latter implies minimality of the 
dLSS associated with \eqref{eq:aslpv}.
Likewise, \cite{LCSSArxive}
shows isomorphism between 
sLSSs in innovation form which satisfy the above sufficient condition for minimality, and for which the image of $B_{\sigma}$
belongs to the image of $K_{\sigma}$. In contrast, Theorem \ref{theo:min} establishes
isomorphism for a much wider class of sLSSs in innovation form. 

\end{Remark}

\section{Minimial covariance realization algorithm }
\label{sect:min:alg}
In this section we present a Ho-Kalman-like algorithm
for computing a minimal sLSS realization innovation
form from $M_{\y,\bu}$. 
Ot this end, we first we recall
the Ho-Kalman realization algorithm for dLSS. 
%We
%then use it to formulate a covariance realization algorithm for sLSSs. 

\subsection{Reduced basis Ho-Kalman algorithm for dLSS}
%Here in this section, 
Below we recall from \cite{CoxLPVSS} an adaptation of the Ho-Kalman-like algorithm. %which 
%uses  small Hankel-matrices.
%uses a selection of rows and columns of the Hankel-matrix 
%to compute a dLSS from
%Markov-parameters. 
%To this end, %we first present the notion of \emph{$n$-selection}. 
%Let us define the set $\Sigma^{n}$ as the set of all words $w \in \Sigma^{*}$ of length less than or equal to $n$, i.e., $\Sigma^{n} =\{ w \in \Sigma^{*} \mid |w| \leq n\} $. 
%\begin{Definition}[Selection]
%following \cite{CoxLPVSS},
Define a \emph{$(n,n_{\mathrm y},n_{\mathrm u})$-selection (selection for short)} as a pair of word-index sets
$\left(\alpha, \beta \right)$
such that,
%each containing $n$ elements,  
%such that 
  %\citep{PetreczkyVidalArxiv}
	%\begin{enumerate}
	%	\item 
  %$\alpha \subseteq \Sigma^{n} \times \{1,2,\cdots, \ny \} $ and $\beta \subseteq \Sigma \times \Sigma^{n} \times \{1,2,\cdots,\un \}$, and 
	%$\alpha$ and $\beta$ have exactly $n$ elements.
 % 
  %= n $, where $\mathrm{card}$ denotes cardinality of the set.
%	\end{enumerate}
	%When $\ny$ and $\un$ are clear from the context, we refer to $(n,\ny,\un)$-selections as $n$-selections,
	%and when $n$ is also clear from the context, we %use the term selection.
%end{Definition}
%In particular, we fix an ordering of
%$\alpha$ and $\beta$:
\begin{equation}
	\label{selection:enum}
	\begin{split}
		&\alpha = \{(u_{i},k_{i}) \}_{i=1}^{n}, ~ \beta = \{(\sigma_{j}, v_{j},l_{j})\}_{j=1}^{n},
	\end{split}
\end{equation}
for some $u_{i} \in \Sigma^{+}$, $|u_i| \le n$, 
$k_{i} \in \{1,2,\cdots,n_{\mathrm y} \}$, $\sigma_{j} \in \Sigma$, $v_{j} \in \Sigma^{+}$, $|v_j| \le n$, 
$l_{j} \in \{1,2,\cdots,n_{\mathrm u}\}$,
$i,j \in \{1,\ldots,n\} $ 
%\begin{Example}\label{eg:selection}
%	Consider $n\!=\!2$, number of outputs and inputs $\ny= \un= \!\!\!=\!\!\!2$, and scheduling signal dimension $\pdim \!\!\!=\!\!\!2$, we have, $\Sigma^{n} = \{\epsilon,1,2,11,12,21,22\}$. Then, one of the \emph{$n$-selection} pair  $\left(\alpha, \beta \right)$  can be chosen as, for e.g.,
%	$\alpha =\{\left(u_{1},k_1 \right), \left(u_2,k_2 \right)  \}= \{\left(\epsilon,1 \right), \left(11,2 \right)  \}$ and $\beta = \{\left(\sigma_{1},v_{1},l_1 \right), \left(\sigma_{2},v_2,l_2 \right)  \}=\{\left(1,21,1 \right), \left(2,22,2 \right)  \}$. 
%\end{Example}
%That is, $\alpha$ and $\beta$ can be thought of 
%as collection of $n$ sequences from $\Sigma^{+}$
%along with some auxiliary indices. 
Intuitively, the components of $\alpha,\beta$
which belong to $\Sigma^{+} \cup \Sigma$ determine the choice of  Markov parameters, while the other components determine a choice of
row and column indices of the chosen Markov parameters, see \cite{MejariLPVS2019} for an example. 
Formally, let $M:\Sigma^{*} \rightarrow \mathbb{R}^{n_{\mathrm y} \times n_{\mathrm u}}$ be a Markov function. 
%potential sub-Markov parameters \eqref{eqn:sub_markov} of an LSS.
%Let us now 
Define the Hankel matrix $\hankredu^{M} \in \mathbb{R}^{n \times n}$ and the matrices $\mathcal{H}_{\sigma, \alpha, \beta}^M \in \mathbb{R}^{n \times n}$, $\mathcal{H}_{\alpha,\sigma}^M \in \mathbb{R}^{n \times n_{\mathrm u}}$ and $\mathcal{H}^M_{\beta} \in \mathbb{R}^{n_{\mathrm y} \times n}$ as
follows as follows:
% in such a way that let us use the ordering of $\alpha$ and $\beta$, and for any 
%for any $i,j=1,\ldots,n$,
%the $(i,j)$-th element of %$\hankredu^{M}$ is of the form
\begin{align}	
	& \left[ \hankredu^{M} \right]_{i,j}  \!=\! \left[M(\sigma_{j}v_{j} u_{i} ) \right]_{k_i,l_j}, \quad i,j=1,\ldots,n \label{eqn:hankel_u} \\
 %\end{equation}
%where 
%Intuitively, the rows of $\hankredu^M$ are indexed by word-index pairs $\left(u_i,k_i \right) \in \alpha$, where $u_i \in \Sigma^{n}$ and $k_i \in \{1,\ldots, \ny\}$ and similarly, the columns of  $\hankredu^M$ are indexed by word-index pairs $\left(\sigma_{j}v_j,l_j \right) \in \beta$, where $\sigma_{j} \in \Sigma$, $v_j \in \Sigma^{n}$ and $l_j \in \{1,\ldots, \un\}$, and 
%the element of $\hankredu^M$ with the row indexed $(u_i,k_i)$ and column index $(\sigma_{j},v_j,l_j)$ is the $(k_i,l_j)$-th entry
%of $M(\sigma_{j}v_ju_i)$.  
%We define the 
%matrices
%µ$\sigma-shifted Hankel-matrix
%as follows: its $i,j$-th entry is given by 
%\begin{align}	 
& \left[\mathcal{H}_{\sigma, \alpha, \beta}^M\right]_{i,j} = \left[ M(\sigma_{j}v_j \sigma u_i) \right]_{k_i,l_j}, ~i,j=1,\ldots,n \label{eqn:hankelu:shift} \\
%\end{equation}
%and the matrices   
%as follows
%\begin{align}
	& \left[\mathcal{H}^M_{\alpha, \sigma }\right]_{i,j} \! = \! \left[ M(\sigma u_i)\right]_{k_i, j},  \ j\!=\!1,\ldots, n_{\mathrm u}, i\!=\!1,\ldots,n \label{eqn:hankel_u_2a}\\
	& \left[\mathcal{H}^M_{\beta}\right]_{i,j} \!=\! \left[ M(\sigma_j v_j) \right]_{i, l_j},  \ i=1,\ldots, n_{\mathrm y}, j=1,\ldots n
 \label{eqn:hankel_u_2b}
\end{align}
%where the subscript $(k,l)$
%$\left[M(\sigma_{j}v_{j} u_{i} ) \right]_{k_i,l_j}$ denotes the entry of $M(\sigma_{j}v_{j} u_{i} )$ on the
%$k_i$-th row and $l_j$-th %column, 
where 
$\left(u_{i},k_i\right) \in \alpha, \left(\sigma_{j},v_j,l_j\right) \in \beta$ are as in  the ordering in \eqref{selection:enum}.
Finally, let 
$\mathcal{L}_{\alpha,\beta}$ be the set of all
words $w \in \Sigma^{*}$ such that 
$M(w)$ occurs in one of the 
matrices, 
%Hankel-matrices 
\eqref{eqn:hankel_u} --
%\eqref{eqn:hankel_u_2a},%\eqref{eqn:hankel_u_2b}, 
\eqref{eqn:hankelu:shift}, i.e.,
%for $M=M_{\y,\bu}$.
%formally, 
\[ \mathcal{L}_{\alpha,\beta}=
   \{ u_i,\sigma_jv_j, \sigma_j v_ju_i, \sigma_jv_j\sigma u_i, \sigma u_i\mid \sigma  \in \Sigma\}_{i,j=1}^{n}
\]
%sequences $w \in \Sigma^{+}$ such that  either
% $w=ivu$ or $w=iv$ or $w=iu$ or $w=i v\sigma u$ for some words $v,u \in \Sigma^{*}$, $i,\sigma \in \Sigma$, $(u,k) \in \alphas $, $(v,l) \in \betas $ for some $k,l=1,\ldots,\ny$.
The promised algorithm is presented in Algorithm \ref{algo:Nice_select_true_deter_basic}.
%using Hankel matrices and selections.
\begin{algorithm}[h]
	\caption{Ho-Kalman for dLSSs}
 %Matrix computations  using Hankel matrices and $n$-selection}
	\label{algo:Nice_select_true_deter_basic}
	\textbf{Input}:  $(n,n_{\mathrm y},n_{\mathrm u})$-selection $\left(\alpha, \beta \right)$, 
    $\{M(w)\}_{w \in \mathcal{L}_{\alpha,\beta}}$, $M(\epsilon)$ 
    %Markov parameters
 %Hankel matrix $\hankredu^M$, %\eqref{eqn:hankel_u}; %shifted Hankel-matrix 
    %ankel_u_2a}-\eqref{eqn:hankel_u_2b} respectively, 
  %   and the matrix $M(\epsilon)$.
	\hrule\vspace*{.1cm}
	%\begin{enumerate}[label=\arabic*., ref=\theenumi{}]
		%\item  Defineee
    %\item
      Construct the matrices 
      $\mathcal{H}^M_{\alpha,\beta}$, $\mathcal{H}^M_{\sigma, \alpha, \beta} $ %\eqref{eqn:hankelu:shift}; Hankel matrices   
     $\mathcal{H}_{\alpha,\sigma}^M$ and $\mathcal{H}_{\beta}^M$, and set
     %defined in \eqref{eqn:h 
    %\item
   %Define the matrices
		\begin{align*}
			&  \hat{A}_{\sigma} = (\hankredu^M)^{-1}  \mathcal{H}_{\sigma, \alpha, \beta}^M, ~  \hat{B}_{\sigma} = (\hankredu^M)^{-1} \mathcal{H}^M_{\alpha,\sigma}, ~ \hat{C} =  \mathcal{H}^M_{\beta}. %\hat{D=M(\epsilon)
		\end{align*}
		%\begin{equation*}
		%\hankredu = U_{n} \Sigma_{n} V^{T}_{n}
		%\end{equation*}
		%\item %Compute $\hat{A}^{d}_{\sigma} \in \mathbb{R}^{n \times n}$, $\hat{B}^{d}_{\sigma} \in \mathbb{R}^{n \times \un}$ and $\hat{C}^{d} \in \mathbb{R}^{\ny \times n}$ %(Eq.~\eqref{eqn:B})
		%s follows: \\
		%Let $\hat{\mathcal{O}}_{\alpha}^{\dagger} = \Sigma_{n}^{-1/2}  U_{n}^{T} $ and $\hat{\mathcal{R}}_{\beta}^{-1} =  V_{n}\Sigma_{n}^{-1/2} $
		%\begin{align*}
		%\hat{A}_{\sigma} &=  \hat{\mathcal{O}}_{\alpha}^{-1} \mathcal{H}_{\sigma, \alpha, \beta} \hat{\mathcal{R}}_{\beta}^{-1} \  &\forall \sigSet  \\
		%\hat{B}_{\sigma} &= \hat{\mathcal{O}}_{\alpha}^{-1} \mathcal{H}_{\alpha,\sigma} \  &\forall \sigSet \\
		%\hat{C} &= \mathcal{H}_{\beta} 
		%\end{align*}
		%\item Compute $\hat{D}=M(\epsilon)$.
		%\begin{equation*}
		%\hat{D}^{d} =\incov{\epsilon}=\expect{\yb(t) \ub^{T}(t)} \varu^{-1}, 
		%\end{equation*}
	%\end{enumerate}
	%\vspace*{.1cm}
	\hrule 
 %\vspace*{.1cm}
	\textbf{Output}: dLSS $(\{\hat{A}_{\sigma}, \hat{B}_{\sigma} \}_{\sigma=1}^{\pdim}, \hat{C}, M(\epsilon) )$
 %\vspace*{}
\end{algorithm}
\vspace{-0.1cm}
\begin{Lemma}[\cite{CoxLPVSS}]
	\label{basis_red:lemma}
	%Let , and 
Assume that there exists a
minimal  dLSS realization of $M$
of dimension $n$.
Then there exists an $(n,n_{\mathrm y},n_{\mathrm u})$-selection $(\alpha,\beta)$ such that $\text{rank}(\hankredu^M) = n$. For any $(n,n_{\mathrm y},n_{\mathrm u})$-selection $(\alpha,\beta)$
for which $\text{rank}(\hankredu^M) = n$, 
%which is a realization of $M$. Then the dLSS  returned by 
        Algorithm \ref{algo:Nice_select_true_deter_basic}
        returns a minimal dLSS realization of $M$.
        %when applied to
	%the matrices $\hankredu^M$, $\mathcal{H}^M_{\sigma, \alpha, \beta}$, $\mathcal{H}_{\alpha,\sigma}^M$, $\mathcal{H}_{\beta}^M$ 
       %(\eqref{eqn:hankel_u}-\eqref{eqn:hankel_u_2b}) %and $M(\epsilon)$, 
       %is a minimal realization of $M$.
	%such that $M_{\hat{\mathscr{S}}}=M$, i.e.
	%\begin{equation}\label{eqn:det_matrix_estimates}
	%\hat{A}_{\sigma} = (\hankredu^M)^{-1}  \mathcal{H}_{\sigma, \alpha, \beta}^M, \ \ \hat{B}_{\sigma} = (\hankredu^M)^{-1} \mathcal{H}^M_{\alpha,\sigma}, \ \hat{C} =  \mathcal{H}^M_{\beta}
	%\end{equation}
	%for $\sigSet$. Then the deterministic LSS representation $\hat{\mathscr{S}}=(\{\hat{A}_{\sigma},\hat{B}_{\sigma}\}_{\sigma \in \Sigma},\hat{C},\hat{B}, M(\epsilon))$ 
	%is such that $M_{\hat{\mathscr{S}}}=M$, i.e., 
	%$M(\sigma w)=\hat{C}\hat{A}_w\hat{B}_{\sigma}$ for all $w \in \Sigma^{*}$. 
\end{Lemma}

\subsection{Covariance realization algorithm}
Theorem \ref{theo:min}, together with Algorithm \ref{algo:Nice_select_true_deter_basic} and Lemma \ref{basis_red:lemma} suggest that in 
order to compute a sLSS realization of $(\y,\bu,\btheta)$, one could apply Algorithm 
\ref{algo:Nice_select_true_deter_basic} to Hankel-matrices corresponding to $M_{\y,\bu}$.
%The drawback of Algorithm \ref{algo:Nice_select_realization} is that it 
 However, the definition of $M_{\y,\bu}$ uses the 
 covariances of the 
 processes $\y^d$ and $\y^s$, which are not directly
 available.
 %computable from $\y$. 
 Below we show how to
 compute $M_{\y,\bu}$ using only the covariances of $\y$ and $\bu$.
 
 To this end, define the Markov function
    $\incov$ such that $\incov(w)$ is formed
    by the first $n_{u}$ columns of $M_{\y,\bu}(w)$, i.e. $\incov(w)=\Lambda_w^{\y^d,\bu}\varu^{-1}$.
    Then $\incov$ has a realization by a dLSS:
    by Lemma \ref{min:col:lem1}
    %it follows that 
    %if \eqref{eq:aslpv}
    %is a realization of $(\y,\bu,\btheta)$, then
    $(\{\sqrt{p_{\sigma}} A_{\sigma}, \sqrt{p_{\sigma}} B_{\sigma}\}_{\sigma=1}^{\pdim},C,D)$ is a dLSS realization of $\incov$.
 %requires the knowledge of covariances of $\yb^s$. While theoretically $\yb^s$ can be constructed from $\yb$, it is difficult to obtain time series of
%$\yb^s$. Below we present a modification of  Algorithm \ref{algo:Nice_select_realization} to deal with the issue. 
%To this end, we need the following lemma.
\begin{Lemma}
\label{final:alg:lemma}
	%Assume that 
    If there exists a sLSS of $(\yb,\ub,\btheta)$. Then 
    %$\Lambda_{w}^{\y^d,\bu}=\Lambda_{w}^{\y,\bu}$
	\begin{equation}	\label{covseq:lemma:eq1:realization1}
	\begin{split}
    & \Lambda_w^{\y^d,\bu}\!\!=\!\! \Lambda^{\y,\bu}_{w},
    ~
      \Lambda_{\sigma w}^{\y^s,\y^s}\!\!=\!\! \Lambda^{\y,\y}_{\sigma w}-\Lambda^{\y^d,\y^d}_{\sigma w}, \\
      & T^{\y^s,\y^s}_{\sigma,\sigma}\!\!=\!\!T_{\sigma,\sigma}^{\y,\y}\!\!- T_{\sigma,\sigma}^{\y^d,\y^d}
    \end{split}
     %& E[\y^d(t)(\z_{w}^{\bu}(t))^{T}] =E[\y(t)(\z_{\sigma w}^{\bu}(t))^T] \\
	 %& \covseq(\sigma w) = E[\yb(t)(\z^{\yb}_{\sigma w}(t))^T] - E[\yb^d(t)(\z_{\sigma w}^{\yb^d}(t))^T] \\
%	   & E[\zwsigs(t)(\zwsigs(t))^{T}] =  E[\z^{\yb}_{\sigma}(t)(\z^{\yb}_{\sigma}(t))^T]- \\
%      & E[\z^{\yb_d}_{\sigma}(t)(\z^{\yb_d}_{\sigma}(t))^T]
%	\end{split}
        \end{equation}    
        for all $w \in \Sigma^{*}$, $\sigma \in \Sigma$. 
    Moreover,  
    %\[ \incov: \Sigma^{*} \ni w \mapsto %\frac{\Lambda_w^{\y,\bu}\varu^{-1}}%{\sqrt{p_w}}
    %\]
	if  $(\{\tilde{A}_{\sigma}^{d}, \tilde{B}_{\sigma}^{d} \}_{\sigma=1}^{\pdim}, \tilde{C}^{d},\tilde{D}^{d})$ is a minimal dLSS realization of $\incov$, then
	\begin{equation}
	\label{covseq:lemma:eq1:realization}
	\begin{split}
				%&  E[\yb^d(t)(\z_{\sigma w}^{\yb^d}(t))^T]
              & \Lambda^{\y^d,\y^d}_{\sigma w} =\frac{1}{p_{\sigma w}} \tilde{C}^d \tilde{A}_{w}^d(\tilde{A}^d_{\sigma} \tilde{\Psig} (\tilde{C}^d)^{T} + \tilde{B}^d_{\sigma} \varu (\tilde{D}^d)^T ) \\
				%E[\yb_d(t)(\z_{\sigma w}^{\yb_d}(t))^T] \\
				% & E[\yb_d(t)(\z_{\sigma w}^{\yb_d}(t))^T]=
				%& E[\zwsigs(t)(\zwsigs(t))]^{T} =E[\z^{\yb}_{\sigma}(t)(\z^{\yb}_{\sigma}(t))^{T}]-T_{\sigma,\sigma,\mathscr{S}} \\
				%E[\z^{\yb_d}_{\sigma}(t)(\z^{\yb_d}_{\sigma}(t))^T] \\
				& %E[\z^{\yb_d}_{\sigma}(t)(\z^{\yb_d}_{\sigma}(t))^T] 
				%E[\z^{\yb_d}_{\sigma}(t)(\z^{\yb_d}_{\sigma}(t))^T] = 
         T_{\sigma,\sigma}^{\y^d,\y^d}=\frac{1}{\psig} (\tilde{C}^d \tilde{\Psig} (\tilde{C}^d)^{T} + \tilde{D}^d\varu (\tilde{D}^d)^T) 
			\end{split}
	\end{equation}
 where $\tilde{P}_{\sigma}=\lim\limits_{\mathcal{I} \rightarrow \infty} \tilde{P}_{\sigma}^{\mathcal{I}}$,
 and $\tilde{P}_{\sigma}^0=0$ and
 $\tilde{P}_{\sigma}^{\mathcal{I}+1}=p_{\sigma}\sum_{\sigma_1 \in \Sigma} \left ( \frac{1}{p_{\sigma_1}} \tilde{A}_{\sigma_1}^d \tilde{P}_{\sigma_1}^{\mathcal{I}} (\tilde{A}_{\sigma_1}^d)^T+%p_{\sigma_1}
 \tilde{B}_{\sigma_1}^d\varu (\tilde{B}_{\sigma_1}^d)^T\right)$.
\end{Lemma}
The proof of Lemma is presented in \cite{rouphael2024minimal}.
The first part of the lemma
relies on orthogonality of $\y^s(t)$ and $\{\z_w^{\bu}(t)\}_{w \in \Sigma^{*}}$
established in \cite{LCSSArxive}. The second part follows by identifying dLSS realizations of $\incov$
with asLSS realizations of $\y^d$ whose noise is $\bu$, and using the equations in \cite[Lemma 4]{PetreczkyBilinear} for
the covariances of $\y^d$.

%\MM{Mihaly: Existence results to be done}
We can then use Lemma \ref{final:alg:lemma}
to compute the necessary values of $M_{\y,\bu}$ as
follows. 
First, we compute  a minimal dLSS realization of $\incov$ using Algorithm \ref{algo:Nice_select_true_deter_basic}, and then use \eqref{covseq:lemma:eq1:realization1} -- \eqref{covseq:lemma:eq1:realization} to compute those values of $\Lambda_w^{\y^s,\y^s}$,
%$T_{\sigma,\sigma}^{\y^s}$, $\sigma \in \Sigma$, $w \in \Sigma^{+}$, 
and hence of $M_{\yb,\ub}$, which are necessary for applying Algorithm \ref{algo:Nice_select_true_deter_basic} to $M_{\yb,\ub}$.
This idea is formalized in Algorithm~\ref{algo:Nice_select_realization}.
%Let $(\alpha,\beta)$-be a $\nx$-selection and let
 %That is, $\mathcal{L}_{\alpha,\beta}$ contains
 %$l=1,\ldots,\ny \un$
%\vspace{-0.2cm}
\begin{algorithm}[h!]
	\caption{Minimal covariance realization algorithm}
	\label{algo:Nice_select_realization}
	\textbf{Input}: $(\nx,\ny,\un\!\! +\! \ny)$-selection $\left(\alpha, \beta \right)$; 
    $(\bar{n},\ny,\un)$-selection \!\!
    $\left(\bar{\alpha}, \bar{\beta} \right)$;
 %$\{\psig\}_{\sigSet}$,   
     %input covariance $\varu =E[\z^{\ub}_{\sigma}(t)(\z^{\ub}_{\sigma}(t))^{T}]	$, 
      covariances
      $\{\Lambda^{\y,\bu}_w \}_{w \in \mathcal{L}_{\alpha,\beta} \cup \mathcal{L}_{\bar{\alpha},\bar{\beta}}}$,
      $\{\Lambda^{\y,\y}_w \}_{w \in \mathcal{L}_{\alpha,\beta}}$,
      $\Lambda_{\epsilon}^{\y,\bu}$,
      $\{T^{\y,\y}_{\sigma,\sigma}\}_{\sigma \in \Sigma}$; integer $\mathcal{I} >0$.

      %$\{E[\zwsigs(t)(\zwsigs(t))^{T}]\}_{\sigma \in \Sigma}$, 
      %$\{E[\y(t)\bu^T(t)]$,
      %$\{E[\yb(t)(\z^{\bu}_{w}(t))^T],  E[\yb(t)(\z^{\yb}_{w}(t))^T]\}_{w \in \mathcal{L}_{\alpha,\beta}}$
	%$W= \{\wordSet{*}\mid w = v_j\sigma u_i, w=u_i, w = v_i, \forall \sigSet, u_i \in \alpha, v_j \in \beta \}$,
%	maximum number of iterations $\mathcal{I} > 0$.
  \vspace*{.1cm}\hrule\vspace*{.1cm}
	\begin{enumerate}[label=\arabic*., ref=\theenumi{}]
        \item Use 
         %the first equation of
         \eqref{covseq:lemma:eq1:realization1}
        to compute $\{\incov(w)\}_{w \in \mathcal{L}_{\bar{\alpha},\bar{\beta}}}$ and $\incov(\epsilon)$.
        \item Run Algorithm \ref{algo:Nice_select_true_deter_basic} with 
	%	the covariances 
  $M=\incov$ and selection $(\bar{\alpha},\bar{\beta})$ and denote the result by  $\mathscr{S}_d=(\{\tilde{A}_i,\tilde{B}_i\}_{i=1}^{\pdim},\tilde{C},\tilde{D})$.

       \item
       \label{step3:Alg:min}
         Compute
         $\{\Lambda^{\y^s,\y^s}_w, \Lambda_w^{\y^d,\bu}\}_{w \in \mathcal{L}_{\alpha,\beta}}$, 
         %$\{}\}_{w \in \mathcal{L}_{\alpha,\beta}}$,
         using \eqref{covseq:lemma:eq1:realization1}, \eqref{covseq:lemma:eq1:realization}
         and $\{\Lambda^{\y,\y^s}_w\}_{w \in \mathcal{L}_{\alpha,\beta}}$.
     \item
     \label{step4:Alg:min}
         %Use this and %the covariances
         %$\incov(\epsilon)$
         Compute %\hspace{-0.1cm} 
         $M_{\y,\bu}(w)$      
         %$\}_{w \in \mathcal{L}_{\alpha,\beta}}$, %$M_{\y,\bu}(\epsilon)$
         from %M\hspace{-0.1cm}
         $\{\Lambda^{\y^s,\y^s}_w, \Lambda_w^{\y^d,\bu}\}$
         for $w \in \mathcal{L}_{\alpha,\beta}$.
         % for $M=M_{\y,\bu}$.

		\item Run Algorithm \ref{algo:Nice_select_true_deter_basic} with $M=M_{\yb,\ub}$ and selection 
        $(\alpha,\beta)$, and denote by  $\mathscr{S}=(\{\hat{A}_i,\left[\hat{B}_i \ \hat{\Bs}_i\right]\}_{i=1}^{\pdim},\hat{C},\hat{D})$ the result.

		%Denote the result returned by  Algorithm~\ref{algo:Nice_select_true_deter}
		%by $\mathscr{S}=(\{\hat{A}^{d}_{\sigma}, \hat{B}^{d}_{\sigma} \}_{\sigSet}, \hat{C}^{d}, \hat{D}^{d})$.

		\item 
          From $\{T_{\sigma,\sigma}^{\y,\y}\}_{\sigma \in \Sigma}$
           compute 
           $\{T_{\sigma,\sigma}^{\y^s,\y^s}\}_{\sigma \in \Sigma}$         
            using \eqref{covseq:lemma:eq1:realization1} and \eqref{covseq:lemma:eq1:realization}.
%\item 
  Let $\hat{K}^{\mathcal{I}}_{\sigma}$, $\hat{Q}^{\mathcal{I}}_{\sigma}$ as in \eqref{statecov:iter}.
\end{enumerate}
	%\vspace*{. 1cm}
	\hrule\vspace*{.1cm}
	\textbf{Output}: 
	%Estimates $\{\tilde{A}_{\sigma}, \tilde{B}_{\sigma}, \tilde{K}^{M}_{\sigma}, \hat{Q}^{M,N}_{\sigma}, \hat{P}^{M,N}_{\sigma} \}_{\sigSet}, \hat{C}^{N}$
	 Matrices $\{\frac{\hat{A}_{\sigma}}{\sqrt{p_{\sigma}}}, 
  \frac{\hat{B}_{\sigma}}{\sqrt{p_{\sigma}}}, \hat{K}^{\mathcal{I}}_{\sigma}, \hat{Q}^{\mathcal{I}}_{\sigma} \}_{\sigma=1}^{\pdim}, \hat{C}, \hat{D}$.
\end{algorithm}
\vspace{-0.1cm}
%It is clear from Theorem \ref{theo:min}  and  Lemma \ref{basis_red:lemma} that Algorithm \ref{algo:Nice_select_realization} is correct.
%%\begin{color}{red}
\begin{Corollary}
	\label{thm:cra:col1}
        Assume that $\y$ is full rank and there exists a minimal  sLSS realization of $(\yb,\ub,\btheta)$ of dimension $\nx$, %in  innovation form,
        and a minimal dLSS realization of $\Psi_{\y,\bu}$ of 
        dimension $\bar{n}$.
     %\begin{description}
     %{\bf (S1)} 
     If the  
     selections $(\alpha,\beta)$ and $(\bar{\alpha},\bar{\beta})$     
     %$(\nx,\ny,\un+\ny)$-selection $(\alpha,\beta)$ and $(\bar{n},\ny,\un)$-selection
     %$(\bar{\alpha},\bar{\beta})$
	satisfy
    \begin{equation}
    \label{thm:cra:col:eq1}  
     \mathrm{rank} H^{M_{\yb,\ub}}_{\alpha,\beta}\!\!=\!\!\nx, \quad
 \mathrm{rank} H^{\incov}_{\bar{\alpha},\bar{\beta}}\!\!=\!\!\bar{n}.
    \end{equation}
    %Moreover, for any selections 
     %{\bf (S2)} If  $\mathrm{rank}  H^{M_{\yb,\ub}}_{\alpha,\beta}=\nx$ and $\mathrm{rank} H^{\incov}_{\bar{\alpha},\bar{\beta}}\!\!=\!\!\bar{n}$, 
     then Algorithm \ref{algo:Nice_select_realization} returns matrices such that
      %\begin{itemize}
 %     \textbf{\bf (S2.1)}
        the sLSS $\tilde{S}=(\{\hat{A}_{\sigma}/\sqrt{p_{\sigma}},\hat{B}_{\sigma}/\sqrt{p_{\sigma}},\hat{K}_{\sigma},\}_{\sigma=1}^{\pdim},\hat{C},\hat{D},\eb^s)$
is a minimal realization of $(\yb,\ub,\btheta)$  
in innovation form, 
where $\hat{K}_{\sigma}=\lim_{\mathcal{I} \rightarrow \infty} \hat{K}_{\sigma}^{\mathcal{I}}$ and
%the innovation noise covariance is
%and
      %\item{\bf (S2.2)}  
%      $\tilde{S}$ is isomorphic to $\mathcal{S}$.
      %in particular, the deterministic input-output behaviors of $\tilde{S}$ and $\mathcal{S}$ coincide. 
      %\item %
%{\bf (S2.3)}
	%$E[\eb^s(t)(\eb^s(t))^T\p_{\sigma}^2(t)]=\lim_{\mathcal{I} \rightarrow \infty} 
 $T^{\e^s,\e^s}_{\sigma,\sigma}p_{\sigma}=\lim_{\mathcal{I} \rightarrow +\infty} \hat{Q}_{\sigma}^{\mathcal{I}}$, $\sigma \in \Sigma$.
 Moreover, there exist an $(\nx,\ny,\un+\ny)$-selection $(\alpha,\beta)$ and $(\bar{n},\ny,\un)$-selection $(\bar{\alpha},\bar{\beta})$
 which satisfy \eqref{thm:cra:col:eq1}.
 %     \end{itemize}
    % \end{description} 
         %\eqref{lemma:stoch-real1:eq} is an asLSS of $(\yb^s,\p)$, and \eqref{lemma:stoch-real1:eq} holds.
%%	%$(\{\frac{1}{\sqrt{\psig}} \hat{A}_{\sigma}, \hat{\Bs}_{\sigma}, \hat{K}^{\mathcal{I}}_{\sigma}, \hat{Q}^{\mathcal{I}}_{\sigma}, \hat{P}^{\mathcal{I}}_{\sigma} \}_{\sigSet}, \hat{C}^{s})$ are such that 
%%%such that with $\hat{K}_{\sigma}=\lim_{\mathcal{I} \rightarrow \infty} \hat{K}_{\sigma}^{\mathcal{I}}$, $\hat{Q}_{\sigma}=\lim_{\mathcal{I} \rightarrow \infty} \hat{Q}_{\sigma}^{\mathcal{I}}$,
%%%	$\hat{P}_{\sigma}=\lim_{\mathcal{I} \rightarrow \infty} \hat{P}_{\sigma}^{\mathcal{I}}$;
%%%	tuple $(\{\hat{A}_{\sigma}^s, \hat{K}_{\sigma} \}_{\sigSet}, \hat{C}^s, I_{\ny}, \hat{\xb},\eb^s)$ is a stationary
%%%	LSS representation of $\yb^s$ without inputs, and 
%%%	$\hat{P}_{\sigma}\!\!=\!\!E[\hat{\xb}(t)\hat{\xb}^T(t)\p_{\sigma}^2(t)]$, $\sigma \in \Sigma$. 
\end{Corollary}
The proof of Corollary \ref{thm:cra:col1} is presented in \cite{rouphael2024minimal}, it  uses
Lemma \ref{basis_red:lemma}, Theorem \ref{theo:min} and the notion
of sLSSs associated with dLSSs.

That is, Algorithm \ref{algo:Nice_select_realization}
returns a minimal sLSS in innovation form based
on the input and output covariances, 
%solves the realization problem 
%\ref{problem0} 
%and the identification problem
%\ref{problem0:alternative},
%if $(\yb,\ub,\btheta)$ are generated by a minimal sLSS in forward innovation form
%and 
if the selections $(\alpha,\beta)$, $(\bar{\alpha},\bar{\beta})$ 
give Hankel-matrices of a correct rank, and there always exists
such selections.
%The latter can always be assumed without loss of generality, due to Theorem \ref{theo:min} and Lemma \ref{basis_red:lemma}, if there exists a sLSS realization of $(\yb,\ub,\btheta)$. 

 That is, by finding a minimal sLSS of $(\yb,\ub,\btheta)$ in innovation form, Algorithm \ref{algo:Nice_select_realization} finds  a sLSS which has the same deterministic behavior as the true system, i.e.
 the same output response as the true system for
 any inputs, noise  and switching
 signals, if the true system is
 assumed to be minimal in innovation form.   %stochastic inputs, noises and scheduling signals, but for all inputs, noises and scheduling signals}, if the latter is assumed to be minimal sLSS in forward innovation form. 

\section{Identification algorithm}
\label{sect:ident}
We formulate an identification algorithm  based on Algorithm \ref{algo:Nice_select_realization}, 
%and selections,  for  $N$-length observation sequence of outputs, inputs and scheduling signals, as detailed in Algorithm~\ref{algo:Nice_select}. 
%Intuitively, the main idea %Algorithm~\ref{algo:Nice_select} 
 by using empirical  covariances
instead of the true ones
$\Lambda_w^{\y,\bu}$, $\Lambda_w^{\y,\y}$, $T^{\y,\y}_{\sigma,\sigma}$ %and  $E[\z^{\yb}_{\sigma}(t)(\z^{\yb}_{\sigma}(t))^{T} ]$
%using the corresponding empirical covariances computed from observed data 
when applying Algorithm \ref{algo:Nice_select_realization}.
%to the thus estimated covariances. 
%In  Algorithm~\ref{algo:Nice_select}, we assume that $\yb$, $\ub$, $\p$, $\eb$ are jointly ergodic.
%More specifically, 
To this end, we make the following assumptions.
\begin{Assumption}
	\label{assum1}
	%\begin{enumerate}
	%       \item 
	%\label{assum1:as1}
	%	\textbf{(1)}
	%	The process $\yb$ is SII, full rank, and $\yb$ has a realization by stationary stochastic LSS representation
	\textbf{(1)}
	%\label{assum1:as2}
	The $(\nx,\ny,\un+\ny)$-selection $\left(\alpha, \beta \right)$ and the $(\bar{n},\ny,\un)$-selection $(\alphas,\betas)$ satisfy \eqref{thm:cra:col:eq1}, 
     $\y$ is full rank, and
     $\nx$ and $\bar{n}$ satisfy the hypothesis of 
 %that they satisfy the hypothesis of part \textbf{(S2)} of 
 Corollary \ref{thm:cra:col1}.
	%$ \mathrm{rank}% \mathcal{H}^{\incov}_{\bar{\alpha},\bar{\beta}} = n$,  
	%$ \mathrm{rank} \mathcal{H}_{\alpha,\beta}^{M_{\yb,\ub}} = \nx$, where $\nx$ is the dimension of a minimal sLSS  realization of $(\yb,\ub,\btheta)$, 
    %    $\nx \le n$ and $n$ is the dimension of a minimal dLSS realization of $\incov$. 
	%\item 
	
	\textbf{(2)}
	%\label{assum1:as3}
	The process $(\yb, \ub,\btheta)$ is ergodic, and the observed sample paths
	$y: \mathbb{Z} \rightarrow \mathbb{R}^{\ny}$, $u: \mathbb{Z} \rightarrow \mathbb{R}^{\un}$ and $q: \mathbb{Z} \rightarrow \Sigma$ of $\yb$, $\ub$ and $\btheta$ respectively 
    satisfy the following.
	%$\{y(t),u(t),q(t) \}_{t=1}^{N}$ is observed and the following holds. 
     %=\sup_{w \in \mathcal{L}_{\alpha,\beta} \cup \mathcal{L}_{\bar{\alpha},\bar{\beta}} |w|$ be such 
    For all $w,v \in \mathcal{L}_{\alpha,\beta} \cup \mathcal{L}_{\bar{\alpha},\bar{\beta}} \cup \Sigma \cup \{\epsilon\}$  define
    the \emph{empirical covariances}
	%for all $w \in \Sigma^{*}, \sigma \in \Sigma$, 
	%a sample path of  $(\yb, \{\mu_{\sigma} \}_{\sigSet})$  such that 
	%for all $ \wordSet{+}$ the following holds. Define
	%Define $\incov^N: \Sigma^{*} \rightarrow \mathbb{R}^{\ny \times \un}$,
	%$\covseq^N: \Sigma^{+} \rightarrow \mathbb{R}^{\ny \times \ny}$,
	%$T_{\sigma,\sigma}^N$, $\sigma \in \Sigma$ as follows 
	\begin{equation}
    \label{emp:cov}
		%\label{assum1:as3:eq}
		%\begin{split}
			%	\expect{\mathbf{y}(t) (\zwy(t))^{T}} &\!\!=\!\! \lim\limits_{N \rightarrow \infty} 
			 %{\incov^N(w)}= \frac{1}{\sqrt{p_w}} 
	%		\left(
   T^{\y,\mathbf{b}}_{v,w,N}=%{\hat{N}}
   \frac{\sum_{t=N_0}^{N} z^{y}_v(t)(z_w^{b}(t))^{T}}{N-N_0}, 
   %\quad T^{N}_{\sigma,\sigma}=\frac{\sum_{t=N_0}^{\hat{N}} z^{y}_{\sigma}(t)(z_{\sigma}^{y}(t))^{T}}{\hat{N}} 
   %T^{r,r}_{\sigma,\sigma}
   %\\
   %  \right)\varu^{-1},  ~ w \in  \Sigma^{*} \\
			%& ~ \incov^N(\epsilon)=\frac{1}{N} \left(\sum_{t=1}^{N}  y(t)u^T(t)\right)\varu^{-1} \\
	%		& \Lambda^{\yb,\yb,N}_{\sigma w}=\frac{1}{N} \sum_{t=|w|}^{N} y(t)(z_{\sigma w}^{y}(t))^{T}, ~  
			%         \frac{1}{N} \sum_{t=|w|}^{N} \bar{y}(t) (\bar{z}_{w}(t))^{T} \\
			%	\expect{\zvy(t)(\zwy(t))^{T}} &\!\!=\!\! \lim\limits_{N \rightarrow \infty} \frac{1}{N}\!\! \sum_{t=\mathrm{max}(|v|,|w|)}^{N}\!\! \bar{z}_{v}^{}(t) \bar{z}_{w}^{T}(t)
			%\expect{\yb(t)(\zwu(t))^{T} }  &= \lim\limits_{N \rightarrow \infty} 
			%\expect{\yb(t) \ub^{T}(t)} &= \lim\limits_{N \rightarrow \infty} \frac{1}{N} \sum_{t=1}^{N} y(t)({u}(t))^{T} \\
			%\expect{\yb(t)(\z^{\yb}(t))^{T} }  &= \lim\limits_{N \rightarrow \infty} 
		%\end{split}
	\end{equation}
	%Then for all $w \in \Sigma^{*}$, $\sigma \in \Sigma$,
	%\begin{align*}
	%	& \incov(w) =  \lim\limits_{N \rightarrow \infty} \incov^N(w),   ~ E[\z^{\yb}_{\sigma}(t)(\z^{\yb}_{\sigma}(t))^{T}]=\lim\limits_{N \rightarrow \infty} T^{\yb,N}_{\sigma,\sigma} \\
		%\incov(\epsilon) =  \lim\limits_{N \rightarrow \infty} \incov^N(\epsilon),  \\
	%	& E[\yb(t)(\z^{\yb}_{\sigma w}(t))^T] =  \lim\limits_{N \rightarrow \infty} \Lambda^{\yb,N}_{\sigma w}, ~ 
	%\end{align*}
	where $b=\begin{cases} y & \mathbf{b}=\y \\ u & \mathbf{b}=\bu \end{cases}$, and
    for all 
    $\mu_{\sigma}(t)=\chi(q(t)=\sigma)$, $\sigma \in \Sigma$,
    %\begin{cases} 1 & q(t)=\sigma \\ 0 & \mbox{otheriwse} \end{cases}$, and 
	and 
    for all $w=\sigma_1\sigma_2 \ldots \sigma_r \in \Sigma^{+}$, $r > 0$,  
    $\sigma_1,\ldots,\sigma_r \in \Sigma$,
    $b \in \{y,u\}$
	\begin{equation*}
		%\label{assum1:as3:eq}
		\begin{split} 
			& \mu_{w}(t) = \mu_{\sigma_{1}}(t-k+1)\mu_{\sigma_{2}}(t-k+2)\cdots \mu_{\sigma_{r}}(t) \\
			& z^{b}_w(t) = b(t-r) \mu_{w}(t-1)\frac{1}{\sqrt{p_{w}}}, \quad   z^b_{\epsilon}(t)=b(t), \\
			%&  z^{y}_w(t) = y(t-|w|) \mu_{w}%(t-1)\frac{1}{\sqrt{p_{w}}}
		\end{split}
	\end{equation*}
	and $N_0$ is an upper bound on the length of
    words in $\mathcal{L}_{\alpha,\beta} \cup \mathcal{L}_{\bar{\alpha},\bar{\beta}} \cup \Sigma$.
    Then we assume that for all $\mathbf{b} \in \{\y,\bu\}$,
    \begin{equation}
     \label{emp:cov:eq1}
      \Lambda^{\yb,\mathbf{u}}_{w}=\lim_{N \rightarrow \infty} T^{\yb, \mathbf{b}}_{\epsilon,w,N}, \quad 
      T^{\yb,\yb}_{\sigma,\sigma}=\lim_{N \rightarrow \infty} T^{\yb, \yb}_{\sigma,\sigma,N}
    \end{equation}
	%\item $r$- selection $\left(\alpha, \beta \right)$ is such that $ \mathrm{rank} \ \hankred  = n$, where $r$ being the estimate of the state dimension $n$. 
	%\end{enumerate}
\end{Assumption}
The first assumption is not restrictive, such selections always exist, and they can be found by an exhaustive search through all the possible selections.
The assumption \eqref{emp:cov:eq1} says that the observed sample $(y,u,q)$ satisfies the law of large numbers; ergodicity, which is a standard assumption in identification,  implies that almost any sample satisfies \eqref{emp:cov:eq1}. 
%Ergodicity is a standard
%assumption in system identification. 
%\begin{color}{red}
\begin{algorithm}[ht]
	\caption{Identification sLSS}
	\label{algo:Nice_select}
	\textbf{Input}: Data $\{{y}(t),u(t),q(t) \}_{t=1}^N$, $(\nx,\ny,\un+\ny)$-selection $\left(\alpha, \beta \right)$ and $(\bar{n},\ny,\un)$-selection $\left(\alphas, \betas \right)$; %$\{\psig\}_{\sigSet}$, %$\varu$, 
	%$W= \{\wordSet{*}\mid w = v_j\sigma u_i, w=u_i, w = v_i, \forall \sigSet, u_i \in \alpha, v_j \in \beta \}$,
	integer $\mathcal{I} >0$.
	\vspace*{.1cm}\hrule\vspace*{.1cm}
	\begin{enumerate}[label=\arabic*., ref=\theenumi{}]
		%\item \textbf{let} $\mathcal{I} \leftarrow \mathbf{0}_{N}$;
		%\item \textbf{let} $\Clust{i}\leftarrow\emptyset$, $i=1,\ldots,\nmodes$;
		%\item Compute empirical covariances 
        %$T^{\sigma,\sigma,N}^{\y,\y}$,
        %$\sigma \in \Sigma$
		%\begin{equation*}
		%{\incov^N(w)}= \frac{1}{\sqrt{p_w}} \left( \frac{1}{N} \sum_{t=|w|}^{N} y(t)(z_w^{u}(t))^{T}  \right)\varu^{-1}
		%		\end{equation*}
	%	for every $w \in \Sigma^{+}$, such that $w=ivu$ or  $w=i v\sigma u$ or  $w=iv$ or $w=\sigma u$ for some words $v,u \in \Sigma^{*}$, $\sigma \in \Sigma$,
		%$(u,k) \in \alpha$, $(i, v,l) \in \beta$ for some $k=1,\ldots,\ny$, $l=1,\ldots,\un$, 
		%$ \forall w \in \Sigma^{+}$. % ${z}^{u}_{w}(t)$ is as in \eqref{assum1:as3:eq}.
		
		%\item Construct \emph{empirical} Hankel matrices $\mathcal{H}^{N}_{\alpha, \beta}$, $\mathcal{H}^{N}_{\sigma, \alpha, \beta} $,  $\mathcal{H}^{N}_{\alpha,\sigma}$ and $\mathcal{H}^{N}_{\beta}$ based on empirical  covariance ${\incov(w)}^{N}$.
		\item Run Algorithm \ref{algo:Nice_select_realization} with 
        the covariances $\Lambda^{\y,\mathbf{b}}_{w}$, $T^{\y,\y}_{\sigma,\sigma}$ being replaced by the 
        empirical covariances $T^{\y,\mathbf{b}}_{\epsilon,w,N}$, $T^{\y,\y}_{\sigma,\sigma,N}$
        respectively for 
        $\mathbf{b} \in \{\y,\bu\}$ and 
        $w \in \mathcal{L}_{\alpha,\beta} \cup \mathcal{L}_{\bar{\alpha},\bar{\beta}}$, $\sigma \in \Sigma$.
	\end{enumerate}
	%\vspace*{. 1cm}
	\hrule\vspace*{.1cm}
	\textbf{Output}: 
	%Estimates $\{\tilde{A}_{\sigma}, \tilde{B}_{\sigma}, \tilde{K}^{M}_{\sigma}, \hat{Q}^{M,N}_{\sigma}, \hat{P}^{M,N}_{\sigma} \}_{\sigSet}, \hat{C}^{N}$
	The matrices  $\{\tilde{A}_{\sigma}^N, \tilde{B}^N_{\sigma}, \tilde{K}^{N,\mathcal{I}}_{\sigma}, \tilde{Q}^{N,\mathcal{I}}_{\sigma} \}_{\sigma=1}^{\pdim}, \tilde{C}^N, \tilde{D}^N$ 
   returned by Algorithm \ref{algo:Nice_select_realization}.
\end{algorithm}
%\end{color}
\begin{Lemma}[Consistency]
	\label{lm:consistency}
	Under Assumption \ref{assum1}
	the matrices returned by Algorithm \ref{algo:Nice_select} satisfy the following:  
	\begin{align*}
		& \tilde{K}_{\sigma}= \lim_{\mathcal{I} \rightarrow \infty} \lim_{N \rightarrow \infty}\tilde{K}^{N,\mathcal{I}}_{\sigma}, ~~
	p_{\sigma} T^{\e^s,\e^s}_{\sigma, \sigma}=%E[\eb^s(t)(\eb^s(t))^T\p_{\sigma}^2(t)]=
  \lim_{\mathcal{I} \rightarrow \infty} \lim_{N \rightarrow \infty} \tilde{Q}^{N,\mathcal{I}}_{\sigma} 
		\\
		%& \tilde{P}_{\sigma}= \lim_{T \rightarrow \infty} \lim_{N \rightarrow \infty} \tilde{P}^{M,N}_{\sigma} \ \sigSet, \\
		& [\tilde{A}_{\sigma}, \tilde{B}_{\sigma}]=\lim_{N \rightarrow \infty} [\tilde{A}_{\sigma}^N,  \tilde{B}_{\sigma}^{N}], ~  [\tilde{C},\tilde{D}]=\lim_{N \rightarrow \infty} [\tilde{C}^{N}, \tilde{D}^{N}]
        %\tilde{A}_{\sigma}=\lim_{N \rightarrow \infty} %\tilde{A}_{\sigma}^{N},
		%\tilde{D}=\lim_{N \rightarrow \infty} %\tilde{D}^{N}
	\end{align*}
	and $\mathcal{S}=(\{\tilde{A}_{\sigma},\tilde{B}_{\sigma},\tilde{K}_{\sigma},\}_{\sigma=1}^{\pdim},\tilde{C},\tilde{D},\hat{\xb},\eb^s)$
	is a minimal sLSS realization of $(\yb,\ub,\btheta)$ in innovation form.
 %%and the noise covariance satisfies
        %the latter sLSS satisfies
%      for
%	$\sigma \in \Sigma$.
\end{Lemma}
That is, Algorithm \ref{algo:Nice_select}
is a statistically consistent. Moreover, if the true system is a minimal sLSS in innovation form, which, by Theorem \ref{theo:min} we can assume w.l.g, then the sLSS returned by Algorithm \ref{algo:Nice_select} will be 
isomorphic to the true system, as $N \rightarrow \infty$. In particular, in the limit, the identified system has \emph{the same deterministic behavior}
as the true system. That is, while the data used for identification had to be sampled from certain distributions, 
the identified model recreates the behavior of the true one for any switching signal and input. 

Algorithm \ref{algo:Nice_select} 
can be improved as follows.
%by using the extensions below.
\begin{Remark}[Computing empirical covariances]
\label{remark:least_square}
Using \eqref{emp:cov} directly for computing empirical covariances 
$T^{\y,\mathbf{b}}_{\epsilon,w,N}$, $\mathbf{b} \in \{\y,\bu\}$
results in slow convergence. Instead, we propose to compute the empirical covariances
by solving a linear regression problem.  
More precisely, assume
that $\mathcal{L}_{\alpha,\beta} \cup \mathcal{L}_{\bar{\alpha},\bar{\beta}}=\{w_1,\ldots, w_k\}$ and 
$\mathcal{L}_{\alpha,\beta}=\{w_1,\ldots,w_r\}$ for some $k,r>0$. 
%Suppose that we want to find the following covariances $$\frac{1}{N-N_0}\sum_{t=N_0}^{N}\r(t)\zwb(t) \ \text{where} \ N_0=\max_{w\in\mathscr{W}}|w| $$ for some pro
%cesses $\r$ and $\b$ and for some $w \in \mathscr{W}=\{w_1,w_2,\dots,w_k\}$ where $k>0$. 
Let us then define the following matrices 
    \[R=\begin{bmatrix}
        y(N_0)\\ \vdots \\ y(N) \end{bmatrix}, \quad % \text{and} \ 
        \Phi_{\mathbf{b}}=\begin{bmatrix}
            \mathbf{z}^{y}_{w_1}(N_0) &\cdots &\mathbf{z}^{\mathbf{b}}_{w_{l(\mathbf{b})}}(N_0)\\
            \vdots & &\vdots\\
            \mathbf{z}^{\mathbf{b}}_{w_1}(N) &\cdots &\mathbf{z}^{\mathbf{b}}_{w_{l(\mathbf{b})}}(N)
        \end{bmatrix}\]
    where for $\mathbf{b}=\y$, $l(\mathbf{b})=r$ and $b=y$, and
    for $\mathbf{b}=\bu$, $l(\mathbf{b})=k$ and $b=u$.
    %Then $\frac{1}{N-N_0}\Phi^TR$ contains all desired the empirical covariances. However, 
     By the well-known formula, 
     $\Hat{\Theta}^{\mathbf{b}}=(\Phi^T_{\mathbf{b}}\Phi_{\mathbf{b}})^{-1}\Phi_{\mathbf{b}}^TR$ is the 
     least squares solution of the equation $R=\Phi_{\mathbf{b}}\Hat{\Theta}^{\mathbf{b}}$. 
     Hence, $\frac{1}{N-N_0}(\Phi^T_{\mathbf{b}}\Phi_{\mathbf{b}})\Hat{\Theta}^{\mathbf{b}}=\frac{1}{N-N_0} \Phi^T_{\mathbf{b}}R$, and the latter 
     contains all covariances $T_{\epsilon,w,N}^{\y,\mathbf{b}}$, 
     $\mathbf{b} \in \{\bu,\y\}$
     required for Algorithm \ref{algo:Nice_select}.
     In turn, $\Hat{\Theta}^{\mathbf{b}}$ can be computed using standard
     linear regression tools. 
     %by taking 
     %$\mathbf{b}=\bu$ and $\mathbf{b}=\y$.
     %Here $w$ ranges through $\mathcal{L}_{\alpha,\beta} \cup \mathcal{L}_{\bar{\alpha},\bar{\beta}}$ if $\mathbf{b}=\bu$, and 
     %$w$ ranges through $\mathcal{L}_{\alpha,\beta}$ if $\mathbf{b}=\y$.     
     %This method can be used in order to calculate the empirical covariances $T^{\y,\mathbf{b}}_{\epsilon,w,N}$ in Algorithm \ref{algo:Nice_select}.
\end{Remark}

\begin{Remark}[Refinement using gradient descend]
\label{rem:gb}
    In order to enhance the performance of Algorithm \ref{algo:Nice_select}, 
    inspired by \cite{cox2018towards}, we propose to
    %we can use the extension of Gradient-Based (GB) search method to LPV systems \cite{4459815,cox2018towards}.
    %This concept was used to find Linear Parameter-Varying LPV models in \cite{cox2018towards}. 
    %Inspired by the author of the later thesis, we applied the same refinement step to Algorithm \ref{algo:Nice_select}. 
    %To this end, we 
    use the sLSS 
    %$\mathcal{S}=(\{\tilde{A}_{\sigma}^N, \tilde{B}^N_{\sigma},\tilde{K}^{N,\mathcal{I}}_{\sigma}\}_{\sigma=0}^{\pdim},\tilde{C}^N, \tilde{D}^N,)$ 
    returned by Algorithm \ref{algo:Nice_select} as initial value for the Gradient-Based (GB) algorithm of \cite[Algorithm 7.1]{cox2018towards}. 
    %Note that the main difference between Algorithm \ref{algo:Nice_select} and the approach of 
    Note that \cite{cox2018towards} uses zero matrices as initial values for the noise gain matrices,
    %$\tilde{K}_{\sigma}$, 
    whereas Algorithm \ref{algo:Nice_select} returns a non-zero estimate of those matrices.
    %and takes a matrix with entries equal to zeros as the initial value of $\tilde{K}^{N,\mathcal{I}}_{\sigma}$ for the Gradient search. 
\end{Remark}

\section{Numerical example}\label{sec:example}
In this section, we illustrate Algorithm \ref{algo:Nice_select} on a numerical example.
%to test the effectiveness of our algorithm. 
All computations are carried out on an i7 2.11-GHz Intel core processor with 32GB of RAM running MATLAB R2023b.
%The proposed algorithm is compared to the  \emph{Predictor Based Subspace IDentification} (PBSID) algorithm reported in~\cite{van2009subspace}. The underlying idea of PBSID methods is to first estimate the state sequence by finding appropriate projection and matrix decomposition of data matrices and then to compute LPV model matrices using the estimated state sequence via convex optimization. 
%We remark that the PBSID algorithm ~\cite{van2009subspace}, requires the observability assumption of the first local model, i.e., the extended observability matrix $\mathcal{O} = [C^{T} (C\tilde{A}_{1})^{T} \ldots (C\tilde{A}^{l-1}_{1})^{T}]^{T}$ (where $\tilde{A}_{1}= A_{1}- K_{1}C$) has to have full column rank for a certain window length $l$, along with other assumptions on scheduling variable sequence. 
%
%The realization algorithm proposed in this contribution relaxes these assumptions. However, it relies on choosing the correct \emph{selections} $\left( \alpha, \beta \right)$. Specifically, in the following example, we show that the observability assumption is not satisfied by the true data generating system, yet, the proposed algorithm is able to estimate the model matrices such that the  estimated output  matches closely with the observed outputs. 
For data generation we use the following randomly generatated  sLSS \eqref{eq:aslpv} with $\Sigma=\{1,2\}$ %and matrices
%is used for data generation with following matrices:
\[\begin{split}
	%\begin{equation*}
	& A_{1} = \begin{bmatrix}
        0.1039    &0.0255    &0.5598\\
        0.4338    &0.0067    &0.0078\\
        0.3435    &0.0412    &0.0776
		\end{bmatrix}, 
        B_{1} = \begin{bmatrix} 1.6143\\
                                5.9383\\
                                7.3671
                \end{bmatrix}, \\
	%\end{equation*}
	%\begin{equation*}
 %\begin{split}
	&A_{2} = \begin{bmatrix}
        0.1834    &0.2456    &0.0511\\
        0.0572    &0.2445    &0.0642\\
        0.1395    &0.6413    &0.5598 
	\end{bmatrix},
 	B_{2} = \begin{bmatrix} 6.0624\\
                                4.9800\\
                                3.1372
                \end{bmatrix},\\
	%\end{equation*}
	%\begin{equation*}
 %\end{split}
 %\]
 %\[
 %\begin{split}
	&K_{1}\!\! = \!\! 
	\left[0.4942 \ 0.2827 \ 0.8098 \right]^{T}, \\ 
 & K_{2}\!\! =\!\! \left[0.6215 \ 0.1561 \ 0.7780  \right]^{T} \\
	%\end{equation*}
	%\begin{equation*}
        &C = \begin{bmatrix}
		0.1144    &0.7623    &0.0020
	\end{bmatrix},
	D =1,
        F=1.
	%\end{equation*}
\end{split}
\]
%\begin{equation*}
%C = \begin{bmatrix}
%1 & 0 & 0
%\end{bmatrix}, 
%\end{equation*}
%which corresponds to state-dimension $\nx\!\!=\!\!3$, output dimension $\ny\!\!=\!\!1$, and the set of discrete modes %$\pdim\!\!=\!\!2$ with 
%$\Sigma \!=\! \{1,2\}$. 
%Note that, the system corresponding to first local model $\tilde{A}_{1} = A_{1}- K_{1}C$ is \emph{not observable}, i.e.,  $ \mathrm{rank}([C^{T} (C\tilde{A}_{1})^{T} \ldots (C\tilde{A}^{l-1}_{1})^{T}]^{T}) =2< \nx$, which is a particular assumption required in subspace based approaches \cite{vansubspace}. 
 %The operator$\mathrm{var}(\cdot)$ denotes the variance of its argument.

Training and validation data of length $N$ and $N_\mathrm{val}$, respectively, are generated using the data generator sLSS with white-noise input process $\ub$ with uniform distribution $\mathcal{U}(-1,1)$, 
an i.i.d. process $\btheta$ %is an i.i.d. process
sucht that 
$\mathcal{P}(\btheta(t)=q)=0.5$, $q=1,2$, and 
%
%an independent scheduling signal process $\p = [\p_{1} \  \ \p_{2}] $ such that $$\p_{\sigma}(t)=\chi(\btheta(t)=\sigma)=
%  \begin{cases}
%  1 & \text{if $\btheta(t)=\sigma,$}\\ 
%  0 & \text{otherwise.}
%  \end{cases}$$ for all $\sigma \in \Sigma$, $t \in \mathbb{Z}$, where $\btheta$ is an i.i.d. process with uniform distribution $\mathcal{U}(0,1)$ taking values in $\Sigma$. 
%This corresponds to the parameter values $\{\psig \}_{\sigma \in \{1,2\}}$  to be $p_{1} \!=\! p_{2} \!=\! E[\mu^{2}_{1}(t)] \!\!= \!\!0.5$. 
%For generating the training data, we
%The standard deviation of the 
a Gaussian white noise $\mathbf{v}$ with variance $\sigma^2_{\mathbf{v}}$.
The corresponding 
%corrupting the training output is $2.25$, i.e., $\eb \sim 1.5 \ \mathcal{N}(0,1)$. 
%This corresponds to the 
\emph{Signal-to-Noise Ratio} $\mathrm{SNR}$
are shown in Table \ref{Tab:BFR_eg1}. We run the version of Algorithm \ref{algo:Nice_select} using the Least-Square method from Remark \ref{remark:least_square}, with %the following $3$-selections 
$(\alpha, \beta) = (\alphas, \betas)$, 
%with $n=3$, 
%\begin{align*}
	$$\alpha\!\!=\!\! \{ (11,1), (1,1) , (\epsilon,1)\}, ~ ~ \beta \!\! =\!\! \{(2,\epsilon,1), (1,2,1),(1,1,1)\},$$
	%\alphas &= \{ (\epsilon,1), (1,1) , (21,1)\}, \ \betas = \{ (\epsilon,1), (2,1),(21,1)\}, 
%	\alphas &= \{ (1,1), (\epsilon,1) , (21,1)\}, \ \betas = \{ (2,\epsilon,1), (1,1,1),(1,2,1)\}, 
%\end{align*}
%which are used to choose corresponding entries of the Hankel matrices. 
%The mean time taken to run the algorithm is $2.34$ seconds without refinement and $7$ minutes with the additional Gradient-Based (GB) refinement step.
and the GB refinement step from Remark \ref{rem:gb}.
For the latter we
%Note that, in order to implement the GB refinement, 
used the LPVcore toolbox  \cite{DENBOEF2021385}.

For validation, the true output minus the measurement noise $Fv(t)$, denoted by $y(t)$, was compared with the output $\hat{y}(t)$ predicted by the estimated model in innovation form using past inputs and outputs, as explained in Section \ref{sec:min}.
The  quality of the match is evaluated
%on a observation noise-free validation data of length $N_\mathrm{val}$ 
via  \emph{Best Fit Rate} %(BFR)
criterion %defined
%for each output channel $y_{i}$, $i\!=\!1,\ldots,\ny$, as
%\begin{subequations} \label{BFR}
%\begin{align*}
	\( \mathrm{BFR}\!\! = \!\! \max\left\{1\!\! - \!\! \sqrt{\frac{\sum_{t=1}^{N_\mathrm{val}}\left(y(t)-\hat{y}(t)\right)^2}{\sum_{t=1}^{N_\mathrm{val}}\left(y(t)-\bar{y}\right)^2}},0\right\} \times 100 \%   \), where
% and  \emph{Variance Accounted For} %(VAF) 
% \(
%	\mathrm{VAF} =  \max\left\{1- \frac{\mathrm{var}(y(t)- \hat{y}(t))}{\mathrm{var}(y)},0\right\} \times 100 \% \),
%\end{align*}
%\begin{align*}
%$\mathrm{BFR}_{y_{i}} =\max\left\{1-\sqrt{\frac{\sum_{t=1}^{N_\mathrm{val}}\left(y_{i}(t)-\hat{y}_{i}(t)\right)^2}{\sum_{t=1}^{N_\mathrm{val}}\left(y_{i}(t)-\bar{y_{i}}\right)^2}},0\right\} \times 100 \%$,  
%$\mathrm{VAF}_{y_{i}} =  \max\left\{1- \frac{\mathrm{var}(y_{i}- \hat{y}_{i})}{\mathrm{var}(y_{i})},0\right\} \times 100 \%$, 
%\end{align*}
%where $\hat{y}$ denotes the simulated one-step ahead model output  and 
$\bar{y}$ denotes the sample mean of the output in the validation set.
%10\log{\frac{\sum_{t=1}^N\left(y(t)-e(t)\right)^2}{\sum_{t=1}^Ne^{2}(t)}}= 
%6.1 \ \mathrm{dB}.$ 

%\footnote{All computations are carried out on an i5 1.8-GHz Intel core processor with 8 GB of RAM running MATLAB R2018a.} 
%\begin{table}
%	\caption{Best Fit Rate (BFR) and Variance Accounted For (VAF) on a noise-free validation data }\label{Tab:BFR_eg1}
%	\begin{center} \vspace{-0.3cm}
%		\begin{tabular}{|c||c|c|c|}  
%			\hline
%			& Algorithm \ref{algo:Nice_select} & PBSID  & PBSID \\ 
%			\hline
%			BFR & 93.56 \%  & 1.58 \% & 0 \%\\
%			\hline
%			VAF & 99.58 \% & 3.16 \% & 0 \% \\
%			\hline
%		\end{tabular} %\vspace{-0.6cm}
%	\end{center}
%\end{table}
\begin{table}[!h]
\vspace{-0.25cm}
	\caption{BFR on a noise-free validation data ($N_\mathrm{val}=500$)}
 \label{Tab:BFR_eg1}
 \centering
	\begin{center} \vspace{-0.5cm}
		\begin{tabular}{|c||c|c||c|c|}  
			\hline
				\multirow{3}{*}{Method} & \multicolumn{2}{c||}{$N=5000$} & \multicolumn{2}{c|}{$N=10000$}\\
    \cline{2-5}
   &\multicolumn{4}{c|}{ $\mathrm{SNR}=6.1 \ \mathrm{dB}$, $\sigma_{\mathbf{v}}=1.5$}\\
   \cline{2-5}
    & BFR[\%] & time[s] & BFR[\%] & time[s]\\
   \hline
   GB + zero noise gain & 83.13 & 197 & 82.42 & 361\\ \hline
   Algorithm \ref{algo:Nice_select}	& 85.67 & 2.60 & 90.86 & 2.03\\ \hline
   Algorithm \ref{algo:Nice_select} + GB & 90.5 & 256 & 91.71 & 323\\
    \hline
    \multirow{2}{*}{} &\multicolumn{4}{c|}{$\mathrm{SNR}=0.5 \ \mathrm{dB}$, $\sigma_{\mathbf{v}}=3$}\\
    \hline
   GB + zero noise gain  & 63.25 & 188 & 66.11 & 316\\ \hline
   Algorithm \ref{algo:Nice_select} & 76.97 & 1.94 & 87.29 & 1.90 \\ \hline
   Algorithm \ref{algo:Nice_select} + GB & 84.17 & 236 & 88.03 & 376\\
    \hline
    \end{tabular} \vspace{-0.6cm}
    \end{center}
\end{table}
%The validation result using one-step ahead predicted outputs $\hat{y}$  
The results are reported in Table~\ref{Tab:BFR_eg1}.
%In \cite{cox2018towards}, they used subspace identification to estimate the deterministic system, i.e., deterministic component of the system decomposition \cite{LCSSArxive}. 
%Due to the fact that they do not estimate the noise matrix $K$, they initiate the gradient search with their deterministic estimation and a noise matrix with entries equal to zero. 
%They argue that the behavior of the generated system is close to the behavior of the real system since the behavior of the initial parameters for the gradient search is already close enough.
%Table \ref{Tab:BFR_eg1} shows that 
We compare Algorithm \ref{algo:Nice_select} with using
only the GB search algorithm from \cite{cox2018towards}
where the initial value of the noise gain is zero, while the other matrices are initialized by the outcome of Algorithm \ref{algo:Nice_select}. 
This is then equivalent with combining the CRA algorithm
of \cite{CoxLPVSS} with GB.
Algorithm \ref{algo:Nice_select} appears to perform better than using only GB search, especially for noisy data and a larger number of data points.  Moreover, Algorithm \ref{algo:Nice_select} is faster than GB search.
If combined with
GB search, the performance of Algorithm \ref{algo:Nice_select} improves further, but at the expense of computational time. 
%especially
%In addition, the estimation of $K$ enhances the gradient search.
%The  added GB step helps in reducing the amount of the training data at the expense of time.

%The results  show a good match between estimated model output w.r.t. true system output. 
%We compare our algorithm with the PBSID approach reported in \cite{van2009subspace} and implemented using PBSID MATLAB toolbox. For PBSID, we consider past and future window length $8$ with zero external input signal as well as with a white Gaussian noise external input $u \!\sim\! \mathcal{N}(0,1)$ with input matrices $B_{1}= B_{2} = \left[0 \ 0 \ 1 \right]^{T}$.

\section{Conclusion}
In this paper, we a characterization of minimality and uniquenss of LSSs in innovation form, and 
we presented a realization algorithm for computing minimal LSS in innovation form. Using this realization algorithm, we formulated a system identification algorithm proven to be statistically consistent. The latter algorithm  
was evaluated on a numerical example and demonstrated promising performance, both in terms of runtime complexity and estimation accuracy. 
%with work of \cite{cox2018towards}, and we showed an improvement by combining the proposed identification algorithm and the Gradient-Based search. This refinement step helps in reducing the amount of the data points required in order to obtain a precise estimation of the original system.

%\addtolength{\textheight}{-3cm}
\bibliographystyle{plain}
\bibliography{Bib.bib}
\appendix
\subsection{Proof of Lemma \ref{min:col:lem1}}\label{app:proof_lemma2}
In order to present the proof we need some preliminaries definitions, notations and lemmas.
\begin{Notation}[autonomous sLPSS]
As discuss in Notation  \ref{Notation:sLPSS}, a sLPSS is noted as $\mathcal{S}=(\{A_{\sigma},B_{\sigma},K_{\sigma}\}_{\sigma=1}^{\pdim},C,D,F,\vb)$.
However, if $\mathcal{S}$ depend only on the noise, i.e., $B_{\sigma}=0$ and $D = 0$ for all $\sigma \in \Sigma$, then the sLSS is called autonomous (with no inputs) and abbreviated as asLSS. 
We will note an asLPSS as $\mathcal{S}=(\{A_{\sigma},K_{\sigma}\}_{\sigma=1}^{\pdim},C,F)$ and which has the following form:
    \begin{equation}
	\label{eq:aslss}
            \mathcal{S}_{as} \left\{     \begin{aligned}
                 &  \xb(t+1)=A_{\btheta(t)}\xb(t)+K_{\btheta(t)}\vb(t) \\
                 & \y(t)=C\xb(t)+F\vb(t)
		\end{aligned}\right.
    \end{equation}
    where 
	$A_{\sigma} \!\! \in \!\! \mathbb{R}^{\nx \times \nx}, K_{\sigma} \! \in\!  \mathbb{R}^{\nx \times \nn}$,
    $\sigma \!\in\! \Sigma=\{1,\ldots,\pdim\}$,
	$C\! \in \! \mathbb{R}^{\ny \times \nx}$, $F \! \in \! \mathbb{R}^{\ny \times \nn}$, and $\xb$, $\vb$, $\y$, $\btheta$ 
    are the stochastic state, noise, output 
    and switching processes. The process $\btheta$ takes values in the set of discrete states $\Sigma$.
    We also say that $\mathcal{S}_{as}$ is a realization of $(\yb, \btheta)$.
\end{Notation}
%Recall from Notation \ref{hilbert:notation} the definition of the Hilbert-space $\mathcal{H}_1$ of zero mean square integrable random variables. 
%Let us denote by $\mathcal{H}_{t,+}^{\ub}$, the closed subspace of $\mathcal{H}_1$ generated by the components of  $\{ \z^{\ub}_{w}(t) \}_{w \in \Sigma^{+}} \cup \{\ub(t)\}$.

Assume that $\mathcal{S}$ of the form \eqref{eq:aslpv}  is a sLSS of $(\yb,\ub,\btheta)$\footnote{By Assumption of Lemma \ref{decomp:lemma:inv}, such a LSS representation exists.}.
Recall from Notation \ref{hilbert:notation} the Hilber-space $\mathcal{H}_1$. 
% denotes the Hilbert-space of zero mean square integrable random variables.
Denote by $\mathcal{H}_{t,+}^{\vb}$ and $\mathcal{H}_t^{\vb}$  the closed-subspace of $\mathcal{H}_1$ 
generated by the components of $\{\z^{\vb}_w\}_{w \in \Sigma^{+}} \cup \{\vb(t)\}$ and $\{\z^{\vb}_w\}_{w \in \Sigma^{+}}$ respectively.
\begin{Lemma}
	\label{decomp:lemma:pf1}
	%With the assumptions and notations  of Lemma \ref{decomp:lemma}, 
  For all $\sigSet$, $\vb(t)$
 is orthogonal to $\mathcal{H}_{t,+}^{\ub}$. 
\end{Lemma}

\begin{pf}[Proof of Lemma \ref{decomp:lemma:pf1}]
        Define
	$\mathbf{r}(t):=\begin{bmatrix} \vb^T(t) & \ub^T(t) \end{bmatrix}^T$. By the definition of a sLSS, $\mathbf{r}$  is a ZMWSII and it is a 
        white noise process w.r.t. $\p$. Moreover, $\vb(t)$ is the upper $\nn$ block of $\r(t)$.
	%First, we show that  $E[\vb(t)(\ub(t+1))^T \p_{\sigma}(t)]=0$.  
        Since $\mathbf{r}$ is a white noise w.r.t. $\p$, 
	$E[\mathbf{r}(t)(\z^{\mathbf{r}}_{w}(t))^T]=0$, and  $\frac{1}{\sqrt{p_{\sigma}}} E[\vb(t)(\z_{w}^{\ub}(t))^T]$ is the lower-left block of 
        that latter matrix, and hence it is also zero.
        From  Definition \ref{defn:LPV_SSA_wo_u} it follows that
  %      \st{$E[\vb(t)(\ub(t))^T\p_1(t)]=$}
        $E[\vb(t)(\ub(t))^T]=0$.
\end{pf}

\begin{Lemma}
	\label{decomp:lemma:pf2}
	For any $w \in \Sigma^{+}$, the components of $\z_{w}^{\vb}(t)$ are orthogonal to
%	\st{$\mathcal{H}_{t+k,+}^{\ub}$, $k \ge 0$} 
	$\mathcal{H}_{t,+}^{\ub}$.
\end{Lemma}

\begin{pf}[Proof of Lemma \ref{decomp:lemma:pf2}]	Since $\mathbf{r}(t):=\begin{bmatrix} \vb(t)^T & \ub(t)^T \end{bmatrix}^T$ is a ZMWSII process, from \cite[Lemma 7]{PetreczkyBilinear} {and because $\r(t)$ is a white noise, it follows that
	$E[\z_{w}^{\mathbf{r}}(t)(\z^{\mathbf{r}}_{v}(t))^T]=E[\z_{\sigma}^{\mathbf{r}}(t)(\z^{\mathbf{r}}_{\sigma}(t))^T]$ if $v=w=\sigma s$ for all $v,w,s \in \Sigma^{+}$, $\sigma \in \Sigma$, and 
	$E[\z_{w}^{\mathbf{r}}(t)(\z^{\mathbf{r}}_{v}(t))^T]=0$ for any other case.}
	
	Since $E[\z_{w}^{\vb}(t)(\z^{\ub}_{v}(t))^T]$ is the upper right block of $E[\z_{w}^{\mathbf{r}}(t)(\z^{\mathbf{r}}_{v}(t))^T]$, it follows that 
	$E[\z_{w}^{\vb}(t)(\z^{\ub}_{v}(t))^T]=0$ if $v \ne w$ and  $E[\z_{w}^{\vb}(t)(\z^{\ub}_{w}(t))^T]=E[\z_{\sigma}^{\vb}(t)(\z^{\ub}_{\sigma}(t))^T]=\frac{1}{p_{\sigma}} E[\ub(t-1)\vb(t-1)\p_{\sigma}^2(t-1)]$, where $\sigma$ is the first letter of $w$, and from  Definition \ref{defn:LPV_SSA_wo_u}, it follows that the latter expectation is zero. That is, $E[\z_{w}^{\vb}(t)(\z^{\ub}_{v}(t))^T]=0$ for all $v \in \Sigma^{+}$. 
	
	Finally,  as $\mathbf{r}(t)$ is a 
	white noise process w.r.t. $\p$, it follows that  $E[\z_{w}^{\mathbf{r}}(t)(\mathbf{r}(t))^T]=0$, and since  $E[\z_{w}^{\vb}(t)(\ub(t))^T]$ is the upper right block of $E[\z_{w}^{\mathbf{r}}(t)(\mathbf{r}(t))^T]$, it then follows that $E[\z_{w}^{\vb}(t)(\ub(t))^T]=0$. 
	Since $\z_w^{\vb}(t)$ is orthogonal to the components of the random variables which generate $\mathcal{H}_{t,+}^{\ub}$, the statement of the lemma follows for $k=0$.
	
%\st{Consider now the case $k > 0$. As $\p_1=1$ and $p_1=1$ it follows that $\z_{w}^{\vb}(t)=\z_{w\underbrace{1\cdots1}_{k}}^{\vb}(t+k)$ (where $\underbrace{1\cdots1}_{k}$ denotes $k$-lenght word of $1$s), and $\z_{w\underbrace{1\cdots1}_{k}}^{\vb}(s)$ is orthogonal to $\mathcal{H}_{s,+}^{\ub}$, according to the case $k=0$. By taking $s=t+k$ for $k > 0$, the statement of the lemma follows.}
\end{pf}

Let us denote by $\mathcal{H}_{t}^{\ub}$, the closed subspace generated by the components of  $\{ \z^{\ub}_{w}(t) \}_{w \in \Sigma^{+}}$.It is clear that $\mathcal{H}_{t}^{\ub} \subseteq \mathcal{H}_{t,+}^{\ub}$.
\begin{Lemma}
	\label{decomp:lemma:pf3}
	% For any $w \in \Sigma^{+}$, the components of $\z_{w}^{\vb}(t)$ are orthogonal to $\mathcal{H}_{t,+}^{\ub}$.
	The components of $\xb^d(t)$ belong to $\mathcal{H}_{t}^{\ub}$ and 
	\begin{equation}
		\label{decomp:lemma:pf3:eq2}
		\xb^d(t) = \sum_{w \in \Sigma^{*}, \sigma \in \Sigma} \sqrt{p_{\sigma w}} A_w B_{\sigma}\z^{\ub}_{\sigma w}(t),
	\end{equation}
	where the right-hand side of converges in the mean square sense. 
\end{Lemma}
\begin{pf}[Proof of Lemma \ref{decomp:lemma:pf3}]
	It is clear from the definition that the components of $\xb^d(t)$ belong to $\mathcal{H}_{t,+}^{\ub}$. Since 
	$\mathbf{x}(t)\! \! = \! \! \sum_{w \in \Sigma^{*}, \sigma \in \Sigma} \sqrt{p_{\sigma w}} A_w \left(K_{\sigma} \z^{\vb}_{\sigma w}(t)  + B_{\sigma}\z^{\ub}_{\sigma w}(t)\right)$
	and the fact that the map $z \mapsto E_l[z \mid M]$ (where $z \in \mathcal{H}_1$) is a continuous linear operator for any closed subspace $M$, it follows that
	\begin{equation}
	\label{decomp:lemma:pf3:eq1}
	\begin{split}
		\xb^d(t)=  \sum_{w \in \Sigma^{*}, \sigma \in \Sigma} &\sqrt{p_{\sigma w}} A_w \left(K_{\sigma} E_l[\z^{\vb}_{\sigma w}(t) \mid H_{t,+}^{\ub}] \right. \\ 
		& \left. + B_{\sigma}E_l[\z^{\ub}_{\sigma w}(t) \mid H_{t,+}^{\ub}]\right).
	\end{split}
	\end{equation}
	From Lemma \ref{decomp:lemma:pf2} it follows that,  $E_l[\z^{\vb}_{\sigma w}(t) \mid H_{t,+}^{\ub}]=0$, and since the components of $\z^{\ub}_{\sigma w}(t)$ belong to $\mathcal{H}_{t,+}^{\ub}$, it follows that $E_l[\z^{\ub}_{\sigma w}(t) \mid H_{t,+}^{\ub}]=\z^{\ub}_{\sigma w}(t)$, hence \eqref{decomp:lemma:pf3:eq1} implies \eqref{decomp:lemma:pf3:eq2}.
	%\begin{equation}
	%	\label{decomp:lemma:pf3:eq1}
	%	\begin{split}
	%		\xb^d(t) %&=\sum_{w \in \Sigma^{*}, \sigma \in \Sigma} \sqrt{p_{\sigma w}} A_w B_{\sigma}E_l[\z^{\ub}_{\sigma w}(t) \mid H_{t,+}^{\ub}] \\
	%		& = \sum_{w \in \Sigma^{*}, \sigma \in \Sigma} \sqrt{p_{\sigma w}} A_w B_{\sigma}\z^{\ub}_{\sigma w}(t).
	%	\end{split}
	%\end{equation}
	Since the components of $\z^{\ub}_{\sigma w}(t)$ belong to $\mathcal{H}_{t}^{\ub}$, it follows that the components of the right-hand side of 
	\eqref{decomp:lemma:pf3:eq2} belongs to $\mathcal{H}_{t}^{\ub}$ and hence the components of $\xb^d(t)$ belong to $\mathcal{H}_{t}^{\ub}$.
	Note that, the convergence of the right-hand side  of \eqref{decomp:lemma:pf3:eq2} in the mean square sense follows from the convergence of the series $\sum_{w \in \Sigma^{*}, \sigma \in \Sigma} \sqrt{p_{\sigma w}} A_w \left(K_{\sigma} \z^{\vb}_{\sigma w}(t)  + B_{\sigma}\z^{\ub}_{\sigma w}(t)\right)$.
\end{pf}

\begin{Lemma}
	\label{decomp:lemma:pf4}
	The components of $\xb^s(t)$ belong to $\mathcal{H}_{t}^{\vb}$, they are orthogonal to $\mathcal{H}_{t,+}^{\ub}$, and 
	\begin{equation}
		\label{decomp:lemma:pf4:eq2}
		\xb^s(t) = \sum_{w \in \Sigma^{*}, \sigma \in \Sigma} \sqrt{p_{\sigma w}} A_w K_{\sigma}\z^{\vb}_{\sigma w}(t),
	\end{equation}
	where the right-hand side converges in the mean-square sense.
\end{Lemma}
\begin{pf}[Proof of Lemma \ref{decomp:lemma:pf4}] From \eqref{decomp:lemma:pf3:eq2}, $\xb^s(t)=\xb(t)-\xb^d(t)$
	and 
	$\mathbf{x}(t)\! \! = \! \! \sum_{w \in \Sigma^{*}, \sigma \in \Sigma} \sqrt{p_{\sigma w}} A_w \left(K_{\sigma} \z^{\vb}_{\sigma w}(t)  + B_{\sigma}\z^{\ub}_{\sigma w}(t)\right)$, it follows that \eqref{decomp:lemma:pf4:eq2} holds and that its right-hand side converges in the mean square sense. From Lemma \ref{decomp:lemma:pf2}, it follows that for any $w \in \Sigma^{+}$, the components of $\z_{w}^{\vb}(t)$ are orthogonal to $\mathcal{H}_{t+k,+}^{\ub}$, hence
	all the summands of the infinite series of \eqref{decomp:lemma:pf4:eq2} are orthogonal to $\mathcal{H}_{t,+}^{\ub}$. 
\end{pf}

\begin{Lemma}[Decomposition of $\yb$]
	\label{decomp:lemma}
	Assume that there exists a sLSS of $(\yb,\ub,\btheta)$ of the form
	\eqref{eq:aslpv} and that $(\yb,\ub,\btheta)$ satisfy Assumption \ref{asm:main}.
	Define \begin{align}
		&\xb^d(t)=E_l[\xb(t)  \mid \{ \z^{\ub}_{w}(t) \}_{w \in \Sigma^{+}} \cup \{\ub(t)\}]
		\label{decomp:lemma:state:eq1} \\
		& \xb^s(t)=\xb(t)-\xb^d(t)
		\label{decomp:lemma:state:eq2} 
	\end{align}
	Then 	\begin{eqnarray}
	  & \mathcal{S}_d \left\{ \begin{aligned} 
             &{\xb}^{d}(t+1) = \sum_{i=1}^{\pdim} (A_{i} {\xb}^{d}(t) \!+\! B_{i}  {\ub}(t))\p_i(t),  \\ 
	  &	{{\yb}}^{d}(t) = C {\xb}^{d}(t) \!+\! D  {\ub}(t),
          \end{aligned}\right.
	\label{eqn:LPV_SSA_deter} \\
		& \mathcal{S}_s \left\{
		\begin{aligned}
                & {\xb}^{s}(t+1) = \sum_{i=1}^{\pdim} (A_{i} {\xb}^{s}(t) \!+\! K_{i} {\vb}(t))\p_i(t) ,  \\ 
		& {{\yb}}^{s}(t) = C{\xb}^{s}(t) \!+\! F {\vb}(t),
		\end{aligned}\right. \label{eqn:LPV_SSA_stoch}
	\end{eqnarray}
	and  $S_s$ is an asLSS realization of $(\yb^s,\btheta)$
	with the noise process $\vb$, $S_d$ is an asLSS realization of $(\yb^d,\btheta)$ with
	the noise process $\ub$.
	In particular, 
	\begin{align} 
	  &{\yb}_d(t)=\sum_{\sigma \in \Sigma, w \in \Sigma^{*}} \sqrt{p_{\sigma w}} CA_wB_{\sigma}\z^{\ub}_{\sigma w}(t) + D\ub(t)
	   %\nonumber \\
	   % & =E_l[\yb(t)  \mid \{ \z^{\ub}_{w}(t) \}_{w \in \Sigma^{+}} \cup \{\ub(t)\}] 
	  \label{decomp:lemma:eq1} \\
	 & {\yb}_s(t)=\sum_{\sigma \in \Sigma, w \in \Sigma^{*}} \sqrt{p_{\sigma w}} CA_wK_{\sigma}\z^{\vb}_{\sigma w}(t) + F\vb(t) 
	  \label{decomp:lemma:eq2} 
	\end{align}
	\end{Lemma}

 \begin{pf}[Proof of Lemma \ref{decomp:lemma}: Decomposition of $\yb$]
It follows, from the Lemmas ~\ref{decomp:lemma:pf1}, \ref{decomp:lemma:pf2}, \ref{decomp:lemma:pf3}, and Lemma \ref{decomp:lemma:pf4} that,
	\begin{equation}
		\label{decomp:lemma:pf:eq2}
		\begin{split} 
			& \xb^d(t+1)=E_l[\xb(t+1)  \mid \mathcal{H}_{t+1,+}^{\ub}]  = \\
			& = \sum_{\sigma \in \Sigma} (A_{\sigma} E_l[\xb(t)\p_{\sigma}(t)\mid \mathcal{H}_{t+1,+}^{\ub}] \\
			&+ B_{\sigma} E_l[\ub(t)\p_{\sigma}(t) \mid \mathcal{H}_{t+1,+}^{\ub}] 
			+
			K_{\sigma} E_{l}[\vb(t)\p_{\sigma}(t) \mid \mathcal{H}_{t+1,+}^{\ub}])
			% & \sum_{\sigma \in \Sigma} (A_{\sigma} E_l[\xb(t)\p_{\sigma}(t)\mid \mathcal{H}_{t+1,+}^{\ub}] + B_{\sigma} E_l[\ub(t)\p_{\sigma}(t) \mid \mathcal{H}_{t+1,+}^{\ub}]
		\end{split}
	\end{equation}
	Note that, $\ub(t)\p_{\sigma}(t)=\sqrt{p_{\sigma}} \z^{\ub}_{\sigma}(t+1)$, hence the components of $\ub(t)\p_{\sigma}(t)$ belong to $\mathcal{H}_{t+1,+}^{\ub}$ and
	therefore  
	\[ E_l[\ub(t)\p_{\sigma}(t) \mid \mathcal{H}_{t+1,+}^{\ub}]=\ub(t)\p_{\sigma}(t) \]
	
	We claim that, 
	\[ 
	E_l[\xb(t)\p_{\sigma}(t)\mid \mathcal{H}_{t+1,+}^{\ub}] = \xb^d(t)\p_{\sigma}.
	\]
	Since $\xb(t)=\xb^d(t)+\xb^s(t)$, it follows that $\xb(t)\p_{\sigma}(t)=\xb^d(t)\p_{\sigma}(t)+\xb^s(t)\p_{\sigma}(t)$. 
	
	From \cite[Lemma 9]{PetreczkyBilinear} 
	and Lemma \ref{decomp:lemma:pf3}--\ref{decomp:lemma:pf4}, it follows that
	\begin{equation}
		\label{decomp:lemma:pf:eq3}
		\begin{split}
			& \xb^d(t)\p_{\sigma}(t) = \sum_{w \in \Sigma^{*}, \sigma^{'} \in \Sigma} \sqrt{p_{\sigma^{'} w\sigma}} A_w B_{\sigma^{'}}\z^{\ub}_{\sigma^{'} w\sigma}(t+1) \\
			& \xb^s(t)\p_{\sigma}(t) = \sum_{w \in \Sigma^{*}, \sigma^{'} \in \Sigma} \sqrt{p_{\sigma^{'} w\sigma}} A_w K_{\sigma^{'}}\z^{\vb}_{\sigma^{'} w\sigma}(t+1) \\
		\end{split}
	\end{equation}
	
	From Lemma \ref{decomp:lemma:pf2}, it follows that $\z^{\vb}_{\sigma^{'} w\sigma}(t+1)$ is orthogonal to $\mathcal{H}_{t+1,+}^{\ub}$ for all $w \in \Sigma^{+}$,
	$\sigma^{'},\sigma \in \Sigma$, and hence  $\xb_s(t)\p_{\sigma}(t)$ is also orthogonal to $\mathcal{H}_{t+1,+}^{\ub}$.  Moreover, since the components of
	$\z^{\ub}_{\sigma^{'} w\sigma}(t+1)$ belong to $\mathcal{H}_{t+1,+}^{\ub}$, it follows that  $\xb^d(t)\p_{\sigma}(t)$ belongs to $\mathcal{H}_{t+1,+}^{\ub}$.
	Hence,
	\[ 
	\begin{split} 
		& E_l[\xb(t)\p_{\sigma}(t)\mid \mathcal{H}_{t+1,+}^{\ub}] \\ &=E_l[\xb^s(t)\p_{\sigma}(t) \mid  \mathcal{H}_{t+1,+}^{\ub}] 
		+E_l[\xb^d(t)\p_{\sigma}(t) \mid  \mathcal{H}_{t+1,+}^{\ub}]\\
		&=E_l[\xb^d(t)\p_{\sigma} \mid  \mathcal{H}_{t+1,+}^{\ub}]=\xb^d(t)\p_{\sigma}(t). 
	\end{split}
	\]
	
	Finally, from Lemma \ref{decomp:lemma:pf2}, it follows that $\vb(t)\p_{\sigma}(t)=\sqrt{p_{\sigma}} \z_{\sigma}^{\vb}(t+1)$ is orthogonal to $\mathcal{H}_{t+1,+}^{\ub}$, and hence
	\[  E_{l}[\vb(t)\p_{\sigma}(t) \mid \mathcal{H}_{t+1,+}^{\ub}]=0. \]
	
	By collecting all these facts, we can show that 
	\[ \begin{split} 
		& \xb^d(t+1)=E_l[\xb(t+1)  \mid \mathcal{H}_{t+1,+}^{\ub}]  = \\
		& = \sum_{\sigma \in \Sigma} (A_{\sigma} E_l[\xb(t)\p_{\sigma}(t)\mid \mathcal{H}_{t+1,+}^{\ub}] +\\
		& B_{\sigma} E_l[\ub(t)\p_{\sigma}(t) \mid \mathcal{H}_{t+1,+}^{\ub}]
		+
		K_{\sigma} E_{l}[\vb(t)\p_{\sigma}(t) \mid \mathcal{H}_{t+1,+}^{\ub}])\\=
		& \sum_{\sigma \in \Sigma} (A_{\sigma}\xb^d(t)+B_{\sigma}\ub(t)) 
	\end{split}
	\]
	
	That is,  the first equation of \eqref{eqn:LPV_SSA_deter} holds.
	
	As to the second equation of \eqref{eqn:LPV_SSA_deter}, we know from \cite{RouphaelPetreczkyArxiv}, that
	\[
	\begin{split}
		&\yb^d(t)=E_l[\yb(t) \mid \mathcal{H}_{t,+}^{\ub}]=\\&=CE_l[\xb(t) \mid \mathcal{H}_{t,+}^{\ub}]+D E_l[\ub(t) \mid \mathcal{H}_{t,+}^{\ub}]+E_l[\vb(t) \mid \mathcal{H}_{t,+}^{\ub}].
	\end{split}
	\]
	
	Since  the components of $\ub(t)$ are among the generators of $\mathcal{H}_{t,+}^{\ub}$, and by Lemma \ref{decomp:lemma:pf1}, $\vb(t)$ is orthogonal to 
	$\mathcal{H}_{t,+}^{\ub}$, $E_l[\vb(t) \mid \mathcal{H}_{t,+}^{\ub}]=0$ and $E_l[\ub(t) \mid \mathcal{H}_{t,+}^{\ub}]=\ub(t)$, and
	 hence the second equation of \eqref{eqn:LPV_SSA_deter} holds. 
	
	From $\yb^s(t)=\yb(t)-\yb^d(t)$, $\xb^s(t)=\xb(t)-\xb^d(t)$ and \eqref{eqn:LPV_SSA_deter}, \eqref{eqn:LPV_SSA_stoch} follows.
	
	Finally, notice that as $\sum_{\sigma \in \sigma} p_{\sigma} A_{\sigma} \otimes A_{\sigma}$ is stable and $\ub$ is a white noise w.r.t. $\p$, $\mathcal{S}_d$ is a well-defined asLSS  with the noise process $\ub$. Moreover, \eqref{eqn:LPV_SSA_deter}
	implies that it is a realization of $(\yb^d,\btheta)$ with the noise process $\ub$.
	Similarly, since  $\sum_{\sigma \in \sigma} p_{\sigma} A_{\sigma} \otimes A_{\sigma}$ is stable and $\vb$ is a white noise w.r.t. $\p$, it follows that $\mathcal{S}_s$ is as asLSS and from
	\eqref{eqn:LPV_SSA_stoch} it follows
	that $\mathcal{S}_s$ is a realization of 
	$(\yb^s,\btheta)$. 
	
	To conclude, \eqref{decomp:lemma:eq1} and
	\eqref{decomp:lemma:eq2} follow from the fact that $\mathcal{S}_d$ is a realization of $(\yb^d,\btheta)$ and $\mathcal{S}_s$ is a realization of $(\yb^s,\btheta)$, by using
	\cite[Lemma 1]{PetreczkyBilinear}.
	
\begin{comment}
	It is left to show that $\mathcal{S}_d$ and $\mathcal{S}_s$ are asLSS of $(\yb^d,\btheta)$ and $(\yb^s,\btheta)$ respectively, and $\xb^d$, $\xb^s$
        are their unique state processes. 
% $(\{A_{\sigma},B_{\sigma}\}_{\sigma \in \Sigma},C,D,\ub)$ and $(\{A_{\sigma},K_{\sigma}\}_{\sigma \in \Sigma},C,I_{\ny},\xb^s,\vb)$
%	are stationary LSS representations without inputs as per Definition \ref{defn:LPV_SSA_wo_u}. 
        This follows from the fact that	
	$\sum_{\sigma \in \Sigma} p_{\sigma} A_{\sigma} \otimes A_{\sigma}$ is stable and $\ub$ and $\vb$ are both white noise processes w.r.t. $\p$,         and from  \eqref{decomp:lemma:pf3:eq2}, \eqref{decomp:lemma:pf4:eq2} and \cite[Lemma 3]{PetreczkyBilinear}. \hfill $\blacksquare$
	%only thing which needs to be shown is that $\begin{bmatrix} (\xb^d)^T & \ub^T \end{bmatrix}^T$ and $\begin{bmatrix} (\xb^s)^T & \vb^T \end{bmatrix}^T$
	%are ZMWSII. However, the latter follows
\end{comment}

\end{pf}

%$\mathcal{H}_{t}^{\vb}$ the close subspace of $\mathcal{H}_1$  generated by the components
\begin{Lemma}\label{decomp:lemma:inv:pf2.1}
	The  components of  $\yb^s(t), \z^{\yb^s}_v(t)$, $\eb^s(t)$, $\z^{\eb^s}_v(t)$, $v \in \Sigma^{+}$ belong to $\mathcal{H}_{t,+}^{\vb}$.
\end{Lemma}
\begin{pf}[Proof of Lemma~\ref{decomp:lemma:inv:pf2.1}]
	%Assume that $(\{A_{\sigma}, B_{\sigma}, K_{\sigma}\}_{\sigma \in \Sigma},C,D,\xb,\vb)$  is a stationary LSS representation of $\yb$ with input $\ub$. By Assumption of Lemma \ref{decomp:lemma:inv}, such a LSS representation exists. 
	Recall \eqref{decomp:lemma:pf4:eq2} from Lemma \ref{decomp:lemma:pf4},
	%\begin{equation*}
%		\xb^s(t) = \sum_{w \in \Sigma^{*}, \sigma \in \Sigma} \sqrt{p_{\sigma w}} A_w K_{\sigma}\z^{\vb}_{\sigma w}(t),
%	\end{equation*}
	as $\yb^s(t)=C\xb^s(t)+\vb(t)$, 
	%\begin{equation*}
	%	\yb^s(t) =  \sum_{w \in \Sigma^{*}, \sigma \in \Sigma} \sqrt{p_{\sigma w}} C A_w K_{\sigma}\z^{\vb}_{\sigma w}(t)+ \vb(t),
	%\end{equation*}
	the components of $\yb^s(t)$ belong to $\mathcal{H}_{t,+}^{\vb}$. 
	In particular, from \cite[Lemma 11]{PetreczkyBilinear}, it follows that the coordinates of
	$\z^{\yb^s}_w(t)$ belong to $\mathcal{H}_{t,+}^{\vb}$ and hence, $\mathcal{H}_{t}^{\yb^s} \subseteq \mathcal{H}_{t}^{\vb}$.
	Since $\eb^s(t)=\yb^s(t)-E_l[\yb^s(t) \mid \mathcal{H}_{t}^{\yb^s}] $, it follows that the components of 
	$\eb^s(t)$ belong to $\mathcal{H}_{t,+}^{\vb}$. Since $\z_{v}^{\vb}(t)=\z_{v1}^{\vb}(t+1)$, $\vb(t)=\z_1^{\vb}(t+1)$, it follows that
	$\mathcal{H}_{t,+}^{\vb} \subseteq \mathcal{H}_{t+1}^{\vb}$ and from  \cite[Lemma 11]{PetreczkyBilinear} it follows that
	the components of $\z_v^{\eb^s}(t)$ belong to $\mathcal{H}_{t}^{\vb} \subseteq \mathcal{H}_{t,+}^{\vb}$. 
\end{pf}
\begin{Lemma}
	\label{decomp:lemma:inv:pf2.2}
	If there exists a sLSS  $(\yb,\ub,\btheta)$ , then the components of 
	$\yb^s(t), \z^{\yb^s}_v(t), \eb^s(t), \z^{\eb^s}_v(t)$, $v \in \Sigma^{+}$ are orthogonal to $\mathcal{H}_{t,+}^{\ub}$.
	i.e., for all $v,w \in \Sigma^{+}$
	%\begin{equation}
	%	\label{decomp:lemma:inv:pf2:eq1}
	%	\begin{split}
	%		&E[\eb^s(t)(\z^{\ub}_w(t))^T]=0, ~ 
	%		E[\ub(t)(\z^{\eb^s}_v(t))^T]=0, ~  \\
	%		&E[\z_w^{\ub}(t)(\z^{\eb^s}_v(t))^T]=0, \\
	%		&E[\yb^s(t)(\z^{\ub}_w(t))^T]=0, ~ 
	%		E[\ub(t)(\z^{\yb^s}_v(t))^T]=0, ~  \\
	%		& E[\z_w^{\ub}(t)(\z^{\yb^s}_v(t))^T]=0 \\
	%	\end{split}
	%\end{equation}
\end{Lemma}
\begin{pf}[Proof of Lemma~\ref{decomp:lemma:inv:pf2.2}]
	From Lemma \ref{decomp:lemma:pf1}--\ref{decomp:lemma:pf2} and by noticing that $\vb(t)=\vb(t)\p_1(t)$ it follows that the elements of $\mathcal{H}_{t,+}^{\vb}$ are orthogonal to $\mathcal{H}_{t,+}^{\ub}$.
	% and the latter space is generated by the coordinates of  $\ub(t)$, $\z_w^{\ub}(t)$, $w \in \Sigma^{+}$.
	Hence, the coordinates of $\ub(t)$, $\z_w^{\ub}(t)$, $w \in \Sigma^{+}$ are orthogonal to $\mathcal{H}_{t,+}^{\vb}$. 
	Since the coordinates of  $\yb^s(t), \z^{\yb^s}_v(t), \eb^s(t),\z^{\eb^s}_v(t)$ belong to $\mathcal{H}_{t,+}^{\vb}$, it follows 
	that the coordinates of 
	$\yb^s(t), \z^{\yb^s}_v(t), \eb^s(t),\z^{\eb^s}_v(t)$
	are orthogonal to  $\mathcal{H}_{t,+}^{\ub}$. 
	Since $\mathcal{H}_{t,+}^{\ub}$ is generated by the coordinates of  $\ub(t)$, $\z_w^{\ub}(t)$, $w \in \Sigma^{+}$, 
	the statement of the lemma follows. %\eqref{decomp:lemma:inv:pf2:eq1} follows.
\end{pf}
\begin{Lemma}%[Correlation analysis \citep{CoxLPVSS}]
	\label{thm:cra}
	Assume that there exists an sLSS of $(\yb,\ub,\btheta)$.
   Then, $(\{A_{\sigma},B_{\sigma}\},C,D,\ub)$ is an asLSS of $(\yb^d,\btheta)$, if and only if
the  dLSS $(\{A_{\sigma},B_{\sigma}\}_{\sigma \in \Sigma},C,D)$ is a realization of $\incov$ and 
$\sum_{i=1}^{\pdim} p_i A_i \otimes A_i$ is stable. 
%\begin{align}
%&C A_{\sigma_{k}} \cdots A_{\sigma_{2}} B_{\sigma_{1}} =  \incov{w}\Sigma^{-1}_{u}, \\
% &D =\frac{\incov(\epsilon)}{\sqrt{p_{\epsilon}}} \Sigma^{-1}_{u}= \expect{\yb(t) \ub^{T}(t)} \Sigma^{-1}_{u}, 
%& \forall w \in \Sigma^{*}, \sigma \in \Sigma: C A_{w} B_{\sigma } =  \incov(\sigma w), \\
%&D = \incov{\epsilon}\Sigma^{-1}_{u}=\expect{\yb(t) \ub^{T}(t)} \Sigma^{-1}_{u},  
%& D = \incov(\epsilon)=\expect{\yb(t) \ub^{T}(t)} \varu^{-1}, 
%\end{align}
%i.e. 
%$\incov$  in \eqref{eqn:input_covar} equals the sub-Markov function $M_{\mathscr{S}}$ \eqref{eqn:sub_markov} of the deterministic LSS representation
%where  $\zwu$ is as defined in \eqref{eqn:zwy}, 
%where $\Sigma_{u} \!=\! \mathrm{var}(\ub)$, 
%if the assumptions A1-A2 hold.

\begin{comment}
If 
$(\{\hat{A}_{\sigma},\hat{B}_{\sigma}\}_{\sigma \in \Sigma},\hat{C},\hat{D})$ is a minimal dLSS realization of $\incov$, 
%and there exists a  sLSS of $(\yb,\ub,\btheta)$, 
then 
$\sum_{i=1}^{\pdim} p_i \hat{A}_i \otimes \hat{A}_i$ is stable. 
%$(\{\hat{A}_{\sigma},\hat{B}_{\sigma}\},\hat{C},\hat{D},\ub)$ is a asLSS of $\yb^d$. 
\end{comment}

%such that its sub-Markov function $M_{\hat{\mathscr{S}}}$ equals $\incov$ and it is minimal dimensional among such deterministic LSS representations, then 
\end{Lemma}

\begin{pf}[Proof of Lemma \ref{thm:cra}]
%\emph{If $(\{A_{\sigma},B_{\sigma}\}_{\sigma \in \Sigma},C,D,\xb,\ub)$ is a stationary LSS representation without inputs of $\yb^d$ $\implies$ $\incov=M_{\mathscr{S}}$, where $\mathscr{S}=(\{A_{\sigma},B_{\sigma}\}_{\sigma \in \Sigma},C,D)$. }
	
{\textbf{Part(1)}	}Assume that $(\{A_{\sigma},B_{\sigma}\}_{\sigma \in \Sigma},C,D,\ub)$ is an asLSS of $(\yb^d,\btheta)$.  Then $\sum_{i=1}^{\pdim} p_i A_i \otimes A_i$ is stable by definition of asLSS. 
	From Lemma \ref{decomp:lemma:inv:pf2.2} it follows that
	that the he components of $\yb^s(t)$ is orthogonal to $\mathcal{H}^{\ub}_{t,+}$, that is,
	$E[\yb^s(t) (\z^{\ub}_w(t))^T] = 0$ and  $E[\yb^s(t) (\ub(t))^T] = 0$. 
	Thus, $E[\yb(t) (\z^{\ub}_w(t))^T] = E[\yb^d(t)(\z^{\ub}_w(t))^T]$, $E[\yb(t) (\ub(t))^T] = E[\yb^d(t)(\ub(t))^T]$. Hence,
	From \eqref{decomp:lemma:pf3:eq2} of Lemma \ref{decomp:lemma:pf3} and the fact that $\yb^d(t)=C\xb^d(t)+D\ub(t)$ it follows that
	%\begin{equation*}
	%	\yb^d(t) =  \sum_{s \in \Sigma^{*}, \sigma \in \Sigma} \sqrt{p_{\sigma w}} C A_s B_{\sigma}\z^{\ub}_{\sigma s}(t)+ D \ub(t)
	%\end{equation*}
	%Consider, 
	\[ 
	\begin{split}
		&E[\yb^d(t) (\z_w^{\ub}(t))^T] = \\
		&\sum_{s \in \Sigma^{*}, \sigma \in \Sigma} \sqrt{p_{\sigma s}} C A_s B_{\sigma}E[\z^{\ub}_{\sigma s}(t) (\z_w^{\ub}(t))^T] \!\! +\!\! 
		DE[ \ub(t) (\z_w^{\ub}(t))^T] \\
		&= \sqrt{p_{\sigma s}} C A_s B_{\sigma} \Lambda_{\ub} = \Psi_{\yb,\ub}(\sigma s)\Lambda_u\sqrt{p_{\sigma s}}
	\end{split}
	\]
	This follows, as  $\ub$ is white noise process, $E[ \ub(t) (\z_w^{\ub}(t))^T] =0$ and $E[\z^{\ub}_{\sigma s}(t) (\z_w^{\ub}(t))^T] = \Lambda_{\ub}$ if $\sigma s = w$, otherwise $E[\z^{\ub}_{\sigma s}(t) (\z_w^{\ub}(t))^T] = 0$.
	Similarly,  $E[\yb^d(t)\ub(t)]= D \Lambda_{\ub} = \Psi_{\yb,\ub}(\epsilon)\Lambda_{\ub}$. 
{Hence, $\Psi_{\yb,\ub} = CA_\sigma B_\sigma$, i.e., $(\{A_\sigma,B_\sigma\}_{\sigma\in\Sigma},C,D)$ is a realization of $\Psi_{\yb,\ub}$.
 
 \textbf{Part(2)} Assume that $\{(A_\sigma,B_\sigma\}_{\sigma\in\Sigma},C,D)$ is a dLSS realization of $\Psi_{\yb,\ub}$ and $\sum_{i=1}^{\pdim} p_i A_i \otimes A_i$ is stable. 
 Consider the asLSS $\mathcal{S} = (\{A_\sigma,B_\sigma\}_{\sigma\in\Sigma},C,D,\ub)$. 
 $\mathcal{S}$ is a well-defined asLSS, as $\ub$ is a white noise w.r.t. $\p$. 

 Let $\Tilde{\xb}$ be the unique trajectory of $\mathcal{S}$, and let $\Tilde{\yb}$ be the unique output of $\mathcal{S}$. From \cite{PetreczkyBilinear} it follows that \[\Tilde{\yb}(t) = \sum_{N=0}^{\infty}\sum_{\tiny \begin{matrix} |w|=N \\ w\in\Sigma^*\end{matrix}}\sum_{\sigma\in\Sigma}\sqrt{p_{\sigma w}}CA_wB_\sigma\mathbf{z}^{\mathbf{u}}_{\sigma w}(t) + D\ub(t)\]
 using $\Psi_{\yb,\ub}(\sigma s)=CA_sB_\sigma$, it follows that 
 \begin{equation}
 \label{part2:eq1}
     \Tilde{\yb}(t) = \sum_{N=0}^{\infty}\sum_{\tiny \begin{matrix} |w|=N \\ w\in\Sigma^*\end{matrix}}\sum_{\sigma\in\Sigma}\sqrt{p_{\sigma w}}\Psi_{\yb,\ub}(\sigma s)\mathbf{z}^{\mathbf{u}}_{\sigma w}(t) + \Psi_{\yb,\ub}(\epsilon)\ub(t)
 \end{equation}

 By assumption, there exists a sLSS $\mathcal{S} = (\{A_\sigma',B_\sigma'\}_{\sigma\in\Sigma},C',D',\ub)$, where $\mathcal{S}$ is a realization of $(\yb,\ub,\btheta)$. 
 By Lemma \ref{decomp:lemma}, $\mathcal{S}_d = (\{A_\sigma',B_\sigma'\}_{\sigma\in\Sigma},C',D',\ub)$ is an asLSS realization of $(\yb^d,\btheta)$. 
 Hence, by \textbf{Part(1)} of the lemma, $(\{A_\sigma',B_\sigma'\}_{\sigma\in\Sigma},C',D')$ is a dLSS realization of $\Psi_{\yb,\ub}$, where $\Psi_{\yb,\ub}(\sigma s) = C'A_s'B_\sigma'$ and $\Psi_{\yb,\ub}(\epsilon) = D'$.
 
 But from \cite{PetreczkyBilinear}, it follows that \[\yb^d(t) = \sum_{N=0}^{\infty}\sum_{\tiny \begin{matrix} |w|=N \\ w\in\Sigma^*\end{matrix}}\sum_{\sigma\in\Sigma}\sqrt{p_{\sigma w}}C'A_w'B_\sigma'\mathbf{z}^{\mathbf{u}}_{\sigma w}(t) + D'\ub(t)\]
 \begin{equation}
 \label{part2:eq2}
     \yb^d(t) =\sum_{\tiny{N=0}}^{\infty}\sum_{\tiny \begin{matrix} |w|=N \\ w\in\Sigma^*\end{matrix}}\sum_{\sigma\in\Sigma}\sqrt{p_{\sigma w}}\Psi_{\yb,\ub}(\sigma s)\mathbf{z}^{\mathbf{u}}_{\sigma w}(t) + \Psi_{\yb,\ub}(\epsilon)\ub(t)
 \end{equation}

Hence, by combining \eqref{part2:eq1} and \eqref{part2:eq2} it follows that $\Tilde{\yb}(t) = \yb^d(t)$.
That is, $\mathcal{S} = (\{A_\sigma,B_\sigma\}_{\sigma\in\Sigma},C,D,\ub)$ is a realization of $(\yb,\btheta)$.}
\end{pf}
{
\begin{pf}[Proof of Lemma \ref{min:col:lem1}] Assume the sLSS of the form \eqref{eq:aslpv} is a realization of $(\yb,\ub,\btheta)$. 
From Lemma \ref{decomp:lemma}, it follows that $\mathcal{S}_s$ defined in \eqref{eqn:LPV_SSA_stoch} is an asLSS realization of $(\yb^s,\btheta)$, and $\mathcal{S}_d$ is an asLSS realization of $(\yb^d,\btheta)$ as defined in \eqref{eqn:LPV_SSA_deter}.
Consider the dLSS $\mathscr{S}_s$ as defined in \eqref{eqn:LPV_SSA:det}.
From \cite[Lemma 4]{PetreczkyBilinear} it follows that, in the terminology of \cite{PetreczkyBilinear}, the $\mathcal{R} = (\{\sqrt{p_i}A_i\},[G_i \dots G_n],C)$ is a recognizable representation of the formal power series
\begin{equation}
\label{FormalPowerSeries}
    s \mapsto \begin{bmatrix} \Psi_{\yb^s}(1s) \dots \Psi_{\yb^s}(\pdim s) \end{bmatrix}
\end{equation}
Note that $\mathcal{R}$ corresponds to 
Hence by \cite[Lemma 4]{PetreczkyBilinear},
\begin{equation}
\label{Psi_ys}
    \Psi_{\yb^s}(\sigma s) = \sqrt{p_w}CA_sG_\sigma
\end{equation}
where $w=\sigma s \in \Sigma^+$, and $\Psi_{\yb^s}(w) = I$ for $w = \epsilon$. 
From Lemma \ref{thm:cra} it follows that $(\{A_\sigma, B_\sigma\}_{\sigma=1}^{\pdim}, C, D)$ is a dLSS realization of $\Psi_{\yb,\ub}$.
Hence
\begin{equation}
\label{Psi_y_u}
    \Psi_{\yb,\ub}=\left\{\begin{array}{ll}
		CA_sB_{\sigma} \ \ & w=\sigma s, \, \sigma \in \Sigma, \, s \in \Sigma^{*} \\
		D            \ \  & w=\epsilon
	\end{array}\right.
\end{equation}
Combining \eqref{Psi_ys} and \eqref{Psi_y_u} results in having 
\begin{equation}
    \begin{aligned} 
    M_{\yb,\ub}(\sigma s) &= \begin{bmatrix} \sqrt{p_{\sigma s}}\Psi_{\yb,\ub}(\sigma s) & \Psi_{\yb^s}(\sigma s) \end{bmatrix} \\
    &= \sqrt{p_{\sigma s}}CA_s \begin{bmatrix} B_\sigma & G_\sigma \end{bmatrix} \\
    M_{\yb,\ub}(\epsilon) &=  \begin{bmatrix} \sqrt{p_{\sigma s}}\Psi_{\yb,\ub}(\epsilon) & \Psi_{\yb^s}(\epsilon) \end{bmatrix} \\
    &= \begin{bmatrix} D & I \end{bmatrix}
    \end{aligned}
\end{equation}
That is, $\mathscr{S}_s$ is a realization of $M_{\yb,\ub}$.
\end{pf}}

\subsection{Proof of Lemma \ref{min:col:lem2}}
% where, $(\{\hat{A}_{\sigma}, \hat{B}_{\sigma}, \hat{K}_{\sigma}\}_{\sigma \in \Sigma},\hat{C},\hat{D},\hat{\xb},\eb^s)$ are as defined in  \eqref{decomp:lemma:sys:eq1}-\eqref{decomp:lemma:innov2}. 
Assume that $\mathcal{S}$ of the form \eqref{eq:aslpv}  is a sLSS of $(\yb,\ub,\btheta)$  \footnote{By Assumption of Lemma \ref{decomp:lemma:inv}, such a LSS representation exists.}.
Recall from Notation \ref{hilbert:notation} the Hilber-space $\mathcal{H}_1$. % denotes the Hilbert-space of zero mean square integrable random variables.
Denote by $\mathcal{H}_{t,+}^{\vb}$ and $\mathcal{H}_t^{\vb}$  the closed-subspace of $\mathcal{H}_1$ 
generated by the components of $\{\z^{\vb}_w\}_{w \in \Sigma^{+}} \cup \{\vb(t)\}$ and $\{\z^{\vb}_w\}_{w \in \Sigma^{+}}$ respectively.
%$\mathcal{H}_{t}^{\vb}$ the close subspace of $\mathcal{H}_1$  generated by the components 

\begin{Lemma}\label{decomp:lemma:inv:pf2}
	$\r=\begin{bmatrix}  (\eb^s)^T & \ub^T \end{bmatrix}^T$ is a white noise process w.r.t. $\p$ and $E[\eb^s(t)\ub^T(t)\p_{\sigma}^2(t)]=0$ for all $\sigSet$.
\end{Lemma}
\begin{pf}[Proof of Lemma \ref{decomp:lemma:inv:pf2}]
	In order to prove the statement of the lemma, we will first show that $\r$ is a ZMWSII, by showing that $\r$  satisfies the conditions of  Definition \ref{def:ZMWSSI} one by one. 
	First, we show that the processes $\r(t),\z_w^{\r}(t),w \in \Sigma^{+}$ is zero mean, square integrable.
	
	By assumption  $\ub$ is a white noise process w.r.t. $\p$.
	From the fact that $\hat{\mathcal{S}}_s$ is an asLSS of $(\yb^s,\btheta)$ it follows that $\eb^s$ is also
	a white noise process w.r.t. $\p$. 
        In particular, $\ub$, $\p$ are  both also ZWMSII and thus $\eb^s(t),\ub(t), \z_w^{\eb^s}(t),\z_{w}^{\ub}(t)$, $w \in \Sigma^{+}$ is zero mean, square integrable. From this it follows that  $\r(t)$ and
	$\z_w^{\r}(t)$ %= \begin{bmatrix}  (\z^{\eb^s}_w(t))^T & (\z^{\ub}_w(t))^T \end{bmatrix}^T$ 
are zero mean and square integrable.

	From Lemma \ref{decomp:lemma:inv:pf2.1} it follows that the components $\eb^s(t)$ belongs to $\mathcal{H}_{t,+}^{\vb}(t)$, where $\vb$ is a noise process of an asLSS of
	$(\yb,\ub,\btheta)$. From the definition of a sLSS it follows that $\mathbf{w}=\begin{bmatrix} \vb^T & \ub^T \end{bmatrix}^T$ is ZMWSII. 
	Hence, with the notation of 
	Definition \ref{def:ZMWSSI}, the $\sigma$-algebras $\mathcal{F}_t^{\mathbf{w}}$ and $\mathcal{F}_t^{\p,+}$ are conditionally independent  
      w.r.t. $\mathcal{F}_t^{\p,-}$. 
	From the fact that $\eb^s(t)$ belongs to $\mathcal{H}_{t,+}^{\vb}$ it follows that $\eb^s(t)$ is measurable with respect to the $\sigma$-algebra generated by
	$\{\vb(t)\} \cup \{ \z_v^{\vb}(t)\}_{v \in \Sigma^{+}}$ and the latter $\sigma$-algebra is a subset of $\mathcal{F}^{\mathbf{w}}_t \lor \mathcal{F}^{\p,-}_t$, 
	where for two $\sigma$-algebras $\mathcal{F}_i$, $i=1,2$, $\mathcal{F}_1 \lor \mathcal{F}_2$ denotes the smallest $\sigma$-algebra generated by the $\sigma$-algebras $\mathcal{F}_1,\mathcal{F}_2$. 
	That is, $\eb^s(t)$ is measurable w.r.t. the $\sigma$ algebra $\mathcal{F}^{\mathbf{w}}_t \lor \mathcal{F}^{\p,-}_t$.
Since  $\mathcal{F}^{\mathbf{w}}_t$ and $\mathcal{F}^{\p,+}$ are conditionally independent w.r.t.  $\mathcal{F}^{\p,-}_t$, from \cite[Proposition 2.4]{vanputten1985} it follows that
	$\mathcal{F}^{\mathbf{w}}_t \lor \mathcal{F}^{\p,-}_t$ and $\mathcal{F}^{\p,+}_t$ are conditionally independent w.r.t. 
	$\mathcal{F}^{\p,-}_t$, and as $\mathcal{F}^{\r}_t \subseteq \mathcal{F}^{\mathbf{w}}_t \lor \mathcal{F}^{\p,-}_t$, %the $\sigma$-algebras
	it follows that   $\mathcal{F}^{\r}_t$ and $\mathcal{F}^{\p,+}$ are conditionally independent w.r.t. $\mathcal{F}^{\p,-}_t$.

	Finally, we show that $\r(t),\z_w^{\r}(t),w \in \Sigma^{+}$
	are   jointly wide-sense stationary, i.e., for all $s,t,\in \mathbb{Z}$, $s \le t$, $v,w \in  \Sigma^{+}$, 
	\begin{align*}
		\expect{\r(t)(\zwr(s))^{T}} &= \expect{\r(t-s) (\zwr(0))^{T}}, \\
		\expect{\r(t)(\r(s))^{T}} &= \expect{\r(t-s)(\r(0))^{T}}, \\
		\expect{\zwr(t)(\zvr(s))^{T}} &= \expect{\zwr(t-s) (\zvr(0))^{T}}.
	\end{align*}
 
	We will show only the last equality, the proofs of the first two are analogous.
	Indeed, from Lemma \ref{decomp:lemma:inv:pf2} and
       from the fact that $\eb^s$, $\ub$ are ZMWSII it follows the processes $\eb^s(t),\z_w^{\eb^s}(t),w \in \Sigma^{+}$ 
        are jointly wide-sense stationary
	and the processes $\ub(t),\z_w^{\ub}(t),w \in \Sigma^{+}$ are also jointly wide-sense stationary.
	From Lemma \ref{decomp:lemma:inv:pf2}  it follows that $\expect{\z^{\eb^s}_w(t+k)(\z_{v}^{\ub}(s+k))^T}$, and hence
	%\begin{align*}
		%&	\expect{\r(t+k)(\zwr(s+k))^{T}} =\\
		%&=\begin{bmatrix} \expect{\eb_s(t+k)(\z_{w}^{\eb_s}(s+k)} & 0 \\ 0 & \expect{\ub(t+k)(\z_{w}^{\ub}(s+k)} \end{bmatrix} \\
		%&=     \begin{bmatrix} \expect{\eb^s(t)(\z_{w}^{\eb^s}(s)} & 0 \\ 0 & \expect{\ub(t)(\z_{w}^{\ub}(s)} \end{bmatrix} \\
		%&=   \expect{\r(t) (\zwr(s))^{T}},
	%\end{align*}
%\begin{align*}
%		&	\expect{\r(t+k)(\r(s+k))^{T}} = \\
%		&= \begin{bmatrix} \expect{\eb^s(t+k)(\eb^s(s+k))^T} & 0 \\ 0 & \expect{\ub(t+k)(\ub(s+k))^T} \end{bmatrix} \\
%		&=     \begin{bmatrix} \expect{\eb^s(t)(\eb^s(s))^T} & 0 \\ 0 & \expect{\ub(t)(\ub(s))^T} \end{bmatrix} \\
%		&= \expect{\r(t)(\r(s))^{T}}, 
%	\end{align*}
\begin{align*}
& \expect{\zwr(t)(\zvr(s))^{T}} = \\
		& =\mathrm{diag}(\expect{\z^{\eb^s}_w(t)(\z_{v}^{\eb^s}(s))^T}, \expect{\z^{\ub}_w(t)(\z_{v}^{\ub}(s))^T})\\
		&= \mathrm{diag}(\expect{\z^{\eb^s}_w(t-s)(\z_{v}^{\eb^s}(0))^T}, \expect{\z^{\ub}_w(t-s)(\z_{v}^{\ub}(0))^T})\\
		& =  \expect{\zwr(t-s) (\zvr(0))^{T}}.
	\end{align*}
	%Above, we used the fact that if $s < t$, then $\ub(s+k)=\z_{h}^{\ub}(t+k+1)$, $\eb^s(s+k)=\z_{h}^{\eb^s}(t+k+1)$, $\ub(s)=\z_{h}^{\ub}(t+1)$, $\eb^s(s)=\z_{h}^{\eb^s}(t+1)$,
	%$\z^{\ub}_w(s+k)=\z_{wh}^{\ub}(t+k)$, $\z_{v}^{\eb^s}(s+k)=\z_{vh}^{\eb^s}(t+k)$, where $h=\underbrace{1\cdots 1}_{t-s}$. 
	
	That is, we have shown that $\r$ satisfies all the conditions of  Definition \ref{def:ZMWSSI}.

	Next we show that $\r$ is a white noise process w.r.t. $\p$, i.e., $E[\r(t)(\z^{\r}_w(t))^T]=0$ for all $w \in \Sigma^{+}$. From 
	Lemma \ref{decomp:lemma:inv:pf2} it follows that 
	\[
	\begin{split}
	& \expect{\r(t) (\zwr(t))^{T}}=  \\
	& \mathrm{diag}(\expect{\eb^s(t)(\z_{w}^{\eb^s}(t))^T}, \expect{\ub(t)(\z_{w}^{\ub}(t))^T})
        \end{split}
	\]
	Notice $\eb_s$ is a white noise process w.r.t. $\p$, since it is the noise process of  $\hat{\mathcal{S}}_s$, and the latter is an asLSS of $(\yb^s,\btheta)$, 
%$\Sigma_s=(\{\hat{A}^{s}_{i}, \hat{K}^{s}_{i} \}_{i=1}^{\pdim}, \hat{C}^{s}, I_{\ny},\hat{\xb}^s,\eb^{s})$ is a stationary
	%LSS representation of $\yb^s$ without inputs, 
        and hence  $\expect{\eb_s(t)(\z_{w}^{\eb_s}(t)}=0$. 
	Furthermore, $\ub$ is a white noise process w.r.t. $\p$ by assumption, so $\expect{\ub(t)(\z_{w}^{\ub}(t))^T}=0$. Hence, $\expect{\r(t) (\zwr(t))^{T}}=0$.

	It is left to show that $E[\eb^s(t)\ub^T(t)\p_{\sigma}^2(t)]=0$ for all $\sigSet$. Notice that $E[\eb^s(t)\ub^T(t)\p_{\sigma}^2(t)]=E[\z^{\eb^s}_{\sigma}(t+1)(\z^{\ub}_{\sigma}(t+1))^T]p_{\sigma}$ by definition, and from Lemma \ref{decomp:lemma:inv:pf2} it follows that $E[\z^{\eb^s}_{\sigma}(t+1)(\z^{\ub}_{\sigma}(t+1))^T]=0$.
\end{pf}

\begin{Lemma}
	\label{decomp:lemma:inv}
	Assume that there exists a sLSS of $(\yb,\ub,\btheta)$ and that the following holds:
        \begin{itemize}
        \item
	 $\hat{\mathcal{S}}_d=(\{\hat{A}_i, \hat{B}_{i} \}_{i=1}^{\pdim}, \hat{C}, \hat{D},\ub)$ is a 
	asLSS of $(\yb^d,\btheta)$ \footnote{Note that the input process of $\mathcal{S}_d$ is $\ub$.},
        \item
	 $\hat{\mathcal{S}}_s=(\{\hat{A}_{i}, \hat{K}_{i} \}_{i=1}^{\pdim}, \hat{C}, I_{\ny},\eb^{s})$
	is asLSS  of $(\yb^s,\btheta)$ in \emph{forward innovation form}, i.e., the $\eb^s$ is the innovation process of $\yb^s$.
       \end{itemize}
	%Define
	%\begin{equation} 
	%	\label{decomp:lemma:sys:eq1}
	%	\begin{split}
	%		& \hat{\xb}(t)=\begin{bmatrix} (\hat{\xb}^d(t))^T  & (\hat{\xb}^s(t))^T \end{bmatrix}^T \\ 
	%		& \hat{A}_{\sigma} = \mathrm{diag}(\hat{A}^{d}_{\sigma},\hat{A}^{s}_{\sigma}), ~ \hat{B}_{\sigma} = \left[ (\hat{B}^{d}_{\sigma})^{T} \ \mathbf{0}_{\nx \times \udim}^{T} \right]^{T} \\
	%		& \hat{K}_{\sigma} = \left[ \mathbf{0}_{\nx \times \ny}^{T} \  (\hat{K}^{s}_{\sigma})^{T}  \right]^{T}, ~ \hat{C} \!=\! \left[ \hat{C}^{d}  \ \ \hat{C}^{s}  \right], ~ \hat{D} \!=\!  \hat{D}^{d}.
	%	\end{split}
	%\end{equation}
       Then $\hat{\mathcal{S}}=(\{\hat{A}_{i},\hat{K}_{i},\hat{B}_{i}\}_{i=1}^{\pdim},\hat{C},\hat{D},\eb^s)$ is a sLSS of $(\yb,\ub,\btheta)$. 
   %Moreover, if $\hat{\xb}^d$, $\hat{\xb}^s$  are  the state processes of $\mathcal{S}_d$
   %    and $\mathcal{S}_s$ respectively, then $\hat{\xb}=\hat{\xb}^d+\hat{\xb}^s$ is the state process of $\hat{\mathcal{S}}$. 
	Moreover, the innovation process $\eb^s$ satisfies
	\begin{equation} 
		\label{decomp:lemma:innov2}
		\eb^s(t)=\yb(t)-E_l[\yb(t) \mid \{\z^{\yb}_w (t), \z^{\ub}_w(t) \}_{w \in \Sigma^{+}} \cup \{\ub(t)\}] 
	\end{equation}
\end{Lemma}

\begin{pf}[Proof of Lemma \ref{decomp:lemma:inv}]
	%It is clear from Lemmas \ref{decomp:lemma:inv:pf1}, \ref{decomp:lemma:inv:pf2}, \ref{decomp:lemma:inv:pf3} that 
	%From Lemma \ref{decomp:lemma:inv:pf3} it follows that 
	%$\{A_{\sigma}\}_{\sigma \in \Sigma}$ satisfies the conditions of the definition of a stationary LSS representation without inputs. From Lemma \ref{decomp:lemma:inv:pf2}, it follows that $\hat{\xb}$ and $\begin{bmatrix} (\eb^s)^T  & \ub^T \end{bmatrix}^T$ satisfies the conditions of an LSS representation without inputs.
	From Lemma \ref{decomp:lemma:inv:pf2} it follows that the noise process $\eb^s$ and the input $\ub$ satisfy the condition of $E[\eb^s(t)(\ub(t))^T\p_{\sigma}^2(t)]=0$, $\sigma \in \Sigma$.
 Since $\hat{\mathcal{S}}_s$ and $\hat{\mathcal{S}}_d$ are both asLSS, it follows that $\sum_{i=1}^{\pdim} p_{i} \hat{A}_i \otimes \hat{A}_i$ is stable.
Hence $\hat{\mathcal{S}}$ satisfies the conditions of a sLSS and it has a unique state process $\hat{\xb}$. 
 Finally, let $\hat{\xb}^s$ and $\hat{\xb}^d$ be the unique state processes of $\hat{\mathcal{S}}_s$ and $\hat{\mathcal{S}}_d$ respectively. We claim that
 $\hat{\xb}(t)=\hat{\xb}^d(t)+\hat{\xb}^s(t)$. Indeed, from  \eqref{stat:state:eq1}, it follows that 
\begin{equation*}
%\label{decomp:lemma:inv:pf1:eq1}
\begin{split}
& \hat{\xb}^d(t) = \sum_{w \in \Sigma^{*}, \sigma \in \Sigma} \sqrt{p_{\sigma w}} \hat{A}_w \hat{B}_{\sigma}\z^{\ub}_{\sigma w}(t), \\
& \hat{\xb}^s(t) = \sum_{w \in \Sigma^{*}, \sigma \in \Sigma} \sqrt{p_{\sigma w}} \hat{A}_w \hat{K}_{\sigma}\z^{\eb^s}_{\sigma w}(t),
\end{split}
\end{equation*}
from which it follows that
\begin{equation}
\label{decomp:lemma:inv:pf1:eq2}
   \hat{\xb}^d(t)+\hat{\xb}^s(t)=
    \sum_{w \in \Sigma^{*}, \sigma \in \Sigma} \sqrt{p_{\sigma w}} \hat{A}^d_w \begin{bmatrix} \hat{K}_{\sigma} & \hat{B}_{\sigma} \end{bmatrix} \z^{\r}_{\sigma w}(t),
\end{equation}
with $\r(t)=\begin{bmatrix} (\eb^s(t))^T & (\ub(t))^T \end{bmatrix}^T$. Then  \eqref{stat:state:eq1} applied to $\vb=\r$ means that
the right-hand side of \eqref{decomp:lemma:inv:pf1:eq2} equals $\hat{\xb}$. 
Finally, from $\hat{\xb}(t)=\hat{\xb}^s(t)+\hat{\xb}^d(t)$ and the fact that $\yb^d(t)=\hat{C}\hat{\xb}^d(t)+\hat{D}\ub(t)$ (as $\hat{\mathcal{S}}_d$ is an asLSS of $(\yb^d,\btheta)$) and $\yb^s(t)=\hat{C}\hat{\xb}^s(t)+\eb^s(t)$ (as $\hat{\mathcal{S}}_s$ is an asLSS of $(\yb^s,\btheta)$), it follows
that $\yb(t)=\hat{C}\hat{\xb}(t)+\hat{D}\ub(t)+\eb^s(t)$, i.e., $\hat{\mathcal{S}}$ is a sLSS of $(\yb,\ub,\btheta)$. 
%%	Hence, $(\{\hat{A}_{\sigma}, \hat{B}_{\sigma}, \hat{K}_{\sigma}\}_{\sigma \in \Sigma},\hat{C},\hat{D},\hat{\xb},\eb^s)$  is a stationary LSS representation of $\yb$ with input $\ub$. \hfill $\blacksquare$
\end{pf}

{
\begin{pf}[Proof of Lemma \ref{min:col:lem2}]
Using the terminology of \cite{PetreczkyBilinear}, consider the recognizable representation $\mathcal{R} = (\nx,\{A_\sigma\}_{\sigma \in \Sigma},\begin{bmatrix} G_1 \dots G_{\pdim} \end{bmatrix},C)$. 
$\mathcal{R}$ is the representation of the formal power series defined in \eqref{FormalPowerSeries}.
$\mathcal{R}$ is stable and observable since the observability of the representations is equivalent to the observability of the dLSS $(\{A_i, \begin{bmatrix} B_i & G_i \end{bmatrix}\}_{i=1}^{\pdim},C,D)$. 
The the existence of the limit $K_i = \lim_{I \rightarrow \infty} K_i^I$ follows from \cite[Lemma 21]{PetreczkyBilinear}, by noticing that the proof of \cite[Lemma 21]{PetreczkyBilinear} requires only the observability of the representation $\mathcal{R}$.
Moreover, from \cite[Lemma 21]{PetreczkyBilinear}, it follows that the asLSS $(\{A_i, K_i\}_{i=1}^{\pdim},C,I,\eb^s)$ is a realization of $(\yb^s,\btheta)$.
From Lemma \ref{thm:cra} it follows that the asLSS $(\{A_i, B_i\}_{i=1}^{\pdim},C,D,\ub)$ is a realization of $(\yb^d,\btheta)$.
hence by Lemma \ref{decomp:lemma:inv} it follows that the sLSS $(\{A_i, B_i,K_i\}_{i=1}^{\pdim},C,D,I,\eb^s)$ is a realization of $(\yb,\ub,\btheta)$.
\end{pf}
}

\subsection{Proof of Theorem \ref{theo:min}}\label{app:proof_theorem1}

{\begin{Lemma}
\label{Lemma:obs:stable}
    If $(\yb,\ub,\btheta)$ has an sLSS realization, and $\mathscr{S} = (\{\hat{A}_{\sigma},\hat{B}_{\sigma}\}_{\sigma \in \Sigma},\hat{C},\hat{D})$ is a minimal dLSS realization of $M_{\yb,\ub}$ then $\mathscr{S}$ is observable and $\sum_{i=1}^{\pdim} p_i \hat{A}_i \otimes \hat{A}_i$ is stable. 
\end{Lemma}

\begin{pf}[Proof of Lemma \ref{Lemma:obs:stable}]
Let $\mathscr{S}=(\{\hat{A}_\sigma, \hat{B}_{\sigma} \}_{\sigma=1}^{\pdim}, \hat{C}, \hat{D},\ub)$ be any dLSS realization of $M_{\yb,\yb}$. Consider the formal series $\Psi : s \mapsto \begin{bmatrix} M_{\yb,\ub}(1s) & \dots & M_{\yb,\ub}(\pdim s) \end{bmatrix}$, then $R_{\mathscr{S}'} = (\{\hat{A}_\sigma\}_{\sigma=1}^{\pdim}, \begin{bmatrix} \begin{bmatrix} \hat{B}_1 & \hat{G}_1 \end{bmatrix} & \dots & \begin{bmatrix} \hat{B}_{\pdim} & \hat{G}_{\pdim} \end{bmatrix} \end{bmatrix}, \hat{C})$ is a recognisable representation of $\Psi$ if and only if $\mathscr{S}'$ is a realization of $M_{\yb,\ub}$. 
It is easy to see that $R_{\mathscr{S}}$ is reachable (observable), if and only if $\mathscr{S}'$ satisfies the reachability and observability rank conditions of \cite{PetreczkyLPVSS}. 
Hence, $\mathscr{S}$ is minimal, if and only if $R_{\mathscr{S}}$ is minimal. 

%Since $(\yb,\ub,\btheta)$ has an sLSS $\mathcal{S}$ realization, the dLSS $\mathscr{S}_\mathcal{S}$ associated with it is a realization of $M_{\yb,\ub}$. 
We will proceed to show that $R_{\mathscr{S}}$ is stable in the sence of \cite{PetreczkyBilinear}, i.e., $\sum_{\sigma=1}^{\pdim} p_\sigma \hat{A}_\sigma \otimes \hat{A}_\sigma$ is stable. 
To this end, let $\mathcal{S}=(\{{A}_\sigma, {B}_{\sigma}, K_\sigma \}_{\sigma=1}^{\pdim}, {C}, {D},F,\vb)$ be an sLSS realization of $(\yb,\ub,\btheta)$ of the form of \eqref{eq:aslpv}. 
Then by Definition \ref{defn:LPV_SSA_wo_u}, $\sum_{\sigma=1}^{\pdim} p_\sigma {A}_\sigma \otimes {A}_\sigma$ is stable. 
Consider the dLSS $\mathscr{S}_\mathcal{S}$ associated with $\mathcal{S}$. By Lemma \ref{min:col:lem1}, $\mathscr{S}_\mathcal{S}$ is a realization of $M_{\yb,\ub}$. 
Let us consider the formal series $\Psi$ defined above and let us consider the recognisable representation $R_{\mathscr{S}_\mathcal{S}} = (\{\sqrt{p_\sigma}A_\sigma\}_{\sigma \in \Sigma}, \Tilde{G}, {C})$ where 
$\Tilde{G} = \begin{bmatrix} \begin{bmatrix} \sqrt{p_i}{B}_1 & {G}_1 \end{bmatrix} & \dots & \begin{bmatrix} \sqrt{p_{\pdim}}{B}_{\pdim} & {G}_{\pdim} \end{bmatrix} \end{bmatrix}$. 
With the same argument as above, it is easy ti see that $R_{\mathscr{S}_\mathcal{S}}$ is a representation of $\Psi$.
Since $R_{\mathscr{S}_\mathcal{S}}$ is stable in the term of \cite{PetreczkyBilinear} (as $\sum_{i=1}^{\pdim} \sqrt{p_\sigma} {A}_\sigma \otimes \sqrt{p_\sigma} {A}_\sigma = \sum_{i=1}^{\pdim} p_\sigma {A}_\sigma \otimes {A}_\sigma$ is stable), then by \cite[Theorem 6]{PetreczkyBilinear}, $\Psi$ is square summable, i.e., $$\sum_{N=0}^{\infty} \sum_{\tiny \begin{matrix} |w| = N\\ w\in \Sigma^{*} \end{matrix}} \Vert \Psi(w) \Vert_2^2 < +\infty$$
From \cite[Theorem 6]{PetreczkyBilinear} it follows that any minimal representation of $\Psi$ is stable. In particular, as $\mathscr{S}$ is minimal, hence $R_{\mathscr{S}}$ is minimal, therefore, $R_\mathscr{S}$ is stable, i.e., $\sum_{\sigma=1}^{\pdim} p_\sigma \hat{A}_\sigma \otimes \hat{A}_\sigma$ is stable.
The observability of $\mathscr{S}$ follows from its minimality.
\end{pf}}
That is, finding an asLSS of $(\yb^d,\btheta)$  boils down to finding a dLSS realization of $\incov$.

More precisely, consider sLSS of $(\yb,\ub,\btheta)$ of the form \eqref{eq:aslpv} and consider the associated asLSS $\mathcal{S}_s$
of $(\y^s,\btheta)$ as defined in Lemma \ref{decomp:lemma}, \eqref{eqn:LPV_SSA_stoch}.
%From~\cite[Lemma 4]{PetreczkyBilinear} the following Lemma follows easily.
\begin{Lemma}[Lemma 4 from \cite{PetreczkyBilinear}]
\label{lemma:stoch-real}
The dLSS $$(\{\sqrt{p_{\sigma}}A_{\sigma},\Bs_{\sigma}\}_{\sigma \in \Sigma},C,I)$$ is a realization of
$\covseq$.
%such that the sub-Markov function $M_{\mathscr{S}_{det}}$ of $\mathscr{S}_{det}$
% equals $\covseq$, \emph{i.e}, 
%\( \covseq(w) = \sqrt{p_w} CA_{w}\Bs, ~ \forall \wordSet{*} \). 
 %\[ \expect{\zwsigs(t)(\zwsigs(t))^{T}} = \frac{1}{\psig} (C\Psig C^{T} + \Qsig), ~ \forall \sigSet. \]
\end{Lemma}

%The proof of Lemma \ref{lemma:stoch-real1} is reported in Appendix~\ref{app:proof_lem5}].

\begin{pf}[Proof of Theorem \ref{theo:min}]
\textbf{Part \ref{theo:min:1}}.  Assume that $\mathcal{S}$ is span-reachable and observable sLSS of $(\yb,\ub,\btheta)$, 
then by definition the dLSS $\mathscr{S}_{\mathcal{S}}$ is also span-reachable and observable, and hence it is a 
minimal realization of $M_{\yb;\ub}$ by \cite{PetreczkyLPVSS}. Let $\mathcal{S}^{'}$ be another sLSS of $(\yb,\ub,\btheta)$. Then by Lemma \ref{min:col:lem1}
the dLSS  $\mathscr{S}_{\mathcal{S}^{'}}$ associated with $\mathcal{S}^{'}$ has the same sub-Markov parameters.
Hence, since $\mathscr{S}_{\mathcal{S}}$ is minimal dimensional, the dimension of $\mathscr{S}_{\mathcal{S}}$ cannot
exceed that of $\mathscr{S}_{\mathcal{S}^{'}}$, and since the dimension of a sLSS is the same as that of the dLSS associated with it, it follows that the dimension of $\mathcal{S}$
cannot exceed that of $\mathcal{S}^{'}$. 

\textbf{Part \ref{theo:min:2}.} To this end,
notice that the dLSS associated with $\mathcal{S}_m$  is $\mathscr{S}_m$, and since $\mathscr{S}_m$ is minimal, it is span-reachable and observable, hence $\mathcal{S}_m$ is span-reachable and observable, and hence minimal.  
The rest of the statement follows from  Lemma \ref{min:col:lem1}.

\textbf{Part \ref{theo:min:3}} Let Let  $\mathcal{S}=(\{A_{\sigma},B_{\sigma},K_{\sigma}\}_{\sigma=1}^{\pdim},C,D,F,\vb)$ be a sLSS realization of $(\yb,\ub,\btheta)$.
From Lemma \ref{min:col:lem1}, we define $\mathscr{S}_{\mathcal{S}}=(\{\sqrt{p_\sigma} A_\sigma,\begin{bmatrix} \sqrt{p_\sigma} B_\sigma & \Bs_\sigma \end{bmatrix} \}_{\sigma=1}^{\pdim},C,\begin{bmatrix} D & I_{\ny} \end{bmatrix}),$ the associated dLSS, which is a realization of $M_{\yb,\ub}$.
From $\mathscr{S}_{\mathcal{S}}$ as an input of Algorithm \ref{algo:Nice_select_true_deter_basic} and from Lemma \ref{basis_red:lemma}, we can get a minimal dLSS $\hat{\mathscr{S}} = (\{\hat{A}_{\sigma},\hat{B}_{\sigma}\}_{\sigma \in \Sigma},\hat{C},\hat{D})$ realization $M_{\yb,\ub}$.
We know, from Lemma \ref{Lemma:obs:stable}, that $\sum_{i=1}^{\pdim} p_i \hat{A}_i \otimes \hat{A}_i$ is stable.
It follows, from \textbf{part 2} of Theorem \ref{theo:min}, that the associated sLSS $\hat{\mathcal{S}}_{\mathscr{S}}$ is in innovation form. 
What remains is to prove that the latter sLSS is minimal. 
Note, from equations \eqref{min:eq-1}, \eqref{statecov:iter} and \eqref{lemma:stoch-real1:eq1}, that the associated dLSS $\hat{\mathscr{S}}_{\hat{\mathcal{S}}_{\hat{\mathscr{S}}}}$ and the minimal dLSS $\hat{\mathscr{S}}$ are the same when $N \rightarrow \infty$.
It follows that $\hat{\mathcal{S}}_{\hat{\mathscr{S}}}$ is minimal because his associated dLSS is minimal.
Therefore, if $(\yb,\ub,\btheta)$ has a sLSS realization, then there exists a minimal sLSS realization in innovation form.

\textbf{Part \ref{theo:min:4}.}  Let  $\mathcal{S}=(\{A_{\sigma},B_{\sigma},K_{\sigma}\}_{\sigma=1}^{\pdim},C,D,F,\vb)$  and  let $\widetilde{\mathcal{S}}=(\{\tilde{A}_{\sigma},\tilde{B}_{\sigma},\tilde{K}_{\sigma}\}_{\sigma=1}^{\pdim},\tilde{C},\tilde{D},\tilde{F},\vb)$ be of the form \eqref{eq:aslpv}, and 
%(\{\tilde{A}_i,\tilde{K}_i,\tilde{B}_i\}_{i=0}^{\pdim},\tilde{C},\tilde{D},I,\eb^s)$.
assume that $\mathcal{S}$ and $\widetilde{\mathcal{S}}$ are minimal sLSS of $(\yb,\ub,\p)$ in forward innovation form.  Let
$\mathscr{S}_{\mathcal{S}}=(\{\sqrt{p_i} A_i,\begin{bmatrix} \Bs_i,\sqrt{p_i} B_i \end{bmatrix} \}_{i=0}^{\pdim},C,D)$ 
and  $\mathscr{S}_{\widetilde{\mathcal{S}}}=(\{\sqrt{p_i} \tilde{A}_i,\begin{bmatrix} \tilde{\Bs}_i & \sqrt{p_i} \tilde{B}_i \end{bmatrix} \}_{i=0}^{\pdim},\tilde{C},\tilde{D})$
be the corresponding dLSSs. By Lemma \ref{min:col:lem1}, $\mathscr{S}_{\mathcal{S}}$ and $\mathscr{S}_{\widetilde{\mathcal{S}}}$ are realizations of the same sub-Markov function,  
they are both span-reachable and observable. Hence, $\tilde{D}=D$, and by \cite{PetreczkyLPVSS} there exists a non-singular matrix $T$ such that
\eqref{lemma:min1:eq1} holds, except for $\tilde{K}_iQ_i=TK_i$, $i=1,\ldots,\pdim$.
%$\tilde{A}_i=TA_iT^{-1}$, $\tilde{B}_i=TB_i$ and $\tilde{\Bs}_i=T\Bs_i$, $i=1,\ldots,\pdim$ and $\tilde{A}_i=TA_iT^{-1}$.
%In order to show that $\mathcal{S}$ and $\widetilde{\mathcal{S}}$ are isomorphic, it is enough to show that $TK_i=\tilde{K}_i$ for $i=1,\ldots,\pdim$. 
In order to show the latter, notice that $\tilde{\Bs}_i=T\Bs_i$, $i=1,\ldots,\pdim$ holds, and use the
the proof of \cite[Part (ii), Theorem 2]{PetreczkyBilinear}, by applying it to the sLSS $\mathcal{S}_{as}=(\{A_i,K_i\}_{i=1}^{\pdim},C,\eb^s)$ and $\widetilde{\mathcal{S}}=(\{\tilde{A}_i,\tilde{K}_i,\}_{i=1}^{\pdim},\tilde{C},I,\eb^s)$. More precisely, %using the terminology of  \cite[Part (ii), Theorem 2]{PetreczkyBilinear},  the representations associated with $\mathcal{S}_s$ and
%$\widetilde{\mathcal{S}_{s}}$ and $R_{\mathcal{S}}=(\nx,\pdim \ny, \{\sqrt{p_i} A_i\}_{i=1}^{\pdim},\Bs,C)$ and $R_{\widetilde{\mathcal{S}}}=(\nx,\pdim \ny, \{\sqrt{p_i} \tilde{A}_i\}_{i=1}^{\pdim},\tilde{\Bs},\tilde{C})$,
%where $\Bs=\begin{bmatrix} \Bs_1 & \ldots & \Bs_{\pdim} \end{bmatrix}$ and  $\hat{\Bs}=\begin{bmatrix} \hat{\Bs}_1 & \ldots & \hat{\Bs}_{\pdim} \end{bmatrix}$ and $\nx$ is the dimension of $\mathcal{S}$.
%It is then easy to see that $T$, $T\Bs=\hat{\Bs}$, $\tilde{A}_i=TA_iT^{-1}$, $i=1,\ldots,\pdim$, $\tilde{C}=CT^{-1}$. That is, 
$T$ satisfies the properties of the matrix $T$ defined in the  proof of \cite[Part (ii), Theorem 2]{PetreczkyBilinear}. Then, the rest of the proof of \cite[Part (ii), Theorem 2]{PetreczkyBilinear} can be repeated verbatim.
% leading to $TK_iQ_i=\tilde{K}_iQ_i$, $i=1,\ldots,\pdim$.
\end{pf}

\subsection{Proof of Lemma \ref{final:alg:lemma}}
\begin{pf} By definition \cite{LCSSArxive}, we know that $\yb(t) = \yb^s(t) + \yb^d(t)$. From that, we have $\zwy(t) = \zwys(t) + \zwyd(t)$. 
From \cite[Lemma 11]{PetreczkyBilinear}, we know that the entries of $\yb^d(t)$ and the entries of $\zwyd(t)$ are in $\Htplusu$. 
From Lemma \ref{decomp:lemma:inv:pf2.2}, it follows that $\zwys(t)$ is orthogonal to $\Htplusu$ and hence to $\zwyd(t)$.
That is, $\Tsigy = \expect{\yb(t)\zwy(t)} = \expect{\yb^s(t)\zwys(t)} + \expect{\yb^s(t)\zwyd(t)} + \expect{\yb^d(t)\zwys(t)} + \expect{\yb^d(t)\zwyd(t)}$. And since $\yb^s(t)$ is orthogonal to $\Htplusu$ and the entries of $\zwyd(t)$ belong to $\Htplusu$, then $\expect{\yb^s(t)\zwyd(t)}=0$.
Similary, since the components of $\yb^d(t)$ are in $\Htplusu$ and $\zwys(t)$ is orthogonal to $\Htplusu$ by Lemma \ref{decomp:lemma:inv:pf2.2}, then $\expect{\yb^d(t)\zwys(t)}=0$.
Then $\Lwy = \expect{\yb(t)\zwy(t)} = \expect{\yb^s(t)\zwys(t)} + \expect{\yb^d(t)\zwyd(t)} = \Lwys + \Lwyd$.

Also, note that $\Lambda_{w}^{\yb^d,\ub} = \expect{\yb^d(t)\zwu(t)} = \expect{(\yb(t) - \yb^s(t))\zwu(t)} = \expect{\yb(t)\zwu(t)} - \expect{\yb^s(t)\zwu(t)} = \expect{\yb(t)\zwu(t)}$ since $\yb^s(t)$ is orthogonal to $\Htplusu$ and hence to $\zwu(t)$. Then $\Lambda_{w}^{\yb^d,\ub} = \Lambda_{w}^{\yb,\ub}$.

Finally, note that $\Tsigy = \expect{\zsigy(t)(\zsigy(t))^T} = \expect{\zsigys(t)(\zsigys(t))^T} + \expect{\zsigys(t)(\zsigyd(t))^T} + \expect{\zsigyd(t)(\zsigys(t))^T} + \expect{\zsigyd(t)(\zsigyd(t))^T}$. 
By \cite[Lemma 11]{PetreczkyBilinear}, the components of $\zsigyd(t)$ belong to $\Htplusu$ and $\zsigys(t)$ is orthogonal to $\Htplusu$ hence $\expect{\zsigys(t)(\zsigyd(t))^T} = \expect{\zsigyd(t)(\zsigys(t))^T} = 0$. Subsequently, $\Tsigy = \Tsigyd + \Tsigys$. 
Then equation \eqref{covseq:lemma:eq1:realization1} of Lemma \ref{final:alg:lemma} holds. 
Equation \eqref{covseq:lemma:eq1:realization} of Lemma \ref{final:alg:lemma} follows from \cite[Lemma 4]{PetreczkyBilinear}.
\end{pf}

\subsection{Proof of Corollary \ref{thm:cra:col1}}
\begin{pf}
Let us recall the equation \eqref{min:markov:eq1}.
\begin{equation*}
       M_{\yb,\ub}(w)=\begin{cases} \begin{bmatrix} \Lambda^{\y^d,\bu}_{w} Q_{\bu}^{-1},  & \Lambda_w^{\y^s} \end{bmatrix} & w \ne \epsilon 
   \vspace{0.2cm} \\
   %& \\
   \begin{bmatrix} \Lambda^{\y^d,\bu}_{\epsilon} Q_{\bu}^{-1},  & I_{\ny} \end{bmatrix}
   & w=\epsilon 
   \end{cases}
\end{equation*}
From Lemma \ref{basis_red:lemma}, it follows that $\mathscr{S}_d=(\{\tilde{A}_i,\tilde{B}_i\}_{i=1}^{\pdim},\tilde{C},\tilde{D})$ is a minimal dLSS realization of the covariance sequence $\Psi = \Lambda^{\y^d,\bu}_{\epsilon} Q_{\bu}^{-1}$.
Then from Lemma \ref{final:alg:lemma}, it follows that \eqref{covseq:lemma:eq1:realization} holds, and hence, step \ref{step3:Alg:min} and step \ref{step4:Alg:min} of Algorithm \ref{algo:Nice_select_realization} result in computing $\left\{M_{\yb,\ub}(w)\right\}_{w \in \mathcal{L}_{\alpha,\beta}}$.
Then by Lemma \ref{basis_red:lemma}, $\mathscr{S}=(\{\hat{A}_i,\left[\hat{B}_i \ \hat{\Bs}_i\right]\}_{i=1}^{\pdim},\hat{C},\hat{D})$ is a minimal dLSS realization of $M_{\yb,\ub}$.
From \textbf{part \ref{theo:min:2}} of Theorem \ref{theo:min}, it follows that the LSS $\tilde{\mathcal{S}}$ is a minimal sLSS realization of $(\yb,\ub,\btheta)$ in innovation form. 
The last statement of the Corollary \ref{thm:cra:col1} follows from Lemma \ref{basis_red:lemma}.
\end{pf}

\subsection{Proof of Lemma \ref{lm:consistency}}
\begin{pf}
The matrices of the sLSS returned by Algorithm \ref{algo:Nice_select} are continuous functions of the entries of the Hankel matrices $\mathcal{H}^{M_{\alpha,\beta},N}$, $\mathcal{H}^{M_{\sigma, \alpha, \beta},N}$,$\mathcal{H}_{\alpha,\sigma}^{M,N}$ and $\mathcal{H}_{\beta}^{M,N}$. 
As $N$ goes to infinity, the empirical covariances $T^{\y,\mathbf{u}}_{\epsilon,w,N}$,$T^{\y,\mathbf{y}}_{\epsilon,w,N}$ and  $T^{\y,\y}_{\sigma,\sigma,N}$ go, respectively, to $\Lambda^{\y,\mathbf{u}}_{w}$,$\Lambda^{\y,\mathbf{y}}_{w}$ and $T^{\y,\y}_{\sigma,\sigma}$.
In other words, the empirical covariances are equal to the true covariances when $N \rightarrow +\infty$. 
In this case, the Hankel matrices computed from the empiraqual covariances are equal to the true Hankel matrices $\mathcal{H}^{M_{\alpha,\beta}}$, $\mathcal{H}^{M_{\sigma, \alpha, \beta}}$,$\mathcal{H}_{\alpha,\sigma}^{M}$ and $\mathcal{H}_{\beta}^{M}$. 
And because of the continuous relationship in Algorithm \ref{algo:Nice_select_true_deter_basic}, the matrices $\{\tilde{A}_{\sigma}^N, \tilde{B}^N_{\sigma}, \tilde{K}^{N,\mathcal{I}}_{\sigma}, \tilde{Q}^{N,\mathcal{I}}_{\sigma} \}_{\sigma=1}^{\pdim}, \tilde{C}^N, \tilde{D}^N$ returned by Algorithm \ref{algo:Nice_select} are equal to true matrices $\{\tilde{A}_{\sigma}, \tilde{B}_{\sigma}, \tilde{K}_{\sigma}, \tilde{Q}_{\sigma} \}_{\sigma=1}^{\pdim}, \tilde{C}, \tilde{D}$. 
Then, note also that the returned sLSS is minimal beacause the dLSS is returned by Algorithm \ref{algo:Nice_select_true_deter_basic} so the dLSS is minimal, hence, the associated sLSS is also minimal by Lemma \ref{min:col:lem2}.
\end{pf}

\end{document}